\documentclass[11pt]{article}
\usepackage{exscale,makeidx,amsmath,amssymb,psfrag,epsf,epsfig,comment}
\usepackage{graphicx}
 \DeclareGraphicsExtensions{.eps, .ps}
\begin{document}

\title{\vspace{-0.9cm}Restricted exchangeable partitions and embedding\\ of associated hierarchies in continuum random trees\vspace{-0.1cm}}

\author{Bo Chen\thanks{University of Oxford; email: chen@stats.ox.ac.uk; 1 South Parks Road, Oxford OX1 3TG, UK; } \and Matthias Winkel\thanks{University of Oxford; 1 South Parks Road, Oxford OX1 3TG, UK; tel: +44-1865-2-72875; fax: +44-1865-2-72595,  email: winkel@stats.ox.ac.uk; Corresponding author}}



\newtheorem{definition}{Definition}
\newtheorem{examples}[definition]{Examples}
\newtheorem{example}[definition]{Example}
\newtheorem{prop}[definition]{Proposition}
\newtheorem{lemm}[definition]{Lemma}
\newtheorem{thm}[definition]{Theorem}
\newtheorem{cor}[definition]{Corollary}
\newtheorem{conj}[definition]{Conjecture}
\newtheorem{rem}[definition]{Remark}
\newtheorem{proc}{Procedure}
\newenvironment{pf}{\begin{trivlist}\item[] \textbf{Proof.}}
                     {\hspace*{\fill} $\square$\end{trivlist}}
\newenvironment{pfofthm1}{\begin{trivlist}\item[] \textbf{Proof of Theorem \ref{thm1}.}}
                     {\hspace*{\fill} $\square$\end{trivlist}}

\newenvironment{pfofthm2}{\begin{trivlist}\item[] \textbf{Proof of Theorem \ref{thm2}.}}
                     {\hspace*{\fill} $\square$\end{trivlist}}

\newenvironment{pfofthm3}{\begin{trivlist}\item[] \textbf{Proof of Theorem \ref{thm3}.}}
                     {\hspace*{\fill} $\square$\end{trivlist}}
                     
\newenvironment{pfofprop8}{\begin{trivlist}\item[] \textbf{Proof of Proposition \ref{propagpd}.}}
                     {\hspace*{\fill} $\square$\end{trivlist}}
                     
\newenvironment{pfofprop10}{\begin{trivlist}\item[] \textbf{Proof of Proposition \ref{skewsc}.}}
                     {\hspace*{\fill} $\square$\end{trivlist}}
                     
\newenvironment{pfoffill}{\begin{trivlist}\item[] \textbf{Proof of Proposition \ref{fill}.}}
                     {\hspace*{\fill} $\square$\end{trivlist}}

\newcommand{\rd}{\mathrm{d}}
\newcommand{\fs}{\mathbf{s}}
\newcommand{\fS}{\mathbf{S}}
\newcommand{\bP}{\mathbb{P}}
\newcommand{\bE}{\mathbb{E}}
\newcommand{\bN}{\mathbb{N}}
\newcommand{\bZ}{\mathbb{Z}}
\newcommand{\bR}{\mathbb{R}}
\newcommand{\bH}{\mathbb{H}}
\newcommand{\cP}{\mathcal{P}}
\newcommand{\cS}{\mathcal{S}}
\newcommand{\cE}{\mathcal{E}}
\newcommand{\cC}{\mathcal{C}}
\newcommand{\cF}{\mathcal{F}}
\newcommand{\cG}{\mathcal{G}}
\newcommand{\cK}{\mathcal{K}}
\newcommand{\cL}{\mathcal{L}}
\newcommand{\cB}{\mathcal{B}}
\newcommand{\cO}{\mathcal{O}}
\newcommand{\cU}{\mathcal{U}}
\newcommand{\cT}{\mathcal{T}}
\newcommand{\bT}{\mathbb{T}}
\newcommand{\cH}{\mathcal{H}}
\newcommand{\cR}{\mathcal{R}}
\newcommand{\ft}{\mathbf{t}}

\maketitle

\vspace{-1cm}

\begin{abstract} We introduce the notion of a restricted exchangeable partition of
  $\bN$. We obtain 
  integral representations, consider associated fragmentations, 
  embeddings into continuum random trees and convergence to such limit trees. 
  In particular, we deduce from the general theory developed here a limit result conjectured previously for 
  Ford's alpha model and its extension, the alpha-gamma model, where 
  restricted exchangeability arises naturally. 
  
  \emph{AMS 2000 subject classifications: 60G09, 60J80.\newline
Keywords: Exchangeability, hierarchy, coalescent, fragmentation, continuum random tree, renewal theory}\vspace{-0.1cm}

\end{abstract}

\section{Introduction}

This paper introduces the concept of restricted exchangeability, which captures a weak form
of exchangeability that occurs naturally in models such as the alpha-gamma tree model of \cite{CFW}.  

\subsection{Motivating example: alpha-gamma trees as random hierarchies}\label{intro1}

An important motivation for this paper is the study of the limiting behaviour of the alpha-gamma tree-growth model \cite{CFW}, which is based on a simple stochastic growth rule to build a tree $T_{n+1}$ from a tree $T_n$ by adding a leaf (degree-1 vertex) labelled $n+1$. Let us specify this rule in a framework of 
hierarchies (also called total partitions or fragmentations in the literature). 

Following \cite{Sch-1870,Sta-99,MPW,HaP-09}, we call 
\em hierarchy on $B\subseteq\bN$ \em any subset $\ft_B$ of the power set of $B$ such that $B\in\ft_B$ and $\{j\}\in\ft_B$ for all $j\in B$, and so that for every $A,A^\prime\in\ft_B$, either $A\subseteq A^\prime$ or $A^\prime\subseteq A$ or $A\cap A^\prime=\varnothing$. To avoid trivialities, we also require $\varnothing\in\ft_B$. We say that a strict subset $A\in\ft_B$ of $A^\prime\in\ft_B$ is \em a maximal subset of $A^\prime$ in $\ft_B$ \em if for all $A^{\prime\prime}\in\ft_B$ with $A\subseteq A^{\prime\prime}\subseteq A^\prime$ either $A=A^{\prime\prime}$ or $A^{\prime\prime}=A^\prime$. For finite $B\subset\bN$ with $\#B\ge 2$,
the maximal subsets $A_1,\ldots,A_k$ of $B$ in $\ft_B$ form a partition of $B$ and the \em restrictions \em $\ft_{A_i}=\ft_B\cap A_i=\{A\cap A_i\colon A\in\ft_B\}$ are hierarchies on $A_i$, $i\in[k]:=\{1,\ldots,k\}$; a hierarchy $\ft_B$ fully encodes a \em rooted tree\em, i.e. a connected acyclic graph, with vertex set $\ft_B$ and edge 
relation linking each set to its maximal non-empty subsets, with \em root \em $\varnothing$ related to $B$; hierarchies $\ft_{A_i}$ are the \em subtrees \em of $\ft_B$ above the \em first branchpoint \em $B$ of $\ft_B$. We call $A\in\ft_B$ \em branchpoint \em or \em internal vertex \em if $\#A\ge 2$. Denote by $\bT_n$ the set of all hierarchies on $[n]$, $n\ge 1$. We say that $\ft_n\in\bT_n$ and $\ft_{n+1}\in\bT_{n+1}$ are \em consistent \em if $\ft_n=\ft_{n+1}\cap[n]$. 

\begin{figure}[h]
\begin{center}
\vspace{-0.2cm}
\hspace{-0.0cm}\includegraphics[scale=0.3]{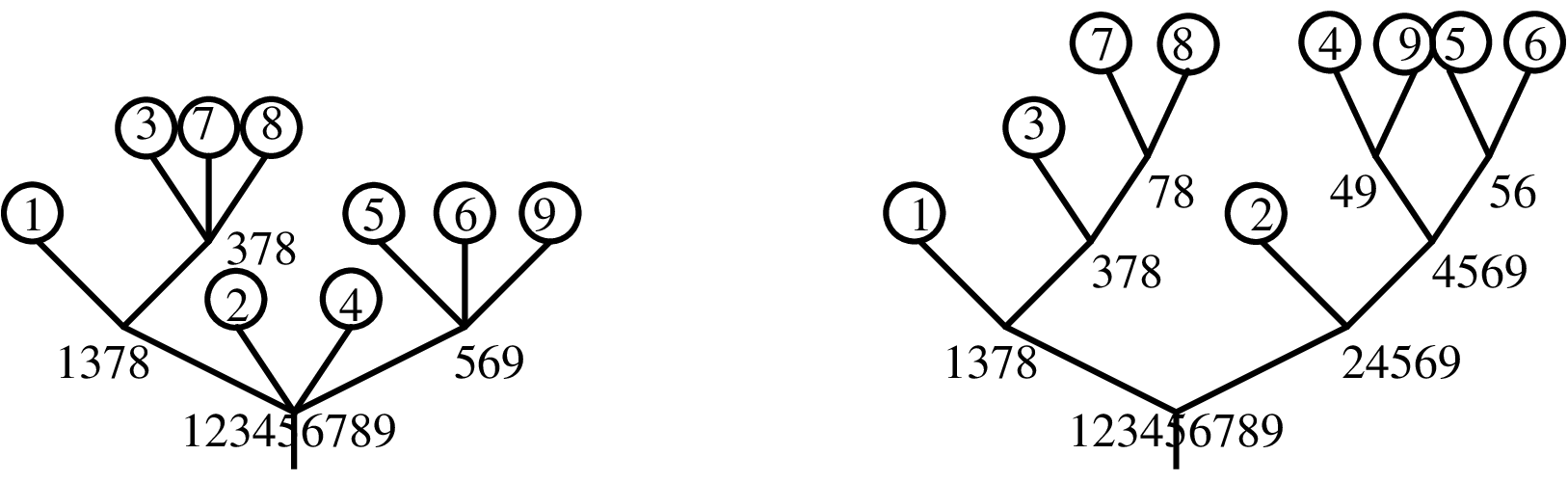}
\end{center}
\vspace{-0.7cm}
\caption{\small Two hierarchies on $B=\{1,2,3,4,5,6,7,8,9\}$ illustrated as rooted trees.}
\label{figfrag}
\end{figure}

\noindent The alpha-gamma model \cite{CFW} is a consistent family $(T_n,n\ge 1)$ of random hierarchies on $[n]$, for
which the conditional distributions of $T_{n+1}$ given $T_n$ are particularly simple. In 
terms of trees, passing from $T_n$ to $T_{n+1}$ means identifying the random place in $T_n$ where
$\{n+1\}$ connects to $T_n$:
for parameters  $0\le\gamma\le\alpha\le 1$ and for $n\ge 1$, vertex $\{n+1\}$ connects to 
\vspace{0.2cm}

\noindent\begin{minipage}[t]{0.13cm}$\;$\end{minipage}\begin{minipage}[t]{0.32cm}$\bullet$\end{minipage}\begin{minipage}[t]{15.55cm}a new vertex $\{j,n+1\}$ inserted (in the edge) below $\{j\}\in T_n$ with probability $(1-\alpha)/(n-\alpha)$;\end{minipage}\\[0.1cm]
\noindent\begin{minipage}[t]{0.13cm}$\;$\end{minipage}\begin{minipage}[t]{0.32cm}$\bullet$\end{minipage}\begin{minipage}[t]{15.55cm}a new vertex $B\cup\{n+1\}$ inserted below branchpoint $B\in T_n$ with probability $\gamma/(n-\alpha)$;\end{minipage}\\[0.1cm]
\noindent\begin{minipage}[t]{0.13cm}$\;$\end{minipage}\begin{minipage}[t]{0.32cm}$\bullet$\end{minipage}\begin{minipage}[t]{15.55cm}an existing branchpoint $B\in T_n$ with probability $((k-1)\alpha-\gamma)/(n-\alpha)$, where $k+1$ is the degree of vertex $B$ in the tree $T_n$, or equivalently $k$ is the number of blocks of the partition into maximal subsets $A_1,\ldots,A_k$ of $B$ in the hierarchy $T_n$;\end{minipage}\\[0.2cm] 
now $T_{n+1}$ is built from $T_n$ by adding $n+1$ to all vertices on the path
between $\{n+1\}$ and $\varnothing$.

\begin{figure}[h]
\begin{center}
\hspace{-0.3cm}\includegraphics[scale=0.4]{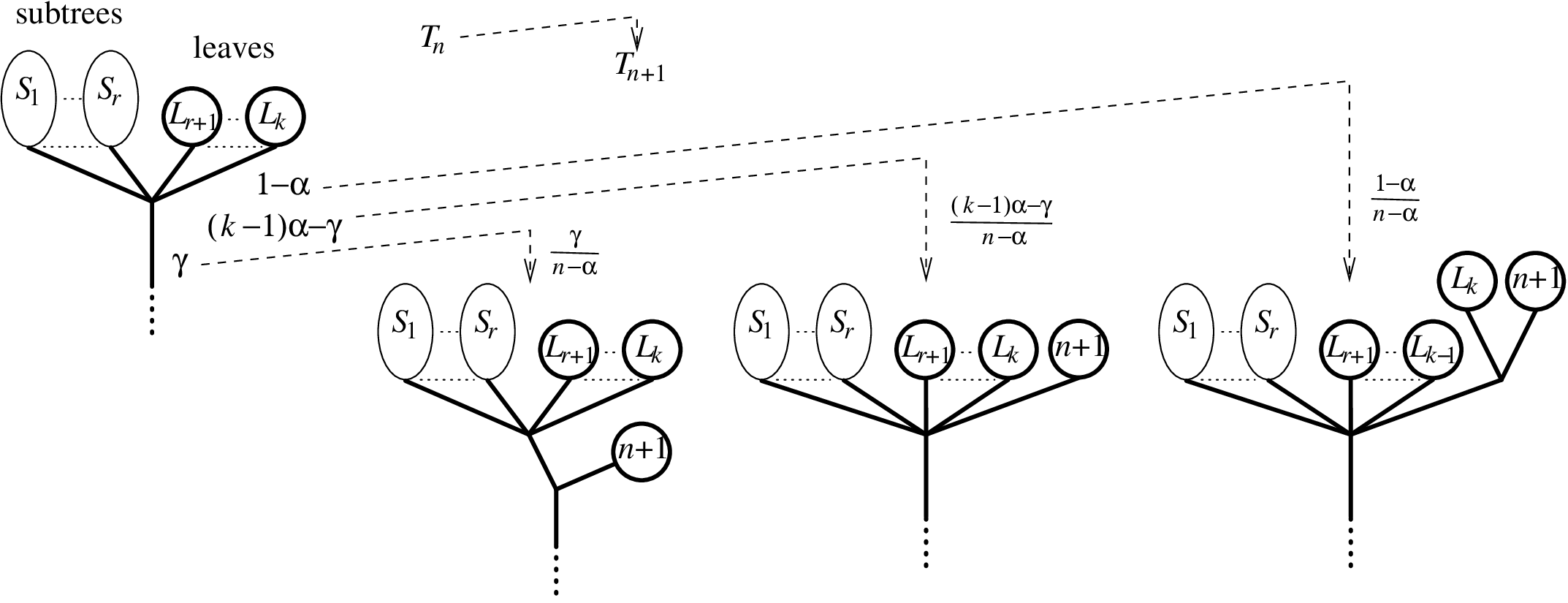}
\vspace{-0.3cm}

\begin{minipage}{15cm}\caption{\small Alpha-gamma growth rule: displayed is one internal vertex, $B$ say, of $T_n$ with 
degree $k+1$, hence vertex weight $(k-1)\alpha-\gamma$, with $k-r$ leaves 
$L_{r+1},\ldots,L_k\in[n]$ and $r$ bigger subtrees $S_1,\ldots,S_r$; all edges also carry weights, weight $1-\alpha$ and $\gamma$ are 
displayed here for the \em leaf edge \em below $\{L_k\}$ and the \em inner edge \em below $B$ only; the three associated 
possibilities for $T_{n+1}$ are displayed.}
\end{minipage}
\end{center}
\label{fig1}
\end{figure}
\vspace{-0.6cm}

\noindent A random hierarchy $T_B$ on $B$ is called \em exchangeable \em \cite{HaP-09} if for every bijection $\beta\colon B\rightarrow B$, the 
hierarchy $\beta(T_B)=\{\{\beta(j)\colon j\in A\},A\in T_B\}$ obtained by permuting labels by $\beta$ is distributed like $T_B$. An alpha-gamma tree $T_n$ for $n\ge 3$ is exchangeable iff $\gamma=1-\alpha$; note for instance that $$\bP\left(T_3=\parbox{0.9cm}{\includegraphics[height=1cm]{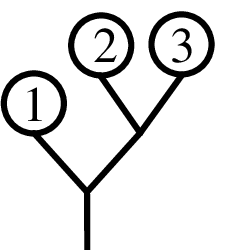}}\right)\,{=}\,\bP\left(T_3=\!\!\parbox{0.9cm}{\includegraphics[height=1cm]{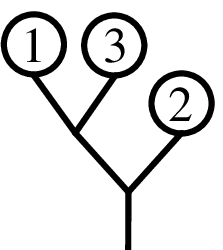}}\right)={\frac{1-\alpha}{2-\alpha}\quad\mbox{while}\quad\frac{\gamma}{2-\alpha}}=\bP\left(T_3=\!\!\parbox{0.9cm}{\includegraphics[height=1cm]{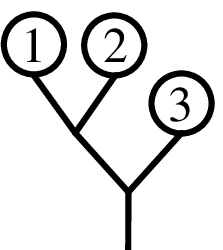}}\right).$$
However, for $\gamma\neq1-\alpha$ there is still some exchangeability. To capture this, we introduce the partition $\Pi_n$ of $T_n$ into maximal strict subsets of $[n]$ and refer to its distribution $P_n$ on the set $\cP_n$ of partitions of $[n]$ as a \em splitting rule\em. We say that $(T_n,n\ge 1)$ is a \em labelled Markov branching model \em if conditionally given $\Pi_n=\{A_1,\ldots,A_k\}$, the hierarchies $T_n\cap A_i$, $i\in[k]$, are independent and distributed as $\beta_i(T_{\#A_i})$, where $\beta_i$ is the unique increasing bijection from $[\#A_i]$ to $A_i$. Then $(P_n,n\ge 2)$ determines the distributions of $T_n$, $n\ge 1$. We will show in Section \ref{alphagamma} that the alpha-gamma model is a labelled Markov branching model with splitting rules $P_n^{\alpha,\gamma}$ satisfying
$$P_n^{\alpha,\gamma}(\pi)=P_n^{\alpha,\gamma}(\beta(\pi))\qquad\mbox{for all bijections $\beta\colon[n]\rightarrow[n]$ with $\pi\cap\{1,2\}=\beta(\pi)\cap\{1,2\}$,}$$
where $\beta(\pi)=\{\{\beta(j)\colon j\in A\},A\in\pi\}$. 
Equivalently, $P_n^{\alpha,\gamma}$ satisfies 
$P_n^{\alpha,\gamma}(\pi)=P_n^{\alpha,\gamma}(\pi^\prime)$ if $\pi\cap\{1,2\}=\pi^\prime\cap\{1,2\}$ and if $\pi=\{A_1,\ldots,A_k\}$ and $\pi^\prime=\{A_1^\prime,\ldots,A_{k}^\prime\}$ have the same multiset of block sizes $\#A_i$, $i\in[k]$, and $\#A_j^\prime$, $j\in[k]$. 

Alpha-gamma trees $(T_n,n\ge 1)$ give rise to a random hierarchy   $\cH=\{A\subset\bN\colon A\cap[n]\in T_n\linebreak\mbox{for all }n\ge 1\}$ on $\bN$. We
studied the limiting behaviour of $T_n$ and identified a scaling
limit in \cite{CFW}, but only obtained convergence in distribution. The crucial tool to strengthen
to convergence in probability is restricted exchangeability, which
we will use to embed $\cH$ and more general 
hierarchies of (restricted exchangeable) Markov branching models into suitable limit trees.

\subsection{Restricted exchangeable partitions and integral representations}

For a partition $\pi=\{\pi_i,i\in\bN\}$ of $B\subseteq\bN$ with disjoint $\pi_i$, $i\in\bN$, each non-empty $\pi_i\subseteq B$ is called a \em block \em of $\pi$. When $\pi$ has only finitely many blocks, we often omit $\varnothing$ from $\pi$. To be definite, we arrange the blocks of $\pi$ in the order of least element, i.e. $\min\pi_i<\min\pi_j$ for every $i<j$, followed by $\varnothing$ with the convention $\min\varnothing=\infty$. For finite $\pi_i$, we consider the \em block size \em $\#\pi_i$. 
We denote the set of all partitions of $B$ by $\cP_B$. Recall $[n]=\{1,\ldots,n\}$ for $n\in\bN$. Note that for $\Gamma\in\cP=\cP_\bN$, the restrictions $\Gamma|_n=\Gamma\cap[n]=\{\Gamma_i\cap[n],i\in\bN\}$ are partitions of $[n]$, $n\in\bN$. On $\cP$, consider the metric $d(\Gamma,\Gamma^\prime)=2^{-\inf\{n\ge 1\colon \Gamma|_n\neq\Gamma^\prime|_n\}}$ and the associated Borel $\sigma$-algebra.

Following de Finetti and Kingman, we call a Borel measure on the space $\cP_B$ of partitions
of $B\subseteq\bN$ \em exchangeable\em, if it is invariant under the natural action on $\cP_B$ of the symmetric group on $B$; and a random partition is called exchangeable if its distribution is exchangeable. 
Then a measure $\mu$ on $\cP$ is exchangeable if and only if the discrete measures $\mu_n$ on $\cP_n=\cP_{[n]}$, given by
\begin{equation}\mu_n(\{\pi\})=\mu(\cP^\pi),\quad\mbox{$\pi\in\cP_n$, where $\cP^\pi=\{\Gamma\in\cP\colon \Gamma|_n=\pi\}$,}
\label{restmu}
\end{equation}
are exchangeable for all $n\ge 1$. Furthermore, a measure $\mu_n$ on $\cP_n$ is exchangeable if $\mu_n(\{\pi\})=\mu_n(\{\pi^\prime\})$ for all $\pi,\pi^\prime\in\cP_n$ with the same multiset of block sizes. 
 
Several weaker forms of exchangeability have been studied in the literature, 
notably Pitman's
partial exchangeability \cite{Pit-95} and Gnedin's constrained exchangeability
\cite{Gne-06}. We introduce here a new weak form of exchangeability and discuss in Section \ref{partialconstrained} how these notions interact. 
  

\begin{definition}\label{df1}\rm For $\pi\in\cP_n$, we call a measure $\mu$ on $\cP^\pi$ \em exchangeable on $\cP^\pi$ \em if 
  $\mu(\cP^{\pi^\prime})=\mu(\cP^{\pi^{\prime\prime}})$ for all $\pi^\prime,\pi^{\prime\prime}\in\bigcup_{m\ge 1}\cP_{n+m}$ with the 
  same multiset of block sizes and with $\pi^\prime\cap[n]=\pi^{\prime\prime}\cap[n]=\pi$.
  
  A measure $\mu$ on $\cP$ is called \em restricted exchangeable \em (RE) if there is $\cC\subset\cK:=\bigcup_{n\ge 1}\cP_n$ s.th.\vspace{-0.1cm}
  \begin{itemize}\item no $\pi\in\cC$ is the restriction of another $\pi^\prime\in\cC$,
      \vspace{-0.2cm}
    \item the measure $\mu$ is carried by $\bigcup_{\pi\in\cC}\cP^\pi$, i.e. $\mu(\cP\setminus\bigcup_{\pi\in\cC}\cP^\pi)=0$,
      \vspace{-0.2cm}
    \item and for each $\pi\in\cC$, the restriction of $\mu$ to $\cP^\pi$ is finite and exchangeable on $\cP^\pi$. 
  \end{itemize}
\end{definition}

\begin{rem}\label{rm2}\rm A measure on $\cP^\pi$ is exchangeable on $\cP^\pi$ if and only if
  $\mu(\cP^{\pi^\prime})=\mu(\cP^{\beta(\pi^\prime)})$ for all $\pi^\prime\in\cP_{n+m}$ and all bijections 
  $\beta\colon[n+m]\rightarrow[n+m]$ with $\pi^\prime\cap[n]=\beta(\pi^\prime)\cap[n]=\pi$, $m\ge 1$.
  
  Note that the set of admissible bijections $\beta$ depends on $\pi^\prime$, and while $\beta(j)=j$, $j\in[n]$, makes 
  $\beta$ admissible, there are many other admissible bijections. The point is that the specific blocks containing 
  $\pi_i$ in $\pi^\prime$ and $\beta(\pi^\prime)$ may have different sizes (while the multisets of all block sizes of $\pi^\prime$ and $\beta(\pi^\prime)$ coincide). This is
  an important feature of our definition of 
  restricted exchangeability. The apparently more natural but strictly weaker concept obtained by restricting the admissible 
  bijections to the subgroup of those with $\beta(j)=j$, $j\in[n]$, is less convenient to work with, since integral 
  representations of such measures -- which we might call \em weakly \em RE -- no longer just involve measures on 
  decreasing sequences, cf.\ Theorem \ref{thm1}. 
\end{rem}


Let $S^\downarrow=\{\fs=(s_i,i\ge 1)\colon s_1\ge s_2\ge\cdots\ge 0,\sum_{i\ge 1}s_i\le 1\}$.
For $\fs\in S^\downarrow$, \em Kingman's paintbox \em \cite{Kin-78rep} is obtained from independent random variables $(\xi_r,r\ge 1)$ with
respective distributions
$$\bP(\xi_r=i)=s_i,\quad i\ge 1,\qquad \bP(\xi_r=-r)=s_0:=1-\textstyle\sum_{i\ge 1}s_i,$$
as the distribution $\kappa_\fs$ on $\cP$ of the exchangeable 
partition $\Pi=\{\{n\ge 1\colon \xi_n=i\},i\in\bZ\}$, which puts any two $m,n\in\bN$ into the same block if and only if $\xi_m=\xi_n$. By the Strong Law of Large Numbers, the vector of block sizes $\#(\Pi\cap[n])^\downarrow$ in decreasing order of size has \em asymptotic frequencies \em\vspace{-0.3cm} $$|\Pi|^\downarrow=\lim_{n\rightarrow\infty}\frac{1}{n}\#(\Pi\cap[n])^\downarrow=(s_i,i\ge 1)=\fs.$$
It is well-known \cite{Kin-78rep,Ald-85,ker89} that exchangeable measures on $\cP$ admit integral representations $\mu=\int_{S^\downarrow}\kappa_\fs\nu(d\fs)$. To establish integral representations for RE measures here, we introduce \em modified paintboxes \em $\kappa_\fs^\pi$, $\pi\in\cK=\bigcup_{n\ge 1}\cP_n$, by conditioning $\kappa_\fs$ on the cylinder set $\cP^\pi=\{\Gamma\in\cP\colon \Gamma|_n=\pi\}$ of $\pi$ in $\cP$, but note that this conditioning is degenerate in some cases; see Section 2 for details.\pagebreak
   

\begin{thm}[Integral representation]\label{thm1} Let $\mu$ be a measure on $\cP$. Then $\mu$ is RE
  if and only if there are a subset $\cC\subset\cK$ such that no $\pi\in\cC$ is the restriction of another $\pi^\prime\in\cC$, and for 
  each $\pi\in\cC$ a finite measure $\nu_\pi$ on $S^\downarrow$ such that\vspace{-0.2cm}
$$\mu=\sum_{\pi\in\cC}\int_{S^\downarrow}\kappa_\fs^\pi\nu_{\pi}(d\fs).$$
\end{thm}
Note that a RE measure $\mu$ can be infinite, if $\cC$ is infinite. However, as $\cC\subset\cK$ is countable, such infinite measures will still be $\sigma$-finite, because 
they are finite on $\cP^\pi$, $\pi\in\cC$. 


\begin{examples}\rm\begin{enumerate}

\item[(i)] For $B\subseteq\bN$, let ${\bf 1}_B$ be the trivial partition of a single block $B$. \em Dislocation measures \em are
measures on $\cP$ carried by $\cP\setminus\{{\bf 1}_\bN\}$, finite on $\cP^\pi$, $\pi\in\cK\setminus\{{\bf 1}_{[n]},n\ge 1\}$. We set
$\cC=\{\{[j],\{j+1\}\},j\ge 1\}$ to naturally decompose
$\cP\setminus\{{\bf 1}_\bN\}=\bigcup_{\pi\in\cC}\cP^\pi$. Bertoin's \cite{Ber-hom} possibly infinite \em exchangeable \em dislocation measures, in the sense of (\ref{restmu}), are exchangeable and finite on $\cP^\pi$, $\pi\in\cC$, so they satisfy Definition \ref{df1}. 
See Sections \ref{sec13} and \ref{fragcoal}.  
\item[(ii)] We can associate dislocation measures with Ford's alpha model \cite{For-05} and the alpha-gamma Markov branching model \cite{CFW}, defined in Section \ref{intro1}, so that $\mu_{\alpha,\gamma}(\cP^\pi)=\lambda_n^{\alpha,\gamma}P_n^{\alpha,\gamma}(\pi)$, $\pi\in\cP_n\setminus\{{\bf 1}_{[n]}\}$, for consistent rates 
$\lambda^{\alpha,\gamma}_n$, $n\ge 2$. These dislocation measures $\mu_{\alpha,\gamma}$ are RE, but not exchangeable, as we illustrated in terms of splitting rules $P_n^{\alpha,\gamma}$ at the end of Section \ref{intro1}. See Section \ref{fragcoal} for an exploration of
the relationship between splitting rules and dislocation measures in a general RE framework.
\end{enumerate}
\end{examples}
From Theorem \ref{thm1} we deduce an integral representation for restricted 
exchangeable dislocation measures. For simplicity we only allow as decomposition of $\cP$ in Definition \ref{df1} the most relevant and natural $\cC=\{\{[j],\{j+1\}\},j\ge 1\}$. 

\begin{cor}\label{cor5} Let $\kappa$ be a RE measure with $\cC=\{\{[j],\{j+1\}\},j\ge 1\}$.
  Then for each $j\ge 1$,
  there are constants $c_j\ge 0$ and $k_j\ge 0$, and a measure $\nu_j$ on $S^\downarrow$ with 
  $$\nu_j(\{(0,0,\ldots)\})=\nu_j(\{(1,0,\ldots)\})=0\quad\! \mbox{\rm and}\!\quad\int_{S^\downarrow}\left(s_01_{\{j=1\}}+\sum_{i\ge 1}s_i^j(1-s_i)\right)\nu_j(d\fs)<\infty,$$
  such that, for $\varepsilon^{(j)}=\{\{j\},\bN\setminus\{j\}\}$ and $\omega^{[j]}=\{[j],\{j+1\},\{j+2\},\ldots\}$, $j\ge 1$,
  $$\kappa=c_1\delta_{\varepsilon^{(1)}}+\sum_{j\ge 1}\left(c_j\delta_{\varepsilon^{(j+1)}}+k_j\delta_{\omega^{[j]}}+\int_{S^\downarrow}\kappa_\fs(\cdot\cap\cP^j)\nu_j(d\fs)\right),\qquad\mbox{\rm where $\cP^j=\cP^{\{[j],\{j+1\}\}}$.}$$
\end{cor}
In the exchangeable case, we have $(c_j,k_j,\nu_j)=(c,0,\nu)$, $j\ge 1$, as was shown by Bertoin \cite{Ber-hom}.

\subsection{RE hierarchies and continuum random trees}\label{sec13}

In the context of our motivating example, the alpha-gamma model, we demonstrated how consistent Markov branching trees give rise to a random hierarchy $\cH$ of $\bN$. Let us investigate this in the context of Bertoin's systematic studies \cite{Ber-book} of exchangeable homogeneous and exchangeable self-similar $\cP$-valued fragmentation processes $(F^*(t),t\ge 0)$ in continuous time, and of Haas and Miermont's \cite{HM} associated self-similar continuum random trees (CRTs).

Bertoin described exchangeable homogeneous fragmentation processes in terms of an exchangeable dislocation measure $\kappa=\sum_{j\ge 1}c\delta_{\varepsilon^{(j)}}+\int_{\cS^\downarrow}\kappa_\fs\nu(d\fs)$ on $\cP$. Informally, blocks fragment independently; for each $\pi\in\cP_n\setminus\{{\bf 1}_{[n]}\}$, there is a competing rate $\kappa(\cP^\pi)$ at which a given block $F_i^*(t)$ undergoes a split whose effect on the first $n$ block members is a partition according to $\pi$. For an $\alpha$-self-similar fragmentation process, this rate is increased (in the case $\alpha>0$) by a factor $|F_i^*(t)|^{-\alpha}$ depending on the asymptotic frequency $|F_i^*(t)|$ of the block. The rate increase is such that singleton blocks and indeed the all-singleton state ${\bf 0}_\bN$ are obtained in finite time.\pagebreak

\noindent Under some regularity conditions, \cite{HM} constructed self-similar CRTs $(\cT_{(\alpha,\nu)},\mu)$ with characteristic pair $(\alpha,\nu)$, i.e.\ random path-connected compact metric spaces $(\cT_{(\alpha,\nu)},\rd)$ equipped with a root $\rho\in\cT_{(\alpha,\nu)}$ and a probability measure $\mu$ on $\cT_{(\alpha,\nu)}$, and with the tree property that there are no cyclic paths. Self-similarity here means that conditionally given the tree up to height $t$ above the root and given subtree masses $\mu(S_i(t))=m_i(t)$ above height $t$, the subtrees $S_i(t)$, $i\ge 1$,  above height $t$, are like independent copies of $\cT_{(\alpha,\nu)}$, with masses rescaled by $m_i(t)$ and distances rescaled by $(m_i(t))^{\alpha}$. These CRTs can be considered as genealogical trees of Bertoin's fragmentation processes; for a $\mu$-distributed i.i.d.\ sample $\Sigma_n^*$, $n\ge 1$, in $\cT_{(\alpha,\nu)}$, we obtain an $\alpha$-self-similar fragmentation process by considering the partition-valued process that has $\{n\ge 1\colon\Sigma_n^*\in S_i(t)\}$, $i\ge 1$, as non-singleton blocks and all other integers in singleton blocks at time $t$, $t\ge 0$.

To any exchangeable $\cP$-valued fragmentation process we associate the exchangeable hierarchy $\cH^*=\{F_i^*(t),i\ge 1,t\ge 0\}$ of all blocks ever visited, equivalently $\cH^*=\{\cL^*(\cT^v)\colon v\in\cT_{(\alpha,\nu)}\}$, where $\cL^*(\cT^v)=\{n\in\bN\colon\Sigma_n^*\in\cT^v\}$ and $\cT^v$ is the subtree of $\cT_{(\alpha,\nu)}$ above $v\in\cT_{(\alpha,\nu)}$. We say that the hierarchy $\cH^*$ is \em embedded in the CRT \em $\cT_{(\alpha,\nu)}$ by the sample $\Sigma_n^*\in\cT_{(\alpha,\nu)}$, $n\ge 1$. 

We now associate with any RE dislocation measure $\kappa$ a \em RE fragmentation process \em $F$, in which each block fragments independently, 
with rates $\kappa(\cP^\pi)$, $\pi\in\cP_n$, affecting the $n$ smallest block members by
partitioning according to $\pi$. We call $\cH=\{F_i(t)\colon i\ge 1,t\ge 0\}$ the associated \em RE hierarchy\em. Alternatively (see Section \ref{fragcoal}), \em RE splitting rules
\em $P_n(\pi)=\kappa(\cP^\pi)/\kappa(\cP\setminus\cP^{{\bf 1}_{[n]}})$, $\pi\in\cP_n\setminus\{{\bf 1}_{[n]}\}$, give rise to consistent \em RE labelled Markov branching trees \em $(T_n,n\ge 1)$ with splitting rules $(P_n,n\ge 2)$ that induce
a RE hierarchy $\{A\subset\bN\colon A\cap[n]\in T_n \mbox{ for all $n\ge 1$}\}$.
Embedding a non-exchangeable hierarchy $\cH$ into a CRT $\cT$ means finding $\Sigma_n\in\cT$, $n\ge 1$, with a non-trivial dependence structure, such that $\cH$ is embedded in $\cT$ by $\Sigma_n$, $n\ge 1$. 

\begin{thm}\label{thm2} Let $\alpha>0$, and let $\kappa$ be a RE dislocation measure of the form
  \begin{equation}\label{thm6form}\kappa=\sum_{j\ge 1}\int_{S^\downarrow}\kappa_\fs(\cdot\cap\cP^j)\nu_j(d\fs),\qquad\mbox{with}\quad\nu(d\fs):=\sum_{j\ge 1}\left(\sum_{i\ge 1}s_i^j(1-s_i)\right)\nu_j(d\fs)\end{equation}
  satisfying $\int_{S^\downarrow}(1-s_1)\nu(d\fs)<\infty$ and $\nu(s_0>0)=0$. Then we can
  construct $(\cT_{(\alpha,\nu)},(\Sigma_i,i\ge 1))$ such that
  $\cH=\{\cL(\cT^v_{(\alpha,\nu)})\colon v\in\cT_{(\alpha,\nu)}\}$ is a RE hierarchy with
  dislocation measure $\kappa$, embedded in a
  self-similar CRT $\cT_{(\alpha,\nu)}$ with characteristic pair $(\alpha,\nu)$, 
  where $\cL(\cT^v_{(\alpha,\nu)})=\{i\in\bN\colon\Sigma_i\in\cT^v_{(\alpha,\nu)}\}$.
\end{thm}
Our proof of Theorem \ref{thm2} in Section \ref{sec4} gives an explicit sampling procedure for
leaves $\Sigma_i\in\cT_{(\alpha,\nu)}$, $i\ge 1$, based on the self-similarity of $\cT_{(\alpha,\nu)}$ and recursive spinal decompositions of subtrees.

Theorem \ref{thm2} partly generalises \cite[Theorem 4]{PW2}. However, apart from the alpha model (the alpha-gamma model with $\gamma=\alpha$, which produces only binary trees), that theorem treats models that are not RE in the sense of Corollary \ref{cor5} nor for other decompositions of $\cP$.

It requires no extra work to also construct hierarchies associated with RE dislocation measures $\kappa$ based on different decompositions of $\cP$. However, in those more general cases, a RE measure still qualifies as a dislocation measure if and only if it is finite on $\cP^{\{[j],\{j+1\}\}}$, $j\ge 1$, and this is necessary for hierarchies to be well-defined. Hence, the decomposition in Corollary \ref{cor5} is the most natural decomposition in the context of fragmentation processes. 

Exchangeable hierarchies $\cH^*$ derived from fragmentation processes (or from Markov branching trees) have been used to construct
CRTs as scaling limits \cite{HMPW}. We carry out a similar programme here for RE hierarchies $\cH$, starting from
a RE dislocation measure of the form identified in Corollary \ref{cor5}. 
We can delabel trees $T_n=\cH\cap[n]$, but retain the root, to obtain rooted combinatorial trees $T_n^\circ$, i.e. connected acyclic graphs with no degree-2 vertex, but some degree-1 vertices, only one of which is distinguished, as the root. We can
regard $T_n^\circ$ as a metric space with unit distance between adjacent vertices and with adjacent vertices connected by unit length line segments. We use notation $T_n^\circ/a$ to scale the length of the line segments and to obtain a metric space with all connecting line segments of length $1/a$, where $a\in(0,\infty)$.

\noindent In the
exchangeable case, \cite{HMPW} obtain CRT convergence under a regular variation condition
\begin{equation}\label{regvar}\nu(s_1\le 1-\epsilon)=\epsilon^{-\alpha}\ell(1/\epsilon)\qquad\mbox{as $\epsilon\downarrow 0$; for some $\alpha\in(0,1)$ and slowly varying $\ell$}
\end{equation}
and a log-moment condition
\begin{equation}\label{logmom}\int_{S^\downarrow}\sum_{i\ge 2}s_i|\log(s_i)|^\varrho\nu(d\fs)<\infty\qquad\mbox{for some $\varrho>0$.}
\end{equation}
%

\begin{thm}\label{thm3} If in the setting of Theorem \ref{thm2}, the measure $\nu$ 
  satisfies {\rm(\ref{regvar})} and {\rm(\ref{logmom})}, and if $\nu_j=\nu_m$ for some $m\ge 1$ and all $j\ge m$, then 
  $$\frac{T_n^\circ}{n^\alpha\ell(n)\Gamma(1-\alpha)}\rightarrow\cT_{(\alpha,\nu)}\qquad\mbox{in probability, in the Gromov-Hausdorff sense.}$$
\end{thm}
Returning to the alpha-gamma model, we can now show that Theorem \ref{thm3} applies to give a scaling limit in probability. The identification of $\nu_j$, $j\ge 1$, in the
parameterisation of Corollary \ref{cor5} finally sheds some light on the peculiar splitting rules and $\nu$-measures in Ford's alpha model and the alpha-gamma model \cite{For-05,HMPW,CFW,PW2}. To do this, we follow \cite{Mie-03,HPW,MPW} and introduce Poisson-Dirichlet dislocation measures ${\rm PD}_{\alpha,\theta}^*(d\fs)$ as $\sigma$-finite measures on $S^\downarrow$ given by
  $$\bE[\sigma_1^{\theta};\sigma_1^{-1}\Delta\sigma_{[0,1]}\in d\fs],\qquad \theta>-2\alpha,\alpha\in(0,1),$$
on the interior of the parameter range, where $(\sigma_t,t\ge 0)$ is a stable subordinator with Laplace transform $\bE[e^{-\lambda\sigma_t}]=e^{-t\lambda^\alpha}$ and where $\Delta\sigma_{[0,1]}$ is the decreasing rearrangements of the jumps $\Delta\sigma_t=\sigma_t-\sigma_{t-}$, $t\in[0,1]$. For $\theta=-2\alpha$, the binary case, ${\rm PD}_{\alpha,-2\alpha}^*(d\fs)$ is defined as the ranked beta measure on $\{(x,1-x,0,\ldots),x\in(1/2,1)\}\subset S^\downarrow$ with density $x^{-\alpha-1}(1-x)^{-\alpha-1}1_{(1/2,1)}(x)$; the associated Markov branching model is Aldous's \cite{Ald-93} beta-splitting model, for $\alpha<1$.

As the references demonstrate, Poisson-Dirichlet dislocation measures give rise to some of the nicest and best-studied parametric families of exchangeable fragmentation processes, while alpha and alpha-gamma models have as their dislocation measure what we have previously written as linear combinations of Poisson-Dirichlet measures of different parameters \cite{CFW}. With the notion of restricted exchangeability, we can now obtain a stronger and more satisfactory connection. 

\begin{prop}\label{propagpd} The alpha-gamma model for $\alpha\in(0,1)$ and $\gamma\in[0,\alpha]$ is a RE Markov branching model with dislocation measure of the
  form identified in Corollary \ref{cor5} with $\nu_1=(1-\alpha){\rm PD}_{\alpha,-\alpha-\gamma}^*$ and $\nu_j=\gamma{\rm PD}_{\alpha,-\alpha-\gamma}^*$, $j\ge 2$.
\end{prop}
The boundary case $\alpha=1$ degenerates \cite{CFW} and leads to RE Markov branching models with\vspace{-0.1cm}
\begin{itemize}
\item for $\gamma=0$ star trees corresponding to $(\nu_1,c_1,k_1)=(0,0,1)$ and $(\nu_2,c_2,k_2)=(0,0,0)$;\vspace{-0.2cm}
\item for $\gamma=1$ comb trees corresponding to $(\nu_1,c_1,k_1)=(0,0,0)$ and $(\nu_2,c_2,k_2)=(0,1,0)$;\vspace{-0.2cm}
\item for $\gamma\in(0,1)$ bushy combs corresponding to $(\nu_1,c_1,k_1)=(0,0,1)$, $(c_2,k_2)=(0,0)$ and 
  $$\nu_2(s_2>0)=0\quad\mbox{and}\quad\nu_2(s_1\in dx)=\gamma x^{-2}(1-x)^{-1-\gamma}1_{(0,1)}(x)dx.$$
\end{itemize}

\subsection{Sampling consistency and the skewed Poisson-Dirichlet model}

Proposition \ref{propagpd} suggests to introduce a three-parameter family of restricted 
exchangeable fragmentation trees that we call the \em skewed Poisson-Dirichlet model\em, by setting 
$$\nu_1=\lambda{\rm PD}_{\alpha,\theta}^*,\qquad\nu_j=(1-\lambda){\rm PD}_{\alpha,\theta}^*,\quad j\ge 2,$$
for $\alpha\in[0,1]$, $\theta\ge-2\alpha$ and $\lambda\in[0,1]$. When $\lambda=(1-\alpha)/(1-\theta-2\alpha)$ and $\theta=-\alpha-\gamma$, this is the alpha-gamma model; when $\lambda=1/2$, this is the exchangeable Poisson-Dirichlet model studied in 
\cite{MPW,HPW}. We will use parameterisations by $(\alpha,\theta,\lambda)$ and $(\alpha,\gamma,\lambda)$, where $\gamma=-\alpha-\theta$. We can apply Theorem \ref{thm3} to obtain a convergence result in probability:

\begin{cor} Let $(T_n,n\ge 1)$ be a consistent family of skewed Poisson-Dirichlet trees for parameters 
  $0<\alpha<1$, $0<\gamma=-\alpha-\theta\le\alpha$ and $0\le\lambda<1$. Then\vspace{-0.1cm} $$\frac{T_n^\circ}{n^\gamma}\rightarrow\cT_{(\gamma,\nu)}\qquad\mbox{in probability, in the Gromov-Hausdorff sense,}$$\vspace{-0.4cm}
  
\noindent where $\cT_{(\gamma,\nu)}$ is a $\gamma$-self-similar CRT associated with measure\vspace{-0.2cm}
$$\nu(d\fs)=\frac{\gamma\Gamma(1-\alpha)}{(1-\lambda)\alpha\Gamma(1-\gamma/\alpha)}\left(\lambda+(1-2\lambda)\sum_{i\ge 1}
s_i^2\right){\rm PD}^*_{\alpha,\theta}(d\fs)$$\vspace{-0.3cm}

\noindent  for $\gamma<\alpha$, while in the binary case $\gamma=\alpha$ (i.e. $\theta=-2\alpha$), we have $\nu(s_1+s_2<1)=0$ and 
  $$\nu(s_1\in dx)=\frac{\alpha}{(1-\lambda)\Gamma(1-\alpha)}\left((1-\lambda)+(4\lambda-2)x(1-x)\right)x^{-\alpha-1}(1-x)^{-\alpha-1}dx.$$\vspace{-0.5cm}
\end{cor}
Regarding the alpha model, $\alpha\in(0,1)$, $\theta=-2\alpha$, $\lambda=1-\alpha$, this confirms in part a conjecture formulated in \cite{PW2}; specifically, the setting of the conjecture was the two-parameter $(\alpha,\theta)$-model that contains the alpha model as
a special case, and the conjecture claims almost sure convergence, while we only obtain 
convergence in probability here. 

Another interesting 
feature of the skewed Poisson-Dirichlet model relates to sampling consistency. Here
we say that a family of unlabelled random trees $(T_n^\circ,n\ge 1)$ is \em sampling
consistent \em if the tree $T_n^\circ$ with a uniformly chosen leaf removed is 
distributed as $T_{n-1}^\circ$. For consistent trees with exchangeable labels such as the exchangeable Poisson-Dirichlet model this is trivially so, but also and non-trivially for the alpha-gamma model that includes non-exchangeable trees \cite{CFW}. Geometrically, this gives sampling consistency for two
two-dimensional subsets of the three-dimensional parameter space (intersecting in the one-parameter family of stable trees \cite{Mie-03} for $\gamma=1-\alpha$), but somewhat surprisingly, sampling consistency does
not extend any further:\vspace{-0.05cm}

\begin{prop}\label{skewsc} The skewed Poisson-Dirichlet model is sampling consistent only for 
  parameters that reduce it to the exchangeable Poisson-Dirichlet model or to the alpha-gamma 
  model.
\end{prop}\vspace{-0.05cm}
This shows that while Theorem \ref{thm2} and \ref{thm3} always refer to Markov branching trees $T_n^\circ$ in the sense of \cite{HMPW}, they typically do not, however, satisfy the sampling consistency property of \cite{HMPW}, so that the theory developed in \cite{HMPW} does not even yield convergence in distribution for these trees, where we here establish convergence in probability.

\subsection{Structure of this paper}

In addition to proofs of main results already formulated, the content of this paper is as follows.
\begin{itemize}\vspace{-0.1cm}\item Section 2 proves Theorem \ref{thm1} and Corollary \ref{cor5} by combining
    approaches of Vershik and Kerov, and of Aldous, both in the exchangeable case.
  \vspace{-0.2cm}\item Section 3 includes a discussion of the relationship between restricted exchangeability, partial exchangeability and constrained exchangeability, and a discussion of RE dislocation measures, RE splitting rules, RE hierarchies and RE fragmentations. 
  \vspace{-0.2cm}\item 
    In Section 4, we develop a new technique to sample leaves in general  
    self-similar CRTs. We make explicit the embedding that we use to prove Theorem 
    \ref{thm2}, and we obtain
    decomposition results along subtrees spanned by the first $k$ sampled leaves 
    (Corollary \ref{cor24}).
  \vspace{-0.2cm}\item In Section 5 we prove Theorem \ref{thm3}. 
    Our approach is similar in spirit to \cite{HMPW}, but with added technical difficulties. 
    We analyse the RE embedding of Theorem \ref{thm2} in detail. While in \cite{HMPW} consideration of a single 
    $\Sigma^*\in\cT_{(\alpha,\nu)}$ gives relevant estimates for all $\Sigma_n^*$, $n\ge 1$, we
    here need individual estimates for each $\Sigma_n$, $n\ge 1$. Methods include 
    Gnedin's constrained paintboxes and renewal theory. 
    We also establish almost sure convergences of rescaled subtrees of 
    $T_n$ spanned by $k$ leaves, as first $n\rightarrow\infty$ in Proposition \ref{fdm} and then 
    also $k\rightarrow\infty$ in (\ref{doublelim}). 

  \vspace{-0.2cm}\item Section 6 provides proofs for Propositions \ref{propagpd} and \ref{skewsc}.

  \vspace{-0.2cm}\item An appendix contains the proof of a technical lemma.
\end{itemize}


\section{Integral representations, proof of Theorem \ref{thm1} and Corollary \ref{cor5}}

Our first aim is to understand exchangeability on subsets of the form $\cP^\pi\subseteq\cP$, for some  $\pi\in\cK$. Let us formally define modified paintboxes. 
For $\fs\in S^\downarrow$ let $m\ge 0$ such that $s_m>s_{m+1}=0$ (or $m=\infty$ if $s_i>0$ for all $i\ge 1$), suppose $\pi\in\cK$ has $k$ blocks $\pi_j\neq\varnothing$, $1\le j\le k$, of which $\ell$ with $\#\pi_j\ge 2$. For the paintbox $\kappa_\fs$ associated with $\fs$, we have $\kappa_\fs(\cP^\pi)>0$ iff either $s_0>0$ and $\ell\le m$, or $s_0=0$ and $k\le m$. In these cases, set $\kappa_\fs^\pi=\kappa_\fs(\cdot|\cP^\pi)$. Then $\kappa_\fs^\pi$ is a modified paintbox:
\begin{enumerate}\item Randomly assign ``colours'' $c(\pi)=(c(\pi_1),\ldots,c(\pi_k))$ to 
  the blocks $\pi_1,\ldots,\pi_k$ using the following rule (with $Z_\fs^\pi$ as normalisation constant)
  \begin{equation}\bP(c(\pi)=(i_1,\ldots,i_k))=\frac{1}{Z_\fs^\pi}\prod_{1\le j\le k}s_{i_j}^{\#\pi_j},\label{kingm}\end{equation}
  where $i_j$ is allowed to be equal to 0 iff $\#\pi_j=1$, and the $i_j$ with $i_j\ge 1$
  are pairwise distinct.
  \item Let $n$ be such that $\pi\in\cP_n$. Conditionally given $c(\pi)=(i_1,\ldots,i_k)$, set for $1\le r\le n$ and $r\in\pi_j$, 
	$$\xi_r=i_j\quad\mbox{if $i_j\ge 1$,}\qquad\mbox{and}\qquad\xi_r=-\min\pi_j\quad\mbox{if $i_j=0$},$$
    and for $r\ge n+1$, consider independent $\xi_r$ with
    $\bP(\xi_r=i)=s_i$, $i\ge 1$, $\bP(\xi_r=-r)=s_0$.
    Then $\kappa_\fs^\pi$ is the distribution of the partition $\Pi=\{\{n\ge 1\colon\xi_n=i\},i\in\bZ\}$, which puts any two $n,n^\prime\in\bN$ into the same block if and only if $\xi_n=\xi_{n^\prime}$.
\end{enumerate}
In the degenerate case when $\kappa_\fs(\cP^\pi)=0$, the numerator of (\ref{kingm}) always vanishes. Roughly speaking, we use all colours $1,\ldots,m$ for the largest blocks of $\pi$.
Formally, we replace 1. by 1$^\prime$.:

\begin{enumerate}\item[1$^\prime$.] Randomly assign ``colours'' using the following rule (with $Z_\fs^\pi$ as normalisation constant): 
$$\bP(c(\pi)=(i_1,\ldots,i_k))=\frac{1}{Z_\fs^\pi}\prod_{1\le j\le k\colon i_j\neq 0}s_{i_j}^{\#\pi_j},$$
if $\{i_1,\ldots,i_k\}\!=\!\{0,\ldots,m\}$, the $i_j\!\ge\!1$ are pairwise distinct and $\sum_{j=1}^k\#\pi_j1_{\{i_j\neq 0\}}$ is maximal. 
\end{enumerate}
Step 2. is applied as before to construct $\Pi$ and hence $\kappa_\fs^\pi$. Note that $\Pi_j=\pi_j$ if $i_j=0$, while $\Pi_j\supset\pi_j$ will have
limiting frequency $s_{i_j}>0$ if $i_j\ge 1$. 
  

Now $\kappa_\fs^\pi(\cP\setminus\cP^\pi)=0$ and, for $\pi^\prime=(\pi_1^\prime,\ldots,\pi_{k^\prime}^\prime)\in\cK^\pi:=\{\pi^\prime\in\cK\colon \pi^\prime\cap[n]=\pi	\}$, 
$$\kappa_\fs^\pi\left(\cP^{\pi^\prime}\right)=\frac{1}{Z^\pi_\fs}\sum_{(i_1,\ldots,i_{k^\prime})\ \textrm{admissible for }(\pi,\pi^\prime,\fs)}s_0^{\#\{J_\fs^\pi\le j\le k^\prime\colon i_j=0\}}\prod_{1\le j\le k^\prime\colon i_j\neq 0}s_{i_j}^{\#\pi_j^\prime1_{\{i_j\ge 1\}}},$$
where $J_\fs^\pi=k+1$ in the degenerate case, $J_\fs^\pi=1$ otherwise, and where $(i_1,\ldots,i_{k^\prime})$ is admissible for $(\pi,\pi^\prime,\fs)$ if $(i_1,\ldots,i_k)$ is as in 1$^\prime$.\ or 1.\ above, respectively, and if for $k+1\le j\le k^\prime$, we allow $i_j$ equal to 0 iff $\#\pi^\prime_j=1$, and the $i_j$ with $i_j\ge 1$, $1\le j\le k^\prime$, and pairwise distinct. 

For $\pi=\{\{1\}\}$, this is a well-known formula
for Kingman's paintbox $\kappa_\fs=\kappa_\fs^\pi$, with $Z^\pi_\fs=1$.
It is easy to show that, in the general case, the modified paintboxes $\kappa_\fs^\pi$ are exchangeable on $\cP^\pi$. 

\begin{prop}\label{propexchpi1} For any $n\ge 1$ and $\pi\in\cP_n$, the modified paintbox $\kappa_\fs^\pi$ can be expressed in terms of any 
  $\Gamma\in\cP^\pi$ with asymptotic frequencies $\fs$, provided that any blocks of 
  $\Gamma$ with zero asymptotic frequency are either subsets of $[n]$ or singletons, 
  as
$$\kappa_\fs^\pi(\cP^{\pi^\prime})=\lim_{r\rightarrow\infty}\frac{\#\{\pi^{\prime\prime}\in\cK^{\pi^\prime}\colon \pi^{\prime\prime}\approx\Gamma|_r\}}{\#\{\pi^{\prime\prime}\in\cK^\pi\colon \pi^{\prime\prime}\approx\Gamma|_r\}}\qquad\mbox{for all $\pi^\prime\in\cK^\pi=\{\pi^\prime\in\cK\colon \pi^\prime\cap[n]=\pi	\}$,}$$
where we write $\pi^\prime\approx\pi^{\prime\prime}$ if $\pi^\prime$ and $\pi^{\prime\prime}$ have the same multiset of block sizes. 
\end{prop}
\begin{pf} This proof is a refinement of the relevant part of the proof of \cite[Theorem 3.1]{ker89}, Kerov's proof of Kingman's paintbox representation of exchangeable partitions in $\cP$, where we need to take into account the restriction to $\cP^\pi$. We evaluate the right-hand side. Numerator and denominator are easily calculated, e.g. for $\pi^\prime=(\pi_1^\prime,\ldots,\pi^\prime_{k^\prime})\in\cP_{n^\prime}^\pi:=\cK^\pi\cap\cP_{n^\prime}$ as
$$\#\{\pi^{\prime\prime}\in\cP^{\pi^\prime}_r\colon \pi^{\prime\prime}\approx\Gamma|_r\}=\sum {r-n^\prime\choose\#\Gamma_{i_1}|_r-\#\pi^\prime_1,\ldots,\#\Gamma_{i_{k^\prime}}|_r-\pi^\prime_{k^\prime},\#\Gamma_{\rm others}|_r}\frac{1}{\prod_{j\ge 1}p_j!},$$
where $\sum$ is over indices $(i_1,\ldots,i_{k^\prime})$ such that $\#\Gamma_{i_j}|_r-\#\pi^\prime_j\ge 0$ for all $j\in[k^\prime]$, $\Gamma_{\rm others}|_r$ is the vector of all $\Gamma_i|_r$, $i\ge 1$, except $\Gamma_{i_1}|_r,\ldots,\Gamma_{i_{k^\prime}}|_r$, and $p_j$ is the number of blocks of $\Gamma|_r$ with $j$ elements, $j\!\ge\! 1$. First assume $s_0\!=\!1\!-\!\sum_{i\ge 1}s_i\!=\!0$,
then 
the limit exists and is 
$Z_\fs^{\pi,\pi^\prime}/Z_\fs^{\pi,\pi}$, where
$$Z_\fs^{\pi,\pi^\prime}=\lim_{r\rightarrow\infty}\frac{\sum{r-n\choose\#\Gamma_{i_1}|_r-\#\pi_1^\prime,\ldots,\#\Gamma_{i_k}|_r-\#\pi_{k^\prime}^\prime,\#\Gamma_{\rm others}|_r}r^d}{{r\choose\#\Gamma_1|_r,\#\Gamma_2|_r,\ldots}}=\sum_{(i_1,\ldots,i_{k^\prime})\ \textrm{admissible for }(\pi,\pi^\prime,\fs)} \prod_{j\colon i_j\neq 0}s_{i_j}^{\#\pi_j^\prime},$$
with $d$ the minimal $\sum_{j=1}^{k^\prime}\#\pi_j^\prime1_{\{i_j=0\}}$, so that $d>0$ only in the degenerate case; this power $d$ is such that terms with higher than the minimal sum vanish as $r\rightarrow\infty$, and we identify $\kappa_\fs^\pi(\cP^{\pi^\prime})$. 

If $s_0>0$, blocks of zero limiting frequency need to be treated differently, because
their union $\widetilde{\Gamma}_0$ now has a limiting frequency, and a union $\widetilde{\pi}_0^\prime$ of blocks of $\pi^\prime$ can indeed be associated with $\widetilde{\Gamma}_0$. Specifically, we calculate a first factor as
$$\lim_{r\rightarrow\infty}\frac{\sum{r-n\choose\#\widetilde{\Gamma}_0-\widetilde{\pi}^\prime_0,\#\widetilde{\Gamma}_{i_1}|_r-\#\widetilde{\pi}_1,\ldots,\#\widetilde{\Gamma}_{i_{\widetilde{k}^\prime}}|_r-\#\widetilde{\pi}^\prime_{\widetilde{k}^\prime},\#\widetilde{\Gamma}_{\rm others}|_r}}{{r\choose\#\widetilde{\Gamma}_0|_r,\#\widetilde{\Gamma}_1|_r,\#\widetilde{\Gamma}_2|_r,\ldots}}
=\sum_{(i_1,\ldots,i_{\widetilde{k}^\prime})\ \textrm{admissible for }(\widetilde{\pi},\widetilde{\pi}^\prime,\fs)} s_0^{\#\widetilde{\pi}^\prime_0}\prod_{j=1}^{\widetilde{k}^\prime}s_{i_j}^{\#\widetilde{\pi}_j^\prime},$$
but then need to also count the further partitions of the block of size $\#\widetilde{\Gamma}_0|_r$. This yields for $\widetilde{\pi}_0^\prime=\pi_{j_1}^\prime\cup\cdots\cup\pi_{j_b}^\prime$ a positive limit factor if $d=\#\widetilde{\pi}_0^\prime-b$ is minimal, which we then calculate as
$$\lim_{r\rightarrow\infty}\frac{\sum{\#\widetilde{\Gamma}_0|_r-\#\widetilde{\pi}_0^\prime\choose\#\Gamma_{i_1}|_r-\#\pi_{j_1}^\prime,\ldots,\#\Gamma_{i_b}|_r-\#\pi_{j_b}^\prime,1,\ldots,1}r^d}{
(\#\widetilde{\Gamma}_0|_r)!}=s_0^{-d};$$
the number of available indices is asymptotically equivalent to $\#\widetilde{\Gamma}_0|_r\sim s_0r$, so that the sum contains $\sim(\#\widetilde{\Gamma}_0|_r)^b$ terms, and this contributes to the asymptotics of the numerator. Finally we sum over the different choices of $\widetilde{\pi}_0^\prime$ with $\#\widetilde{\pi}_0^\prime-b=d$ to identify
$\kappa_\fs^\pi(\cP^{\pi^\prime})$.  
\end{pf}

With these representations of the modified paintboxes, we now obtain the integral
representation of general measures that are exchangeable on $\cP^\pi$ for some $\pi\in\cP_n$. 

\begin{prop}\label{propexchpi2} Let $\mu$ be a finite measure, exchangeable on $\cP^\pi$ for some $\pi\in\cP_n$. Then there is
  a finite measure $\nu$ on $S^\downarrow$ such that $\mu=\int_{S^\downarrow}\kappa^\pi_\fs\nu(d\fs)$.
\end{prop}
\begin{pf} This proof uses a combination of the martingale method due to Vershik and Kerov \cite[Theorem 2]{vk81} and the de Finetti method used by Aldous \cite{Ald-85}. W.l.o.g., $\mu$ is a
  probability measure. 
  Let $\Pi\sim\mu$ for an exchangeable probability measure on $\cP^\pi$. For $n^\prime\ge n$ and $\pi^\prime\in\cK^\pi\cap\cP_{n^\prime}$, consider the process \vspace{-0.3cm}
$$X_r=\frac{\#\{\pi^{\prime\prime}\in\cK^{\pi^\prime}\colon \pi^{\prime\prime}\approx\Pi|_r\}}{\#\{\pi^{\prime\prime}\in\cK^\pi\colon \pi^{\prime\prime}\approx\Pi|_r\}},\quad r\ge n^\prime,$$
  in the decreasing filtration $\cF_r$ generated by the block sizes of $\Pi|_u$, $u\ge r$. By exchangeability, $X_r$ depends only on
  the block sizes $B_r$ of $\Pi|_r$ and is hence $\cF_r$-measurable
  and $\bE[X_r|\cF_{r+1}]$ only depends on $X_{r+1}$. For a multiset 
  $b$ of block sizes, denote by $m(b)$ (resp. $m^\prime(b)$) the number of partitions in $\cK^\pi$ (resp. in $\cK^{\pi^\prime}$) with
  block sizes $b$. By exchangeability, each of these is equally likely. For block sizes $B_{r+1}=b_{r+1}$, we denote by 
  $m(b_r,b_{r+1})$ the number of partitions in $\cK^{\widetilde{\pi}}$ with block sizes $b_{r+1}$, where $\widetilde{\pi}$ is any
  specific partition with block sizes $b_r$. Then there are $m(b_r)m(b_r,b_{r+1})$ partitions in $\cK^\pi$ with block sizes $b_{r+1}$ 
  that restrict to block sizes $b_r$. With this
  notation, we have $X_r=m^\prime(B_r)/m(B_r)$. Then 
  $$\bE[X_r|B_{r+1}=b_{r+1}]=\sum_{b_r}\frac{m(b_r)m(b_r,b_{r+1})}{m(b_{r+1})}\frac{m^\prime(b_r)}{m(b_r)}=\frac{1}{m(b_{r+1})}\sum_{b_r}m(b_r,b_{r+1})m^\prime(b_r)=\frac{m^\prime(b_{r+1})}{m(b_{r+1})}$$
  for all admissible $b_{r+1}$ shows that $(X_r,r\ge n^\prime)$ is a bounded martingale
  and hence converges a.s.

  On the other hand, de Finetti's theorem yields that asymptotic frequencies exist $\mu$-a.s. Specifically, consider a partition $\Pi$ with distribution $\mu$ and, independently,
a sequence $U_i$, $i\ge 1$, of auxiliary independent uniform random variables. Then the random variables
$$\Xi_j=U_i\quad\mbox{if $j\in\Pi_i$,}\qquad j\ge n+1,$$
are exchangeable. By de Finetti's theorem, they are conditionally i.i.d. and the atom
sizes $S_i$ of the random limiting distribution in random (``size-biased'') order satisfy
$$S_i=\lim_{r\rightarrow\infty}\frac{\#\{j\in\{n+1,\ldots,n+r\}\colon \Xi_j=U_i\}}{r}=\lim_{r\rightarrow\infty}\frac{\#\Pi_i\cap[r]}{r}.$$
Clearly, the latter limit does not depend on the auxiliary variables $(U_i,i\ge 1)$, so asymptotic frequencies exist $\mu$-a.s. Furthermore, $\mu$-a.e. partition is 
such that blocks with zero asymptotic frequency either only involve elements of $[n]$ or are singletons. Denote by $\nu$ the distribution on $S^\downarrow$ of the asymptotic frequencies $\fS=(S_i,i\ge 1)$ rearranged into decreasing order of $\Pi$.  

  This means that $\mu$ is concentrated on those partitions for which Proposition
  \ref{propexchpi1} yields modified paintbox representations, and we see that
  $X_r\rightarrow\kappa_\fS^\pi(\cP^{\pi^\prime})$ a.s., where $\fS\sim\nu$; but $(X_r,r\ge n^\prime)$ is 
  a bounded martingale, so exchangeability on $\cP^\pi$ yields  
  $$\!\!\!\!\int_{S^\downarrow}\kappa_\fs^\pi(\cP^{\pi^\prime})\nu(d\fs)=\bE[\kappa_\fS^\pi(\cP^{\pi^\prime})]=\bE[X_{n^\prime}]=\sum_{\widetilde{\pi}\in\cP^\pi_{n^\prime}\colon \widetilde{\pi}\approx\pi^\prime}\mu(\cP^{\widetilde{\pi}})\frac{1}{\#\{\pi^{\prime\prime}\in\cP^\pi_{n^\prime}\colon \pi^{\prime\prime}\approx\widetilde{\pi}\}}=\mu(\cP^{\pi^\prime}).\vspace{-0.65cm}$$
%
\end{pf}

This proof raises the question whether we could have done without the martingale 
method or without the de Finetti argument, as can be done in the exchangeable case.
To avoid the de Finetti argument, we would have to generalise Proposition \ref{propexchpi1}
to ensure that all $\Gamma$ for which the limits in Proposition \ref{propexchpi1}
exist converge to modified paintboxes, which seems more difficult given the exceptional
non-singleton sets of zero limiting frequency. On the other hand, our de Finetti argument only
identifies the distribution of $\Pi$ restricted to $\{n+1,n+2,\ldots\}$ and gives little information about the conditional distribution of how the blocks of 
$\pi$ attach themselves to such paintboxes. We have not found a simple and direct
argument to see why the modified paintboxes describe the only way to attach $\pi$ in an exchangeable way. 

Now recall that Theorem \ref{thm1} states that RE measures on $\cP$ are
precisely those of the form $\mu=\sum_{\pi\in\cC}\int_{S^\downarrow}\kappa_\fs^\pi\nu_\pi(d\fs)$.

\begin{pfofthm1} First consider $\mu=\sum_{\pi\in\cC}\int_{S^\downarrow}\kappa_\fs^\pi\nu_\pi(d\fs)$ with $\cC$ such that no 
  $\pi\in\cC$ is a restriction of another $\pi^\prime\in\cC$. Since $\kappa_\fs^\pi$ only charges $\cP^\pi$, the measure $\mu$ only
  charges $\bigcup_{\pi\in\cC}\cP^\pi$. Furthermore, the restrictions of $\mu$ are finite and
  exchangeable on $\cP^\pi$. Hence $\mu$ is RE.

  Conversely, let $\mu$ be any RE measure on $\cP$ with $\cC$ such that the three bullet points of Definition
  \ref{df1} hold. Then the sets $\cP^\pi$, $\pi\in\cC$, are disjoint and the restrictions of $\mu$ to 
  $\cP^\pi$ are finite and exchangeable on $\cP^\pi$. By Proposition \ref{propexchpi2}, the restrictions of $\mu$ to $\cP^\pi$ can
  be represented as $\int_{S^\downarrow}\kappa_\fs^\pi\nu_\pi(d\fs)$. Since furthermore
  $\mu(\cP\setminus\bigcup_{\pi\in\cC}\cP^\pi)=0$, we have
  $$\mu=\sum_{\pi\in\cC}\mu(\cdot\cap\cP^\pi)=\sum_{\pi\in\cC}\int_{S^\downarrow}\kappa_\fs^\pi\nu_\pi(d\fs).$$
\vspace{-0.7cm}

\end{pfofthm1} 
The proof of the Corollary \ref{cor5} is now straightforward. Note, however, that $\nu_j$ is not $\nu_\pi$ for $\pi=\{[j],\{j+1\}\}$, $j\ge 1$. Instead, we set $k_j=\nu_\pi(\{(0,0,\ldots)\}\ge 0$ and $c_j=\nu_\pi(\{(1,0,\ldots)\})\ge 0$. The corresponding modified paintboxes are $\delta_{\{[j],\{j+1\},\{j+2\},\ldots\}}$ and $\delta_{\{\{j+1\},\bN\setminus\{j+1\}\}}$, respectively, except for $j=1$, where
it is $\frac{1}{2}(\delta_{\{\{1\},\bN\setminus\{1\}\}}+\delta_{\{\{2\},\bN\setminus\{2\}\}})$. We also incorporate the normalisation constants $Z^\pi_\fs$ of the modified paintboxes as densities into $\nu_j$ and use restricted Kingman paintboxes $\kappa_\fs(\cdot\cap\cP^j)$ rather than normalised modified paintboxes $\kappa_\fs^\pi$.

\section{Basic results on restricted exchangeability and related notions}\label{sectprel}

\subsection{Partially exchangeable and constrained exchangeable partitions}\label{partialconstrained}

Let us explore the connections between the RE partitions introduced in this paper and other generalisations of exchangeability studied in the
literature, notably partial exchangeability and constrained
exchangeability.
 \textit{Partially exchangeable
partitions} were introduced by Pitman \cite{Pit-95}. A measure $\mu_n$ on $\cP_n$ is partially exchangeable if 
$\mu_n(\pi)=\mu_n(\pi^\prime)$ for all $\pi,\pi^\prime\in\cP_n$ with the same \em vector \em of block sizes \em in the order of least element. \em Partially
exchangeable measures are not RE, in
general, nor vice versa. Specifically,
$\pi=\left\{\{1, 2\}, \{3, 4\}\right\}$ and $\pi^\prime=\left\{\{1,3\},
\{2,4\}\right\}$ have the same
mass for partially exchangeable measures but not
necessarily for RE measures. Vice versa, consider $\pi=\left\{\{1, 2, 3\}, \{4,5\}\right\}$ and
$\pi^\prime=\left\{\{1,2\},\{3, 4, 5\}\right\}$. In fact, ``the intersection'' of
the two concepts is exchangeability:

\begin{prop}\label{prop10}
A measure $\mu_n$ of $\cP_n$ is exchangeable if and only if
it is both partially exchangeable and RE with $\cC=\{{\bf 0}_{[2]},{\bf 1}_{[2]}\}=\{\{\{1\},\{2\}\},\{\{1,2\}\}\}$.
\end{prop}
\begin{pf}
The ``only if'' part follows straight from the definitions. For the ``if'' part, suppose that $\pi,\pi^\prime\in\cP_n\setminus\{{\bf 1}_{[n]}\}$ have the same multiset of block sizes. Let $\widetilde{\pi}$ be 
such that, for blocks in order of least element, $\widetilde{\pi}_1=(\pi_1\cup\{\min\pi_2\})\setminus\{2\}$ and $\widetilde{\pi}_2=(\pi_2\setminus\{\min\pi_2\})\cup\{2\}$, $\widetilde{\pi}_j=\pi_j$, $j\ge 3$. Similarly construct
$\widetilde{\pi}^\prime$ from $\pi^\prime$. By partial exchangeability $\mu_n(\pi)=\mu_n(\widetilde{\pi})$ and
$\mu_n(\pi^\prime)=\mu_n(\widetilde{\pi}^\prime)$. But $\widetilde{\pi}^\prime,\widetilde{\pi}\in\cP^{\{\{1\},\{2\}\}}$, so
by restricted exchangeability, we have $\mu_n(\widetilde{\pi}^\prime)=\mu_n(\widetilde{\pi})$.
\end{pf}
\em Constrained exchangeable partitions \em were
introduced by Gnedin \cite{Gne-06}. Let
$\varsigma=(\varsigma_k,k\ge 1)$ be a fixed
sequence of integers $\varsigma_k\geq 1$. Consider the set $\cP^{\varsigma-\rm constr}$ of partitions $\Gamma\in\cP$ that are
\textit{constrained with respect to $\varsigma$} in the sense that each block $\Gamma_k$
contains the $\varsigma_k$ least elements of $\bigcup_{j\geq k}\Gamma_j$
for every $k\ge 1$ with $\Gamma_k\neq \varnothing$. A measure $\mu$ on $\cP$ is
\textit{constrained exchangeable} if $\mu(\cP\setminus\cP^{\varsigma-\rm constr})=0$ for some $\varsigma$, and if
$\mu_n(\pi)=\mu_n(\pi^\prime)$ for all $\pi,\pi^\prime\in\{\Gamma|_n\colon \Gamma\in\cP^{\varsigma-\rm constr}\}$ with the same
multiset of block sizes and all $n\ge 1$. 
For $\varsigma=(1,2,1,\ldots)$, under a constrained
exchangeable measure, $\pi=\{\{1,3\}, \{2, 4\},
\{5\}\}$ and $\pi^\prime=\{\{1,2\}, \{3, 4\}, \{5\}\}$ have the same mass, but not necessarily under a 
RE measure. Vice versa, 
restrictions to $\cP^{\{[j],\{j+1\}\}}$ of a RE measure $\mu$ are constrained
exchangeable if we take $\varsigma=(j,1,1,\ldots)$, but as soon as $\mu$ gives positive mass to more than one $\cP^{\{[j],\{j+1\}\}}$, $j\ge 1$, constrained exchangeability in Gnedin's sense fails.

\subsection{RE hierarchies and fragmentation processes}\label{fragcoal}


In Section \ref{intro1}, we defined hierarchies $\cH_B$ on sets $B\subseteq\bN$ and represented hierarchies on finite $B\subset\bN$ as graph-theoretic trees above a root $\varnothing$, with edges between each block $A\in\cH_B$, $\#A\ge 2$, and its maximal subsets in $\cH_B$, which form a partition of $A$. 
For infinite $B\subseteq\bN$, the notion of a maximal subset $A$ of $B$ in $\cH_B$ is more delicate, and it is not always true that there are maximal subsets that form a partition of $B$.

For a hierarchy $\cH_B$ on infinite $B\subseteq\bN$, we say $\cH_B$ is \em closed \em if for all sequences $(A_j,j\ge 1)$ in $\cH_B$ that are increasing for the inclusion partial order, we have $\bigcup A_j\in\cH_B$, and if for all decreasing sequences we have $\bigcap A_j\in\cH_B$. A closed hierarchy $\cH_B$ is uniquely determined by its restrictions $\cH_B\cap[n]$, $n\ge 1$, as $\cH_B=\{A\subseteq B\colon A\cap[n]\in\cH_B\cap[n]\mbox{ for all $n\ge 1$}\}$. For every hierarchy $\cH_B$ there is a \em closure \em $\cH^{\rm cl}_B$, the intersection of all closed hierarchies containing $\cH_B$.  

Recall from Section \ref{intro1} definitions of labelled Markov branching models $(T_n,n\ge 1)$ with splitting rules $P_n$, $n\ge 2$, and that hierarchies $T_n$ on $[n]$ are called consistent if $T_{n+1}\cap[n]=T_n$,  
Both consistency and the labelled Markov branching property can be viewed as properties of the distributions $Q_n$ of $T_n$, $n\ge 1$. This labelled Markov branching property implies a Markov branching property  \cite{HMPW} for rooted delabelled trees $T_n^\circ\sim Q_n^\circ$, as follows: call \em size \em the number of leaves (non-root degree-1 vertices), \em first split \em the decreasing sequence of subtree sizes for the vertex adjacent to the root; conditionally given that the first split of $T_n^\circ$ is $n_1\ge\cdots\ge n_k$, the subtrees are distributed as if they were independent with respective distributions $Q_{n_i}^\circ$, $1\le i\le k$. On the other hand, the associated family $(Q_n^\circ,n\ge 1)$ will not, in general, have the sampling consistency property of \cite{HMPW}, which asserts that a tree $T_n^\circ\sim Q_n^\circ$ with a leaf picked uniformly at random removed (together with any resulting degree-2 vertex) has distribution $Q_{n-1}^\circ$, for $n\ge 3$. 

Let $\cP^j=\cP^{\{[j],\{j+1\}\}}$, $j\ge 1$. We say that a splitting rule $P_n$ is \em RE \em if for all $1\le j\le n-1$ and $\pi,\pi^\prime\in\cP^j_n:=\{\pi\in\cP_n\colon\pi\cap[j+1]=\{[j],\{j+1\}\}\}$ with the same multiset of block sizes, we have $P_n(\{\pi\})=P_n(\{\pi^\prime\})$. The alpha-gamma model of Section \ref{intro1} is an example. 

If $\kappa$ is a RE
dislocation measure as in Corollary \ref{cor5}, then \vspace{-0.1cm} 
\begin{equation}\label{split}P_n(\{\pi\})=\left\{\begin{array}{ll}\kappa(\cP^\pi)/\kappa(\cP\setminus\cP^{{\bf 1}_{[n]}}),&\kappa(\cP\setminus\cP^{{\bf 1}_{[n]}})>0,\\
  \delta_{{\bf 0}_{[n]}}(\{\pi\}),&\kappa(\cP\setminus\cP^{{\bf 1}_{[n]}})=0\mbox{ or $n=2$,}\end{array}\right.
\quad\mbox{ $\pi\in\cP_n\setminus\{{\bf 1}_{[n]}\}$, $n\ge 2$,}\vspace{-0.1cm}
\end{equation}
defines RE splitting rules and hence inductively a consistent Markov branching model $(Q_n,n\ge 1)$ that we also refer to as \em RE\em. More specifically, there is always $n_0\in\{2,3,\ldots\}\cup\{\infty\}$ such that the second line in (\ref{split}) applies for $n<n_0$ but not for $n\ge n_0$. The second line leads to the \em minimal \em hierarchy $Q_n(\{\{[n],\{1\},\ldots,\{n\},\varnothing\}\})=1$ of $[n]$. We have $P_{n_0}(\{[n_0-1],\{n_0\}\})=1$ \em degenerate\em, while for all $n\ge n_0+1$, we have $P_{n}(\{[n-1],\{n\}\})<1$, \em non-degenerate\em, if $n_0<\infty$. 

Let us call a consistent RE Markov branching model $(Q_n,n\ge 1)$ with splitting rules $(P_n,n\ge 2)$ \em regular \em if there is $n_0\ge 2$ such that $Q_n$ is minimal for $n<n_0$,
and if $P_n$ is degenerate for $n=n_0$, non-degenerate for $n\ge n_0+1$. 
\begin{prop}\label{prop11} All regular consistent labelled Markov branching models $(Q_n,n\ge 1)$ with RE splitting rules $(P_n,n\ge 2)$ 
  are of the form {\rm(\ref{split})} for some RE measure $\kappa$ as in Corollary \ref{cor5}.
\end{prop}
\begin{pf} In Pitman's \cite{csp} formalism of exchangeable partition probability functions (EPPFs)
  $$p_n^j(\#\pi_1,\ldots,\#\pi_k)=P_n(\{\pi\}),\qquad\pi\in\cP_n^j=\cP^{\{[j],\{j+1\}\}}\cap[n],j\in[n-1],$$
  consistency in the RE case (extending \cite[Formula (16)]{MPW}) is equivalent to \vspace{-0.2cm} $$p_n^j(n_1,\ldots,n_k)=p_{n+1}^n(n,1)p_n^j(n_1,\ldots,n_k)+\sum_{i=1}^{k+1}p_{n+1}^j(n_1,\ldots,n_{i-1},n_i+1,n_{i+1},\ldots,n_k)\vspace{-0.2cm}$$
  for all $n_1,\ldots,n_k\in\bN$, $k\ge 2$, $n=n_1+\cdots+n_k$, $j\in[n-1]$.
  For $\lambda_n=0$, $n<n_0$, any $\lambda_{n_0}\in(0,\infty)$ and $(1-p_{n+1}^n(n,1))\lambda_{n+1}=\lambda_n$, $n\ge n_0$, we see that 
  $\kappa(\cP^\pi)=\lambda_nP_n(\{\pi\})$, $\pi\in\bigcup_{n\ge 2}\cP_n\setminus\{{\bf 1}_{[n]}\}$, 
  defines a RE measure that has the properties required.
\end{pf}

By Kolmogorov's consistency theorem, we can consider a consistent family $(T_n,n\ge 1)$ of trees $T_n\sim Q_n$ with $T_{n+1}\cap[n]=T_n$, $n\ge 1$, and associate $\cH=\{A\subset\bN\colon A\cap[n]\in T_n\mbox{ for all $n\ge 2$}\}$ as random closed hierarchy on $\bN$, which we call \em RE \em if $(Q_n,n\ge 1)$ is RE. 
In the regular RE case with $n_0=2$, we can consistently embed into continuous time the blocks of $T_n$, $n\ge 2$, using\vspace{-0.1cm}
\begin{itemize}\item exponential holding times $\eta_{[n]}$ of rate $\lambda_n$ for state 
    $[n]$, $\lambda_n$ as in the proof of Proposition \ref{prop11};\vspace{-0.2cm}
  \item recursively and independently as blocks appear from splits, $\eta_{\pi}$ at rate $\lambda_{\#\pi}$ for any 
    $\pi\in T_n$.\vspace{-0.2cm}
\end{itemize}  
With the convention that $\lambda_1=0$ gives infinite holding times, the collection of blocks held at any given time $t\ge 0$ forms a partition $F|_n(t)$ of $[n]$. Indeed, this construction
yields consistent homogeneous fragmentation processes $(F|_n(t),t\ge 0)$ in $\cP_n$, $n\ge 2$, that determine a $\cP$-valued process $(F(t),t\ge 0)$, which we call a RE homogeneous fragmentation process. 

We can also generalise Bertoin's \cite{Ber-book} Poissonian construction to directly obtain RE homogeneous fragmentation processes $(F(t),t\ge 0)$ in $\cP$ from a RE dislocation measure $\kappa$. This provides an alternative construction of the same random closed hierarchy $\cH=\{F_i(t),i\ge 0,t\ge 0\}^{\rm cl}$, but we do not need this alternative construction and leave the details to the reader. 

In the regular case with $n_0\ge 3$, a block $[n_0-1]$ is never split under the exponential-rates construction above; informally $[n_0-1]$ is a limiting block at infinity alongside many other such blocks of size $n_0-1$ that still need splitting to obtain a hierarchy -- they need partitioning into singletons. The simplest kind of irregular model of a RE hierarchy can be obtained here by some intermediate partitioning of these blocks of size $n_0-1$. 
It is possible, but not as natural as in the regular case with $n_0=2$, to incorporate such
further splits in a common embedding, also when other irregularities occur with more degenerate splitting rules. For our next aim of embedding hierarchies into self-similar CRTs, such embeddings do not provide a suitable framework.

\section{Embedding in self-similar CRTs, proof of Theorem \ref{thm2}} \label{sec4}

\subsection{Self-similar CRTs, fragmentation processes and spinal decomposition}\label{fragCRT}

Aldous \cite{Ald-91} called a pair $(\cT,\mu)$ a
\em continuum tree \em if $\cT$ is an $\bR$-tree, $\mu$ a finite measure
on $\cT$, with\vspace{-0.1cm}
\begin{enumerate}
   \item the measure $\mu$ supported by the set ${\rm Lf}(\cT)$ of leaves of $\cT$,\vspace{-0.2cm}
   \item the measure $\mu$ has no atoms,\vspace{-0.2cm}
   \item for every $x\in\cT\backslash{\rm Lf}(\cT)$, positive mass
   $\mu(\cT_x)>0$ in the subtree $\cT_x$ rooted at $x$.\vspace{-0.1cm}
\end{enumerate}
We specify a root vertex $\rho\in\cT$ and distance function $\rd$. For technical simplicity, we follow Aldous \cite{Ald-CRT3} and use CRTs in
$\ell_1=\ell_1(\bN)$. We endow the set of compact subsets of $\ell_1$ with the Hausdorff metric, and the set of finite measures on $\ell_1$ with any
metric inducing the topology of weak convergence, so that the set $\bH$ of
pairs $(T,\mu)$ where $T$ is a rooted $\mathbb{R}$-tree embedded as
a subset of $\ell_1$ and $\mu$ is a finite measure on $T$, is endowed with
the product Borel $\sigma$-algebra. 

A \em Continuum Random Tree (CRT) \em is a random variable with values in the set of continuum trees. To be specific, we call \em distribution of a CRT 
$(\cT,\mu,\rho,\rd)$ \em the distribution on $\bH$ of the particular random isometric embedding of $(\cT,\rd)$ in $\ell_1$ obtained from a random sample $\Sigma_i^*$, $i\ge 1$, of independent leaves with distribution $\mu/\mu(\cT)$, using $0\in\ell_1$ as the root and the $i$th coordinate direction in $\ell_1$ to embed the branch leading to leaf $\Sigma_i^*$, finally passing to the $\ell_1$-closure and the weak limit of the $\mu(\cT)$-multiples of empirical measures of the embedded $\Sigma_1^*,\ldots,\Sigma_i^*$, $i\ge 1$. 

For, $\alpha\in\bR$, $x\in[0,1]$ and $\fs\in S^\downarrow$, we denote by $Q_x^\alpha$ the distribution of the $\alpha$-scaled tree $(\cT,x\mu,\rho,x^\alpha\rd)$ and by $Q_\fs^\alpha$ the distribution of a bush of
independent trees with distributions $Q_{s_i}^\alpha$, $i\ge 1$, all grafted to the
same root. For every $u\geq 0$, consider the bush $\cB(u)$ obtained by grafting the connected components $\mathcal{T}_i(u)$, $i\in I$, of the open set
$\{x\in\cT\colon  \rd(x,\rho)>u\}$ to the same root. Recall that a CRT is called $\alpha$-self-similar
in the sense of \cite{HM}, if for all $u\ge 0$ and conditionally given $(\mu(\mathcal{T}_i(u)),i\in I)^\downarrow=\fs\neq 0$, we have
$\cB(u)\sim Q_\fs^\alpha$.

For $\alpha\in\bR$, a $\mathcal{P}$-valued process $(\Pi(t),t\geq0)$ is an \em exchangeable $\alpha$-self-similar fragmentation process \em if $\Pi=(\Pi(t),t\geq0)$ is exchangeable and if given
$\Pi(t)=\pi$, the partition $\Pi(t+s)$ has the same law as the
random partition whose blocks are those of
$\pi_i\cap\Pi^{(i)}(|\pi_i|^{-\alpha} s),i\geq 1$, where $(\Pi^{(i)},
i\geq1)$ is a sequence of i.i.d.\ copies of $\Pi$. The process
$X=(|\Pi(t)|^\downarrow,t\geq 0)$ is an $S^\downarrow$-valued $\alpha$-self-similar fragmentation. 
Bertoin proved in \cite{Ber-lp} that the
distribution of an exchangeable $\mathcal{P}$-valued self-similar fragmentation is
determined by a triple $(\alpha,c,\nu)$, where $\nu$ is a
dislocation measure on $S^\downarrow$, i.e. $\nu(s_1=1)=0$ and $\int_{S^\downarrow}(1-s_1)\nu(d\fs)<\infty$. In this paper, we take $c=0$ and $\nu$ \em conservative\em, i.e. $\nu(s_0>0)=0$, where $s_0=1-\sum_{i\ge 1}s_i$. We call
$(\alpha,\nu)$ \em characteristic pair\em.

According to
\cite{HM}, 
there exists a self-similar CRT $\cT$ associated with $(\alpha,\nu)$, provided also that $\alpha>0$ (and $\nu$ is infinite, but this is not essential unless it is required that the topological support of $\mu$ is $\cT$). 
Specifically, $Y=((\mu(\cT_i(u)),i\in I_u)^\downarrow,u\ge 0)$ has the same distribution 
as $X$.



Consider $\fs\in S^\downarrow$ and $\fs^{(i)}\in S^\downarrow$, $i\ge 1$. We call \em fragmentation of $\fs$ by $\fs^{(\cdot)}$ \em the mass partition ${\rm Frag}(\fs,\fs^{(\cdot)})$
given by the decreasing rearrangement of $(s_is_j^{(i)}, i,j\in\bN)$. Bertoin
showed that the process $(X(t), t\geq0)$ is Markovian and its
semigroup can be described as follows. For every $t, t^\prime\geq 0,$ the
conditional distribution of $X(t+t^\prime)$ given $X(t)=\fs$
is the law of ${\rm Frag}(\fs, \fS^{(\cdot)}),$ where each
$\fS^{(i)}$ independently is distributed as $X(t^\prime s_i^{-\alpha})$, see
\cite[Proposition 3.7]{Ber-book}. 

Consider an infinite block $B\subseteq\bN$ and $\Gamma\in\cP$. We call \em fragmentation of $B$ by $\Gamma$ \em the partition ${\rm Frag}(B,\Gamma):=\beta_B(\Gamma)\in\cP_B$, where $\beta_B$ is the unique increasing bijection from $\bN$ to $B$. This is a slight variation of Bertoin's \cite{Ber-book} notion, who uses $\Gamma\cap B$, not  $\beta_B(\Gamma)$, but this is useful in Lemma 
\ref{PPP} as it allows to recover $\Gamma$ from ${\rm Frag}(B,\Gamma)$, and it is also
instructive in the RE case.



\begin{figure}[t]
\psfrag{Tvmu}{$\!\!\!\!\!\!\!\!\!\cT_{(v)}(u-)$}
\psfrag{Tvu}{$\!\!\!\!\!\!\!\!\!\cT_{(v)}(u)$}
\psfrag{vvv}{$v$}
\psfrag{vau}{$\begin{array}{rcl}v(u)\!\!\!&\!\!=\!\!&\!\!\!v(\eta_v(t))\\&=&\!\!\!v^{\alpha\rightarrow 0}(t)\end{array}$}
\psfrag{ttt}{$t$}
\psfrag{Bau}{$\cB_{(v)}^{\alpha\rightarrow 0}(t)$}
\psfrag{Xvtm}{$X_{(v)}^{\alpha\rightarrow 0}(t-)$}
\psfrag{zvm}{$\!\!\!\!\zeta_v^{\alpha\rightarrow 0}$}
\psfrag{rho}{$\rho$}
\psfrag{ddd}{${\rm d}(\rho,v(u))=u$}
\psfrag{000}{$0$}

\epsfxsize=13cm\vspace{-0.2cm}
\begin{center}\epsfbox{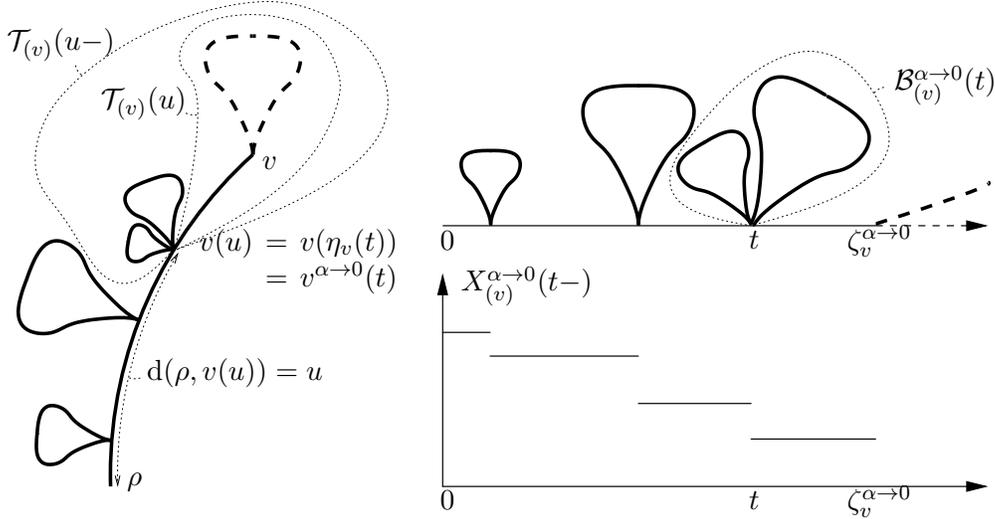}\vspace{-0.1cm}

\begin{minipage}{15cm}\caption{{\small The tree $\cT$ on the left-hand side has its spine to $v\in\cT$ exposed; for a vertex
  $v^\alpha(u)$ on the spine $[[\rho,v]]$, the subtree containing $v$ has been indicated. The 
  right-hand side displays rescaled subtrees and the residual mass process after passage to homogeneous time via $\eta_v$.\vspace{-0.8cm}}}\label{spinefig}\end{minipage}
  \end{center}
\end{figure}

Given a CRT $(\cT,\mu,\rho,\rd)$ and $v\in\cT$, we denote by
$v(u)$ the point on the \em spine \em $[[\rho,v]]$ with ${\rm
d}(\rho,v(u))=u$, $0\le u\le{\rm d}(\rho,v)$, and obtain a parameterisation 
$[[\rho,v]]=\{v(u),0\le u\le{\rm d}(\rho,v)\}$ \em by distance\em, cf. Figure \ref{spinefig}. We consider the subtree $\cT_{(v)}(u)=\{w\in\cT\colon 
{\rm d}(\rho, w\wedge v)>u\}$ of $\cT$ containing
$v$ rooted at $v(u)$, and its mass $X_{(v)}(u)=\mu(\cT_{(v)}(u))$. For $\alpha>0$, let $\eta_v^{\alpha\rightarrow 0}$ be the \em $\alpha$-self-similar time change \em with\vspace{-0.2cm}
\begin{equation}\label{sst}\eta_v^{\alpha\rightarrow 0}(t)=\inf \left\{u\ge 0\colon  \int_0^u
(X_{(v)}(y))^{-\alpha}dy>t\right\},\quad 0\le t<\zeta_v^{\alpha\rightarrow 0}=\int_0^{{\rm d}(\rho,v)}(X_{(v)}(y))^{-\alpha}dy.\vspace{-0.1cm} \end{equation}
Then
$v^{\alpha\rightarrow 0}(t)=v(\eta^{\alpha\rightarrow 0}_v(t))$, $\cT_{(v)}^{\alpha\rightarrow 0}(t)=\cT_{(v)}(\eta^{\alpha\rightarrow 0}_v(t))$ and
$X_{(v)}^{\alpha\rightarrow 0}(t)=\mu(\cT^{\alpha\rightarrow 0}_{(v)}(t))$ are associated time-changed quantities. In particular, $[[\rho,v[[=\{v^{\alpha\rightarrow 0}(t),0\le t<\zeta_v^{\alpha\rightarrow 0}\}$ is a new parameterisation of the spine, which we call parameterisation \em by time\em. Denote by
$S^v(t)=(S_i^v(t),i\ge 1)\in S^\downarrow$ the sequence
such that $X_{(v)}^{\alpha\rightarrow 0}(t-)S^v(t)$ is the decreasing sequence of
$\mu$-masses of the connected components of $\{w\in\cT\colon  v^{\alpha\rightarrow 0}(t)\in
[[\rho, w[[\}$, also $F_v(t)=X_{(v)}^{\alpha\rightarrow 0}(t)/X^{\alpha\rightarrow 0}_{(v)}(t-)$ the component of $S^v(t)$ corresponding 
to the subtree containing $v$.
Moreover, we denote by\vspace{-0.05cm}
 $$\left(\cB_{(v)}^{\alpha\rightarrow 0}(t),\frac{\mu|_{\cB_{(v)}^{\alpha\rightarrow 0}(t)}}{X^{\alpha\rightarrow 0}_{(v)}(t-)},v^{\alpha\rightarrow 0}(t),\frac{\rd|_{\cB_{(v)}^{\alpha\rightarrow 0}(t)}}{(X_{(v)}^{\alpha\rightarrow 0}(t-))^\alpha}\right),\qquad\mbox{where }\cB_{(v)}^{\alpha\rightarrow 0}(t)=\cT_{(v)}^{\alpha\rightarrow 0}(t-)\setminus\cT_{(v)}^{\alpha\rightarrow 0}(t)\vspace{-0.05cm}$$
the associated rescaled spinal bush, of mass $1-F_v(t)$, at time $t\ge 0$.

The following lemma is a description in the CRT framework of Bertoin's tagged particle process that is a bit richer than often stated, but follows from the same arguments.

\begin{lemm}\label{PPP} Let $(\cT,\mu,\rho,\rd)$ be an
  $\alpha$-self-similar CRT with characteristic pair $(\alpha,\nu)$ and $\Sigma^*\sim\mu$. Then
$(S^{\Sigma^*},F_{\Sigma^*})$ is a Poisson point process on
$S^\downarrow\times(0,1)$ with intensity measure $\widetilde{\nu}^*$ given by
$$\widetilde{\nu}^*(d\fs,dx)=\sum_{i\ge 1}s_i\delta_{s_i}(dx)\nu(d\fs).$$
\end{lemm}
\begin{pf} Let $Y$ be the self-similar mass-fragmentation
process associated with the CRT $(\cT, \mu)$ and $Y^{\alpha\rightarrow 0}$ the 
homogeneous mass-fragmentation process obtained by applying the $\alpha$-self-similar
time-change to each block: 
$Y^{\alpha\rightarrow 0}(t)=(\mu(\cT^{\alpha\rightarrow 0}_i(t)),i\in I^{\alpha\rightarrow 0}_t)^\downarrow$, 
where $\{\cT^{\alpha\rightarrow 0}_i(t),i\in I^{\alpha\rightarrow 0}_t\}=\{\cT_{(v)}^{\alpha\rightarrow 0}(t)\colon v\in\cT,\zeta_v^{\alpha\rightarrow 0}>t\}$.	 On an extended probability space, denote by $\Pi$ a homogeneous exchangeable $\cP$-valued
fragmentation process associated with $Y^{\alpha\rightarrow 0}$. Without loss of
generality, we can consider $X_{(\Sigma^*)}^{\alpha\rightarrow 0}(t)=|\Pi_1(t)|$, by 
exchangeability. Since $|\Pi_1(t)|>0$ a.s., the block $\Pi_1(t)$ is infinite and there is a 
unique partition $\Pi^{(1)}(t)$ of $\bN$ such
that $\Pi_1(t)={\rm Frag}(\Pi_1(t-),\Pi^{(1)}(t))$. Furthermore, 
$S^{\Sigma^*}(t)=|\Pi^{(1)}(t)|^\downarrow$. By Bertoin's Poissonian
construction of exchangeable fragmentations, $\Pi^{(1)}$ is
a (time-homogeneous) Poisson point process with intensity measure
$\kappa=\int_{\cS^\downarrow}\kappa_\fs\nu(d\fs)$.
Hence, $S^{\Sigma^*}$ is a Poisson point process on $S^\downarrow$
with intensity measure $\nu$.

As $\Sigma^*$ is distributed according to $\mu$, it is not hard to show that the distribution of 
$(S^{\Sigma^*},F_{\Sigma^*})$ can be obtained by marking $S^{\Sigma^*}$ via the 
size-biased marking kernel $K^*(\fs,\cdot)=\sum_{i\ge 1}s_i\delta_{s_i}$ and so $(S^{\Sigma^*},F_{\Sigma^*})$
is a Poisson point process with intensity $K^*(\fs,dx)\nu(d\fs)=\widetilde{\nu}^*(d\fs,dx)$.
\end{pf}

By the stopping line argument of \cite[Proposition 4]{HPW}, this yields the following joint description of the ordered coarse and unordered fine spinal decompositions along the spine to $\Sigma^*\sim\mu$.

\begin{prop}[Spinal decomposition \cite{BeR,HPW}]\label{spindec} Let $(\cT,\mu,\rho,\rd)$ be an
  $\alpha$-self-similar CRT with characteristic pair $(\alpha,\nu)$ and $\Sigma^*\sim\mu$. Then the process $(S^{\Sigma^*},F_{\Sigma^*},\cB_{(\Sigma^*)}^{\alpha\rightarrow 0})$ is a Poisson point process with
  intensity measure
  $$\widetilde{\nu}^*_{\rm bush}(d\fs,dx,dT)=\sum_{i\ge 1}s_i\delta_{s_i}(dx)Q^\alpha_{(s_1,\ldots,s_{i-1},s_{i+1},\ldots)}(dT)\nu(d\fs).$$
  Conversely, the isometry class of $(\cT,\mu,\rho,\rd)$ is a measurable function of 
  $(S^{\Sigma^*},F_{\Sigma^*},\cB_{(\Sigma^*)}^{\alpha\rightarrow 0})$.  
\end{prop}

\subsection{A generic procedure to sample a leaf from a self-similar CRT}

Our aim is to generalise Lemma \ref{PPP} and Proposition \ref{spindec} to leaves 
other than the $\mu$-sampled leaf $\Sigma^*$ where we are effectively marking a Poisson point process with intensity measure $\nu$ using the size-biased marking kernel $K^*(\fs,\cdot)=\sum_{i\ge 1}s_i\delta_{s_i}$ from $S^\downarrow$ to $(0,1)$. We will now consider other marking kernels. It will be convenient to adopt an idea from Pitman's EPPF formalism and specify the probability
that a \em specific \em part of size $x$ is chosen with probability $P(\fs,x)$ so
that the probability of choosing a mass $x$ is $K(\fs,\{x\})=m_xP(\fs,x)$ where for $\fs=(s_i,i\ge 1)\in S^\downarrow$, we let $m_x=\#\{i\ge 1\colon s_i=x\}$.

\begin{definition}\rm\label{dfspf}
A measurable function $P\colon S^\downarrow\times(0,1)\rightarrow [0,1]$ that
fulfils the two conditions 
$$\bullet\mbox{ $P(\fs,x)=0$ if $x\not\in\{s_i,i\ge 1\}$,}\qquad\qquad\qquad \bullet\mbox{ $\sum_{i\ge 1}P(\fs,s_i)=1$,}\qquad\qquad\qquad\qquad\qquad$$ 
is called a selection probability function (SPF). 
\end{definition}

\begin{example}\rm The SPF
associated with a leaf chosen according to $\mu$ is
$P(\fs,s_i)=s_i$.
\end{example}
We now formulate the procedure to sample a special leaf
$\Sigma$ based on an SPF $P$ from an $\alpha$-self-similar CRT $(\cT, \mu,
\rho, {\rm d})$ with dislocation measure $\nu$, $\cT\sim Q^\alpha_1=Q_1$ for short ($\alpha>0$ fixed).

\begin{proc}\rm\label{proc1} Let $P$ be an SPF as in Definition \ref{dfspf}
fulfilling
\begin{equation}\label{finite}\int_{S^\downarrow}\sum_{i\ge 1}(1-s_i)P(\fs,s_i)\nu(d\fs)<\infty.
\end{equation}
\begin{enumerate}
\item[\rm 0.] We start from $(\cT_1, \mu_1, \rho_1, {\rm d}_1):=(\cT, \mu,
\rho,{\rm d})$ and $i=1$ and proceed inductively.

\item[\rm 1.] Conditionally given $(\cT_i, \mu_i, \rho_i, {\rm d}_i)$, let
  $\Sigma_{(i)}\sim\mu_i$.
\item[\rm 2.] Conditionally given
$(\cT_i,\Sigma_{(i)})$,
we consider the parameterisation in homogeneous time of the spine $[[\rho_i,\Sigma_{(i)}[[=\{\Sigma^{\alpha\rightarrow 0}_{(i)}(t),t\ge 0\}$ and
pick as $\cT_{i+1}$ a subtree $\cS$ off the spine; specifically, if $\cS$ is a subtree rooted at the spinal vertex $\Sigma_{(i)}^{\alpha\rightarrow 0}(t)$, it is selected with probability
$$\bP\left( \cT_{i+1}=\cS \left|
\cT_i,\Sigma_{(i)}\right.\right)
=P\left(S^{\Sigma_{(i)}}(t),\frac{\mu_i\left(\cS\right)}{\mu_i\left(\cT^{\alpha\rightarrow 0}_{(\Sigma_{(i)})}(t-)\right)}\right)\prod_{t^\prime<
t}P\left(S^{\Sigma_{(i)}}(t^\prime),F_{\Sigma_{(i)}}(t^\prime)\right).$$

\item[\rm 3.] Let $\tau_{(i)}=\inf\{t\geq 0\colon 
\cT_{i+1}\cap\cT_{\Sigma_{(i)}}^{\alpha\rightarrow 0}(t)=\varnothing\}$. We turn $\cT_{i+1}$ into a CRT with rescaled mass measure, root and rescaled distance function as follows:
$$\mu_{i+1}=\frac{\mu_i|_{\cT_{i+1}}}{\mu_i(\cT_{i+1})},\quad\rho_{i+1}=\Sigma_{(i)}^{\alpha\rightarrow 0}(\tau_{(i)}),\quad{\rm d}_{i+1}=\frac{{\rm
d}_i|_{\cT_{i+1}\times\cT_{i+1}}}{(\mu_i(\cT_{i+1}))^\alpha}.$$
\item[\rm 4.] Repeat within the subtree $(\cT_{i+1},\mu_{i+1},\rho_{i+1},{\rm d}_{i+1})$ by increasing $i$ by $1$ and proceeding to 1.
 \item[\rm 5.] As $i\rightarrow\infty$, we obtain a sequence $(\Sigma_{(i)}^{\alpha\rightarrow 0}(\tau_{(i)}),i\ge 1)$ in $\cT$
  that increases in the sense that $\Sigma_{(i)}^{\alpha\rightarrow 0}(\tau_{(i)})\in[[\rho,\Sigma_{(i+1)}^{\alpha\rightarrow 0}(\tau_{(i+1)})]]$ and hence 
  converges. Let 
 $\Sigma=\lim_{i\rightarrow\infty}\Sigma_{(i)}^{\alpha\rightarrow 0}(\tau_{(i)})$.

\end{enumerate}
\end{proc}
Note that Step 2.\ is well-defined as $\prod_{t^\prime\ge 0}P(S^{\Sigma_{(i)}}(t^\prime),F_{\Sigma_{(i)}}(t^\prime))=0$, by Proposition \ref{spindec}.


\begin{figure}[h]
\psfrag{vvv}{$\Sigma_{(1)}$}
\psfrag{www}{$\Sigma_{(2)}$}
\psfrag{xxx}{$\!\!\Sigma_{(3)}$}
\psfrag{yyy}{$\Sigma$}
\psfrag{T2}{$\cT_2$}
\psfrag{T3}{$\cT_3$}
\psfrag{rho}{$\rho$}
\psfrag{s1t1}{$\Sigma_{(1)}^{\alpha\rightarrow 0}(\tau_{(1)})$}
\psfrag{s2t2}{$\!\!\!\!\!\!\!\!\!\!\Sigma_{(2)}^{\alpha\rightarrow 0}(\tau_{(2)})$}

\epsfxsize=14cm
\begin{center}\epsfbox{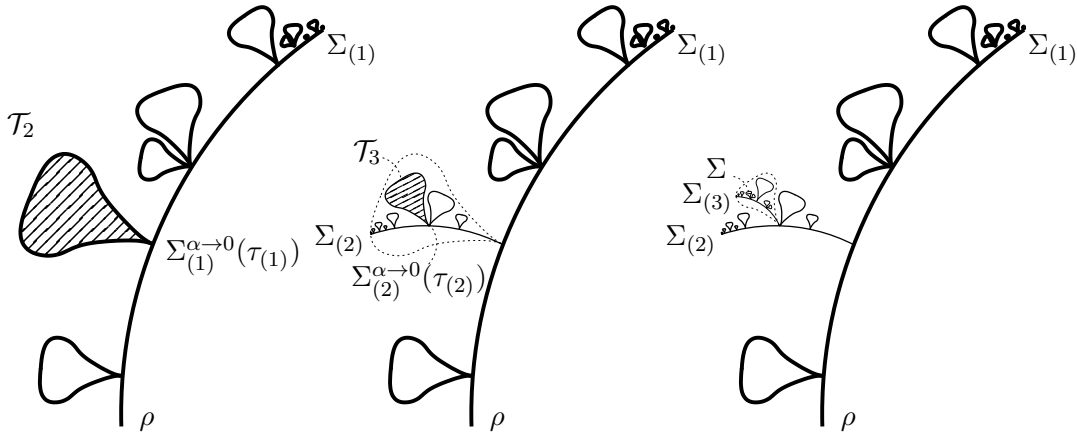}
\begin{minipage}{15cm}\caption{{\small In Procedure 1 we begin by sampling in $\cT_1=\cT$ a leaf $\Sigma_{(1)}\sim\mu$ and pick one of the spinal subtrees as $\cT_2$ according to SPF $P$. Within $\cT_i$ for $i=2$, rescaled, we repeat by sampling a leaf $\Sigma_{(i)}\sim\mu_i$ and pick spinal subtree $\cT_{i+1}$ according to $P$. As $i\rightarrow\infty$ we indicate $\Sigma=\lim_{i\rightarrow\infty}\Sigma_{(i)}^{\alpha\rightarrow 0}(\tau_{(i)})$.}}\label{proc1fig}\vspace{-0.5cm}\end{minipage}\end{center}
\end{figure}

\noindent Roughly speaking, this sampling procedure is that we travel along the
spine $[[\rho, \Sigma_{(1)}]]$ and keep selecting subtrees until the
first time we choose a subtree not containing $\Sigma_{(1)}$ and then
repeat inductively in the subtree until we reach a leaf $\Sigma$ in the limit, see Figure \ref{proc1fig}. 
We show in the following proposition that there is a spinal
subordinator associated with $\Sigma$.

\begin{prop}\label{sleaf} Let $\Sigma$ be sampled according to Procedure \ref{proc1}.
\begin{enumerate}
\item[\rm (i)] The process $(S^\Sigma,F_\Sigma,\cB_{(\Sigma)}^{\alpha\rightarrow 0})$ is a Poisson 
  point process with intensity measure\vspace{-0.1cm}
  $$\widetilde{\nu}^P_{\rm bush}(d\fs,dx,dT)=\sum_{i\ge 1}P(\fs,s_i)\delta_{s_i}(dx)Q_{(s_1,\ldots,s_{i-1},s_{i+1},\ldots)}(dT)\nu(d\fs).\vspace{-0.2cm}$$ 
  Specifically, $\left(\left(S^{\Sigma}(t),F_\Sigma(t),\cB_{(\Sigma)}^{\alpha\rightarrow 0}(t)\right), 0\leq t<\tau_{(1)}\right)$ is a killed Poisson point process with killing rate
  $\int_{S^\downarrow}\sum_{i\ge 1}(1-s_i)P(\fs,s_i)\nu(d\fs)$ and intensity measure\vspace{-0.2cm} 
  $$\widetilde{\nu}^P_{(1),\rm bush}(d\fs,dx,dT)=\sum_{i\ge 1} s_iP(\fs,s_i)\delta_{s_i}(dx)Q_{(s_1,\ldots,s_{i-1},s_{i+1},\ldots)}(dT)\nu(d\fs).\vspace{-0.2cm}$$ 
\item[\rm (ii)]
 Let $\xi^{\Sigma}_t=-\log X_{(\Sigma)}^{\alpha\rightarrow 0}(t)$, $t\ge 0$. Then
$\xi^{\Sigma}$ is a pure jump subordinator with Laplace exponent $\Phi_\Sigma$ and
L\'evy measure $\Lambda_\Sigma$ given by\vspace{-0.2cm}
\begin{equation}\label{sleaf1}\Phi_\Sigma(q)=\int_{\cS^\downarrow}\sum_{i\ge 1}(1- s_i^q)
P(\fs,s_i)\nu(d\fs)\quad\mbox{and}\quad
\Lambda_{\Sigma}=\int_{S^\downarrow}\sum_{i\ge 1}P(\fs,s_i)\delta_{-\log s_i}\nu(d\fs).\vspace{-0.2cm}
\end{equation}
\end{enumerate}
\end{prop}

\begin{pf} 
(i) This proof relies heavily on Poisson point process techniques. We use
the terminology of Kingman \cite{Kin-PP}. By Proposition \ref{spindec}, the process $(S^{\Sigma_{(1)}},F_{\Sigma_{(1)}},\cB^{\alpha\rightarrow 0}_{(\Sigma_{(1)})})$ is a Poisson point process with intensity measure $\widetilde{\nu}^*_{\rm bush}$. Step 2. of Procedure \ref{proc1} can be read and analysed as follows. We mark some points of this Poisson point process with a selected subtree $\cT^{\rm sel}_{(\Sigma_{(1)})}(t)$ using the kernel \vspace{-0.2cm}
$$K(\fs,x,B^\prime;dT^{\prime\prime})=P(\fs,x)\delta_{(\{0\})}(dT^{\prime\prime})
	+\sum_{S\;{\rm connected\;component\;of\; }B^\prime\setminus\{\rho^\prime\}}P(\fs,\mu^\prime(S))\delta_S(dT^{\prime\prime}),\vspace{-0.1cm}$$
where $B^\prime$ is short for $(B^\prime,\mu^\prime,\rho^\prime,\rd^\prime)$ and
$T^{\prime\prime}$ is short for $(T^{\prime\prime},\mu^{\prime\prime},\rho^{\prime\prime},\rd^{\prime\prime})$, also
$S$ for $(S,\mu^\prime|_S,\rho^\prime,\rd^\prime|_{S\times S})$ and $\{0\}$
for $(\{0\},0,0,0)$. 
By standard marking and mapping, we get a new Poisson point process $(S^{\Sigma_{(1)}},F_{\Sigma_{(1)}},\cB_{(\Sigma_{(1)})}^{\rm rem},\cT^{\rm sel}_{(\Sigma_{(1)})})$, where $\cB_{(\Sigma_{(1)})}^{\rm rem}(t)=\cB_{(\Sigma_{(1)})}^{\alpha\rightarrow 0}(t)\setminus\cT_{(\Sigma_{(1)})}^{\rm sel}(t)$ with 
intensity measure\vspace{-0.2cm}
$$\sum_{i\ge 1}s_i\delta_{s_i}(dx)\left(P(\fs,s_i)Q_{\widehat{\fs}^{(i)}}(dB^\prime)\delta_{\{0\}}(dT^{\prime\prime})+\sum_{j\neq i}P(\fs,s_j)Q_{\widehat{\fs}^{(i,j)}}(dB^\prime)Q_{s_j}(dT^{\prime\prime})\right)\nu(d\fs),\vspace{-0.2cm}$$
where $\widehat{\fs}^{(i)}=(s_1,\ldots,s_{i-1},s_{i+1},\ldots)$ is the sequence $\fs$ with $s_i$ removed and similarly $\widehat{\fs}^{(i,j)}$ is the sequence $\fs$ with $s_i$ and $s_j$ removed.
 
In Step 3., we set $\tau_{(1)}=\inf\{t\ge 0\colon \cT^{\rm sel}_{(\Sigma_{(1)})}(t)\neq\{0\}\}$, 
exponentially distributed with rate\vspace{-0.2cm}
$$\int_{S^\downarrow}\sum_{i\ge 1}s_i\sum_{j\neq i}P(\fs,s_j)\nu(d\fs)=\int_{S^\downarrow}\sum_{j\ge 1}(1-s_j)P(\fs,s_j)\nu(d\fs)<\infty,\vspace{-0.2cm}$$
note (\ref{finite}). Standard thinning and projecting yields that $((S^{\Sigma}(t),F_{\Sigma}(t),\cB^{\alpha\rightarrow 0}_{(\Sigma)}(t)),0\le t<\tau_{(1)})=((S^{\Sigma_{(1)}}(t),F_{\Sigma_{(1)}}(t),\cB^{\rm rem}_{(\Sigma_{(1)})}(t)),0\le t<\tau_{(1)})$ is an independently killed Poisson point process with intensity measure $\sum_{i\ge 1}s_iP(\fs,s_i)\delta_{s_i}(dx)Q_{\widehat{\fs}^{(i)}}(dB^\prime)\nu(d\fs)$, as required for the second assertion. The rescaled tree $\cT_2=\cT^{\rm sel}_{(\Sigma_{(1)})}(\tau_{(1)})\sim Q_1$ is independent of this killed
Poisson point process and also jointly independent of the pair formed by the bush $\cB^{\rm rem}_{(\Sigma_{(1)})}$ and the rescaled tree $\cT^{\alpha\rightarrow 0}_{(\Sigma_{(1)})}(\tau_{(1)})$ that has distribution $Q_{s_i}$ for $s_i=F_{\Sigma_{(1)}}(\tau_{(1)})$, using the converse statement in Proposition \ref{spindec}. 

In Step 4., the induction proceeds on $\cT_i\sim Q_1$, $i\ge 2$, all independent of
the past, so this Poisson point process extends indefinitely, but ignores points at
$\tau_{(1)}+\cdots+\tau_{(i)}$, $i\ge 1$.\linebreak[4] 
These are exponentially spaced and i.i.d., hence form an independent Poisson point process. The independence and distributional properties that we noted identify the
distribution of \linebreak[4] $(S^\Sigma(\tau_{(1)}),F_\Sigma(\tau_{(1)}),\cB^{\alpha\rightarrow 0}_{(\Sigma)}(\tau_{(1)}))=(S^{\Sigma_{(1)}}(\tau_{(1)}),\mu^{\rm sel}(\cT^{\rm sel}_{(\Sigma_{(1)})}(\tau_{(1)})),\widetilde{\cB}_1)$, and the
intensity measure\vspace{-0.2cm}
$$\sum_{i\ge 1}s_i\sum_{j\neq i}P(\fs,s_j)\delta_{s_j}(dx)Q_{(s_i,\widehat{\fs}^{(i,j)})}(dB^\prime)\nu(d\fs)=\sum_{j\ge 1}(1-s_j)P(\fs,s_j)\delta_{s_j}(dx)Q_{\widehat{\fs}^{(j)}}(dB^\prime)\nu(d\fs),\vspace{-0.2cm}$$
because we define $\widetilde{\cB}_i$ by grafting to the same root $\cB^{\rm rem}_{(\Sigma_{(i)})}(\tau_{(i)})\sim Q_{\hat{\fs}^{(i,j)}}$ and the rescaled $\cT^{\alpha\rightarrow 0}_{(\Sigma_{(i)})}(\tau_{(i)})$ has distribution $Q_{s_i}$. Standard superposition completes the proof of (i).

(ii) By (i) and standard mapping, $(\Delta\xi_t^\Sigma,t\ge 0)$ is a Poisson point 
process with intensity measure $\Lambda_\Sigma$, hence $\xi_t^\Sigma=\sum_{s\le t}\Delta\xi_s^\Sigma$ is a pure jump subordinator with Laplace exponent $\Phi_\Sigma$.
\end{pf}


\subsection{A procedure to sample a sequence of leaves from a self-similar CRT}\label{sample}

In this section, we formulate a special inductive procedure to sample $k$ leaves
$\Sigma_1,\ldots, \Sigma_k$ from a self-similar CRT $(\cT, \mu)$ with
characteristic pair $(\alpha,\nu)$, where\vspace{-0.2cm}
$$\nu(d\fs)=\sum_{j\ge 1}\left(\sum_{i\ge 1}s_i^j(1-s_i)\right)\nu_j(d\fs)\vspace{-0.2cm}$$
for some measures $\nu_j$, $j\ge 1$, representing a RE dislocation measure as in Corollary \ref{cor5}.
Clearly, the measures $\nu_j$, $j\ge 1$, are absolutely continuous with respect to
$\nu$. We 
denote their Radon-Nikodym derivatives by $f_j=d\nu_j/d\nu$, $j\ge 1$, and
define selection functions\vspace{-0.2cm}
\begin{eqnarray*}P_0(\fs,s_i)=\sum_{\ell\ge 1}s_i^\ell(1-s_i)f_\ell(\fs),\quad&& 
  P^{\rm old}_{k}(\fs,s_i)=\frac{\sum_{\ell\ge k+1}s_i^\ell(1-s_i)f_\ell(\fs)}{\sum_{\ell\ge k}s_i^\ell(1-s_i)f_\ell(\fs)},\quad k\ge 1,\\ \mbox{and for $j\neq i$}&& P^{\rm new}_{k}(\fs,s_i,s_j)=\frac{s_i^ks_j f_k(\fs)}{\sum_{\ell\ge k}s_i^\ell(1-s_i)f_\ell(\fs)},\quad k\ge 1.\vspace{-0.1cm}\end{eqnarray*}

\begin{proc}\label{proc2}\rm\begin{enumerate}\item[(0)] To sample $\Sigma_1$
  in the whole CRT 
$(\cT_{1,\varnothing},\mu_{1,\varnothing},\rho_{1,\varnothing},\rd_{1,\varnothing})=(\cT,\mu,\rho,\rd)$ we use
  step ($k$,$\varnothing$) for $k=1$ and then proceed inductively.\vspace{-0.1cm}

  \item[($k$,$\varnothing$)] Sample leaf $\Sigma_k$ in $\cT_{k,\varnothing}$
    according to Procedure \ref{proc1} using the SPF $P_0$. Then 
    increase $k$ by $1$, set $B=[k-1]$ and 
    $\cT_{k,B}=\cT$, and proceed to step ($k$,$B$).\vspace{-0.1cm}
  
  \item[($k$,$B$)] with $B\neq\varnothing$.\vspace{-0.2cm}\begin{enumerate}\item[1.] Given $\Sigma_i\in\cT_{k,B}$, $i\in B$, denote by $v_{k,B}$ the branch point that separates the labels in $B$ into several subtrees, so that
	$\displaystyle[[\rho_{k,B},v_{k,B}]]=\bigcap_{i\in B}[[\rho_{k,B},\Sigma_i]]$.\vspace{-0.1cm}
    \item[2.] Conditionally given $(\cT_{k,B};\Sigma_i,i\in B)$, with spine $[[\rho_{k,B},v_{k,B}[[=\{v_{k,B}^{\alpha\rightarrow 0}(t),0\le t<\zeta_{v_{k,B}}^{\alpha\rightarrow 0}\}$, pick as $\cT_{k,B^\prime}$ either a new subtree $\cS$ above some $v_{k,B}^{\alpha\rightarrow 0}(t)$ with probability\vspace{-0.2cm}
\begin{eqnarray*}&&\bP\left( \cT_{k,B^\prime}=\cS \Big|\cT_{k,B};\Sigma_i,i\in B
\right)\\
&&\hspace{0.5cm}=P_{\#B}^{\rm new}\left(S^{v_{k,B}}(t),F_{v_{k,B}}(t),\frac{\mu_{k,B}\left(\cS\right)}{\mu_{k,B}\left(\cT^{\alpha\rightarrow 0}_{(v_{k,B})}(t-)\right)}\right)\prod_{t^\prime<
t}P_{\#B}^{\rm old}\left(S^{v_{k,B}}(t^\prime),F_{v_{k,B}}(t^\prime)\right),\vspace{-0.2cm}
\end{eqnarray*}
or, in the case $\#B\ge 2$, a new or old subtree $\cS$ above $v_{k,B}$ with probability\vspace{-0.2cm}
$$\bP\left( \cT_{k,B^\prime}=\cS \Big|\cT_{k,B};\Sigma_i,i\in B\right)
=\frac{\mu_{k,B}\left(\cS\right)}{\mu_{k,B}\left(\!\cT^{\alpha\rightarrow 0}_{(v_{k,B})}(\zeta_{v_{k,B}}\!-)\!\right)}\prod_{t^\prime<\zeta_{v_{k,B}}}\!\!\!\!\!P_{\#B}^{\rm old}\left(
S^{v_{k,B}}(t^\prime),F_{v_{k,B}}(t^\prime)\right),\vspace{-0.2cm}$$
where $B^\prime=\{i\in B\colon \Sigma_i\in\cS\}$ and new/old means without/with any $\Sigma_i$, $i\in B$.\vspace{-0.1cm}  
  \item[3.] Let $\tau_{k,B}=\min\{\zeta_{v_{k,B}}^{\alpha\rightarrow 0},\inf\{t\ge 0\colon \cT_{k,B^\prime}\cap\cT^{\alpha\rightarrow 0}_{v_{k,B}}(t)=\varnothing\}\}$. We turn $\cT_{k,B^\prime}$ into a CRT with rescaled mass measure, root and rescaled distance function as follows:\vspace{-0.2cm}
$$\mu_{k,B^\prime}=\frac{\mu_{k,B}|_{\cT_{k,B^\prime}}}{\mu_{k,B}(\cT_{k,B^\prime})},\quad\rho_{k,B^\prime}=v_{k,B}^{\alpha\rightarrow 0}(\tau_{k,B}),\quad{\rm d}_{k,B^\prime}=\frac{{\rm
d}_{k,B}|_{\cT_{k,B^\prime}\times\cT_{k,B^\prime}}}{\left(\mu_{k,B}(\cT_{k,B^\prime})\right)^\alpha}.\vspace{-0.2cm}$$
  \item[4.] Repeat within the subtree $(\cT_{k,B^\prime},\mu_{k,B^\prime},\rho_{k,B^\prime},{\rm d}_{k,B^\prime})$ by proceeding to step ($k$,$B^\prime$).\vspace{-0.1cm}
  \end{enumerate}
  \end{enumerate}
\end{proc}
Note that the probabilities in Step 2. add up to 1 since $\sum_{j\colon j\neq i}s_i^ks_jf_k(\fs)=s_i^k(1-s_i)f_k(\fs)$. From Proposition \ref{sleaf}, we obtain the following by straightforward arguments.

\begin{figure}[t]
\psfrag{S1}{$\Sigma_1$}
\psfrag{S2}{$\Sigma_2=v_{4,\{2\}}$}
\psfrag{S3}{$\Sigma_3$}
\psfrag{s1}{$(4,[3])$}
\psfrag{s2}{$(4,\{2,3\})$}
\psfrag{s3}{$(4,\{2\})$}
\psfrag{rho}{$\rho$}
\psfrag{s4}{$(4,\varnothing)$}
\psfrag{v43}{$v_{4,[3]}$}
\psfrag{v42}{$v_{4,\{2,3\}}$}
\epsfxsize=7cm\vspace{-0.2cm}
\begin{center}\epsfbox{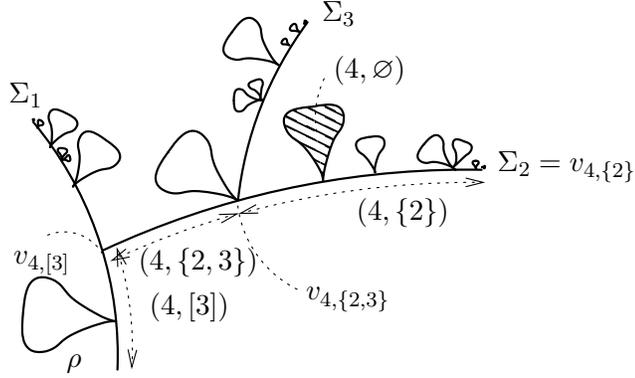}\vspace{-0.2cm}
\begin{minipage}{15cm}\caption{{\small Sampling $\Sigma_4$ in $(\cT;\Sigma_i,i\in[3])$ using Procedure 2: first step is $(4,[3])$, random selection picks the old subtree $\cT_{4,\{2,3\}}$; in step $(4,\{2,3\})$, the old subtree $\cT_{4,\{2\}}$ is picked; in step $(4,\{2\})$, a new subtree, shaded, is picked; step $(4,\varnothing)$ takes place in this subtree, using Procedure 1.\vspace{-0.8cm}}}\label{proc2fig}\end{minipage}\end{center}
\end{figure}

\begin{cor}\label{cor22} Sample $(\Sigma_k,k\ge 1)$ following Procedure \ref{proc2}. Let $v_k$ be the branch point in $\cT$ that separates $[k]$ into different subtrees, $k\ge 1$. Then
$\left(\left(S^{\Sigma_k}(t),F_{\Sigma_k}(t),\cB_{(\Sigma_k)}^{\alpha\rightarrow 0}(t)\right), 0\leq t<\zeta_{v_k}^{\alpha\rightarrow 0}\right)$ is a Poisson point process with killing rate
  $\lambda_k=\int_{S^\downarrow}\sum_{i\ge 1}\sum_{\ell=1}^{k-1}s_i^\ell(1-s_i)\nu_\ell(d\fs)$ and intensity measure\vspace{-0.5cm}
  \begin{equation}\label{intmeask}\widetilde{\nu}^{(k)}_{\rm bush}(d\fs,dx,dT)=\sum_{i\ge 1}\delta_{s_i}(dx)Q_{(s_1,\ldots,s_{i-1},s_{i+1},\ldots)}(dT)\sum_{\ell\ge k}s_i^\ell(1-s_i)\nu_\ell(d\fs).\vspace{-0.1cm}
  \end{equation} 
\end{cor}
Note $\lambda_1=0$, so the Poisson point process is not killed and Corollary \ref{cor22} describes the whole tree $\cT$ jointly with $\Sigma_1$,
decomposed along its spine $[[\rho,\Sigma_1[[$. For $k\ge 2$, Corollary \ref{cor22} describes a spinal decomposition along $[[\rho,v_{k}[[$, but
 not the subtrees above $v_{k}$. This is done in  
Lemma \ref{lm23}.\vspace{-0.1cm}

\begin{pf} The case $k=1$ follows straight from step ($1$,$\varnothing$) of 
  Procedure \ref{proc2} and Proposition \ref{sleaf}. We then proceed by 
  induction in $k$. Assuming that the
  statement is true for $k$, step ($k+1$,$[k]$) 2. and standard thinning with 
  probabilities $P_{k+1}^{\rm old}(\fs,s_i)$ yields\vspace{-0.2cm} 
  $$\widetilde{\nu}^{(k+1)}_{\rm bush}(ds,dx,dT)=\sum_{i\ge 1}P_{k+1}^{\rm old}(\fs,s_i)\delta_{s_i}(dx)Q_{(s_1,\ldots,s_{i-1},s_{i+1},\ldots)}(dT)\sum_{\ell\ge k} s_i^\ell(1-s_i)\nu_\ell(d\fs),\vspace{-0.2cm}$$
  as claimed, and an extra rate $\int_{S^\downarrow}\sum_{i\ge 1}(1-P_{k+1}^{\rm old}(\fs,s_i))\sum_{\ell\ge k}s_i^\ell(1-s_i)\nu_\ell(d\fs)$ is added to the 
  killing rate $\lambda_k$ from the induction hypothesis. This completes the induction\vspace{-0.1cm} step.
\end{pf}

To identify the distribution $Q_1^{[k]}$ of $(\cT;\Sigma_i,i\in[k])$ constructed 
according to Procedure \ref{proc2} run up to some $k\ge 2$, we study its
branching structure recursively by specifying the first branch point $v_k$ that separates $[k]$ into several subtrees denoted by $\cT^{[k]}_\ell$ with label partition $\Pi^{[k]}$ and a remaining bush $\cB_{[k]}$ of unlabelled subtrees, with joint relative subtree sizes $S^{[k]}\in S^{\downarrow}$. For $x\in(0,1]$ and $B=\{b_1,\ldots,b_k\}\subset\bN$ with
$1\le b_1<\cdots<b_k$, it will be convenient to denote by $Q_x^B$ the distribution of
a rescaled and relabelled version of $(\cT;\Sigma_i,i\in[k])$, where the mass measure has been multiplied by $x$, the distance function by $x^\alpha$, and $\Sigma_i$
is renamed to \vspace{-0.1cm} $\Sigma_{b_i}$, $i\in[k]$. %


\begin{lemm}\label{lm23} The first branching of $(\cT;\Sigma_i,i\in[k])$ separating $[k]$ and 
  associated subtrees described in $\cC^{\rm br}_k=(S^{[k]},\Pi^{[k]},\cT^{[k]},\cB_{[k]})$ 
  are independent of $\cC^{\rm pre}_k=((S^{v_k}(t),F_{v_k}(t),\cB^{\alpha\rightarrow 0}_{(v_k)}(t)),0\le t<\zeta_{v_k})$,
  with distribution given by\vspace{-0.2cm}
  \begin{eqnarray*}&&\hspace{-0.5cm}\bP(S^{[k]}\in d\fs,\Pi^{[k]}=\pi,(\cT_1^{[k]};\Sigma_i,i\in\pi_1)\in dT_1,\ldots,(\cT_r^{[k]};\Sigma_i,i\in\pi_r)\in dT_r,\cB_{[k]}\in dB^\prime)\\ 
  &&\hspace{1cm}=\frac{1}{\lambda_k}\left(\sum_{i_1,\ldots,i_r\;\rm distinct}Q_{\widehat{\fs}^{(i_1,\ldots,i_r)}}(dB^\prime)\prod_{\ell=1}^r s_{i_\ell}^{\#\pi_\ell}Q_{s_{i_\ell}}^{\pi_\ell}(dT_\ell)\right)\nu_m(d\fs),\vspace{-0.2cm}
  \end{eqnarray*}
  where $\pi=(\pi_{i_1},\ldots,\pi_r)\in\cP_k$ and $m=\min\pi_2-1$, also $\widehat{\fs}^{(i_1,\ldots,i_r)}$ is $\fs$ with $s_1,\ldots,s_{i_r}$ removed.
\end{lemm}
The kernel $\kappa_{\fs,\pi}(dT_1,\cdots,dT_r,dB^\prime)=\sum_{i_1,\ldots,i_r\;\rm distinct}Q_{\widehat{s}^{(i_1,\ldots,i_r)}}(dB^\prime)\prod_{\ell=1}^r s_{i_\ell}^{\#\pi_\ell}Q_{s_{i_\ell}}^{\pi_\ell}(dT_\ell)$ is a fancy paintbox
that equips each block under $\kappa_\fs$ with a tree and embeds the labels for $\pi\in\cK$.\vspace{-0.1cm}

\begin{pf} For $k=1$, this is trivial since $v_1=\Sigma_1$ is a leaf.
  Now suppose that the result holds for all $[j]\subseteq[k]$, and 
  consider $k+1$. In 
  our use of standard Poisson point process arguments as well as in extracting from
  Procedure \ref{proc2} as from Procedure \ref{proc1}, we build on the proof of Proposition \ref{sleaf}.

  For $\pi\in\cP_{k+1}\setminus\{{\bf 1}_{[k+1]}\}$, let $A_\pi=\{\Pi^{[k+1]}=\pi\}$ be 
  the event that $v_{k+1}$ splits $[k+1]$ into $\pi$. The simplest case
  is for $\pi=\{[k],\{k+1\}\}$. By Corollary \ref{cor22}, the  
  decomposition of $\cT$ along the spine $[[\rho,v_k[[$ 
  is given by the Poisson point process 
  $\left(\left(S^{\Sigma_k}(t),F_{\Sigma_k}(t),\cB_{(\Sigma_k)}^{\alpha\rightarrow 0}(t)\right), 0\leq t<\zeta_{v_k}^{\alpha\rightarrow 0}\right)$ with intensity measure (\ref{intmeask}),
  killed at rate
$\lambda_k=\int_{S^\downarrow}\sum_{i\ge 1}\sum_{\ell=1}^{k-1}s_i^\ell(1-s_i)\nu_\ell(d\fs)$. By comparison with the statement of Corollary \ref{cor22} for $k+1$, we see $\bP(A_{\{[k],\{k+1\}\}})=1-\lambda_k/\lambda_{k+1}$. \em Conditionally given $A_{\{[k],\{k+1\}\}}$\em, the distribution of $(S^{\Sigma_k}(\tau_{k+1,[k]}),F_{\Sigma_k}(\tau_{k+1,[k]}),\cB^{\rm rem}_{(\Sigma_k)}(\tau_{k+1,[k]}),\cT^{\rm sel}_{k+1,[k]}(\tau_{k+1,[k]}))$ is\vspace{-0.3cm}
\begin{eqnarray}&&\hspace{-0.5cm}\frac{1}{\lambda_{k+1}-\lambda_k}\sum_{i\ge 1}\sum_{j\neq i}P_k^{\rm new}(\fs,s_i,s_j)\delta_{s_i}(dx)Q_{\widehat{\fs}^{(i,j)}}(dB^\prime)Q_{s_j}(dT^{\prime\prime})\sum_{\ell\ge k}s_i^\ell(1-s_i)\nu_\ell(d\fs)\nonumber\\
	 &&\hspace{1cm}=\frac{1}{\lambda_{k+1}-\lambda_k}\sum_{i\ge 1}\sum_{j\neq i}\delta_{s_i}(dx)Q_{\widehat{\fs}^{(i,j)}}(dB^\prime)Q_{s_j}(dT^{\prime\prime})s_i^ks_j\nu_k(d\fs),\label{firstspl}
\end{eqnarray}\vspace{-0.3cm}

\noindent   
independently of the rescaled $(\cT^{\alpha\rightarrow 0}_{(\Sigma_k)}(\tau_{k+1,[k]});\Sigma_i,i\in[k])$
that has $Q_1^{[k]}$ as conditional distribution given $A_{\{[k],\{k+1\}\}}$. Note also, that the sampling of $\Sigma_{k+1}$ in the rescaled 
$\cT_{k+1,[k]}^{\rm sel}(\tau_{k+1,[k]})$ yields conditional distribution $Q_1^{\{k+1\}}$ given $A_{\{[k],\{k+1\}\}}$, and that by standard thinning arguments these are conditionally independent of $\left(\left(S^{\Sigma_{k+1}}(t),F_{\Sigma_{k+1}}(t),\cB_{(\Sigma_{k+1})}^{\alpha\rightarrow 0}(t)\right), 0\leq t<\zeta_{v_{k+1}}^{\alpha\rightarrow 0}\right)$ given $A_{\{[k],\{k+1\}\}}$. Multiplying by $\bP(A_{\{[k],\{k+1\}\}})$, this yields the result for $\pi=\{[k],\{k+1\}\}$.

  Now consider any other $\pi=\{\pi_1,\ldots,\pi_r\}\in\cP_{k+1}\setminus\{{\bf 1}_{[k+1]}\}$ and write
  $m=\min\pi_2-1\in[k-1]$. Note that also $m=\min\pi_2\cap[k]-1$. 
  By the induction hypothesis, the collections $\cC^{\rm pre}_k$ describing the spine
  to the branch point separating $[k]$, and $\cC^{\rm br}_k$ describing the branching
  and rescaled subtrees, are independent. We read and analyse Step 2. of Procedure 
  \ref{proc2} by marking $\cC^{\rm pre}_k$ as we marked the Poisson point process in 
  the proof of Proposition 
  \ref{sleaf} \em and \em similarly and independently selecting a new or old subtree
  $\cS$ above $v_k$ with probability\vspace{-0.2cm}
  $$\bP\left( \cT^{\rm sel}=\cS \Big|\cT_{k,B};\Sigma_i,i\in B\right)
=\frac{\mu_{k,B}\left(\cS\right)}{\mu_{k,B}\left(\cT^{\alpha\rightarrow 0}_{(v_{k,B})}(t-)\right)}.\vspace{-0.2cm}$$
  Then $A_\pi$ is an intersection of two independent events $A_\pi=A^{\rm pre}_k\cap A^{\rm br}_\pi$ given by\vspace{-0.2cm}
  $$A_k^{\rm pre}=\{\cT^{\rm sel}_{v_k}=\{0\}\mbox{ for all $0\le t<\zeta_{v_k}$}\}\quad\mbox{and }A^{\rm br}_\pi=\{\cL_k(\cT^{\rm sel})=\pi_{(k+1)}\cap[k]\},\vspace{-0.2cm}$$
  where $\cL_k(\cS)=\{i\in[k]\colon \Sigma_i\in\cS\}$ and $\pi_{(k+1)}$ is the block of $\pi$
  containing $k+1$. By construction, $(\cC^{\rm pre}_k,A_k^{\rm pre})$ and 
  $(\cC^{\rm br}_k,A_\pi^{\rm br})$ are also independent and, since the random 
  variables used to sample $\Sigma_{k+1}$ in $\cT^{\rm sel}$ are conditionally 
  independent of $(\cC^{\rm pre}_k,A_k^{\rm pre})$ given $\cT^{\rm sel}$, also
  $\cC^{\rm br}_{k+1}$ is independent of $(\cC^{\rm br}_k,A_\pi^{\rm br})$, hence
  of $\cC^{\rm br}_{k+1}$, since on $A_\pi^{\rm br}$, we have $\cC^{\rm br}_{k+1}=\cC^{\rm br}_k$. The distribution of $\cC^{\rm br}_{k+1}$ now follows from
  the conditional distribution of $\cT^{\rm sel}$ given $\cC^{\rm br}_k$, the
  recursive nature of Procedure \ref{proc2} and the stability of
  the procedure under increasing bijections from $[j]$ to other
  sets $B\subset\bN$ with $\#B=j$ that allows us to apply the induction
  hypothesis to obtain that the sampling of $\Sigma_{k+1}$ in the rescaled 
  $\cT^{\rm sel}\sim Q_1^B$ yields a tree with rescaled distribution $Q_1^{B\cup\{k+1\}}$, as required. 
\end{pf}
Inductively, Lemma \ref{lm23} yields a subtree decomposition of $(\cT;\Sigma_1,\ldots,\Sigma_k)$. 
For $\varnothing\neq B\subset[k]$, consider $[[\rho_{k,B},v_{k,B}[[=\!\!\{v_B^{\alpha\rightarrow 0}(t),0\!<\!t\!<\!\zeta_B^{\alpha\rightarrow 0}\}$ in $\cT_{k,B}\!\subset\!\cT$, \em branch $B$\em,  as in Procedure 2 (cf. Figure \ref{proc2fig}, where $T_3\!=\!\{[3],\{2,3\},\{1\},\{2\},\{3\}\}$). For the rescaled $\cT_{k,B}\sim Q^B_1$, i.e.\  $Q^{[\#B]}_1$ pushed forward under the increasing bijection $[\#B]\rightarrow B$, Corollary \ref{cor22} gives the distribution of the analogous point process $((S^B(t),F_B(t),\cB^{\alpha\rightarrow 0}_B(t)),0\!\le\! t\!<\!\zeta_B^{\alpha\rightarrow 0})$ that captures the spinal subtrees off $[[\rho_{k,B},v_{k,B}[[$. The remaining split at $v_{k,B}$ into relative sizes $S^B$, of which $F^B_i$ is the size corresponding to the $i$th block of the split $\Pi^B$ of $B$ and $\cB_B$ is the subbush of unlabelled subtrees of the remaining sizes in $S^B$, can be read from Lemma \ref{lm23}, as $Q^B_1$ is just a push-forward of $Q^{[\#B]}_1$.

\begin{cor}[Subtree decomposition]\label{cor24} The discrete tree shapes $T_k$, $k\ge 1$, of the reduced trees $R(\cT;\Sigma_1,\ldots,\Sigma_k):=[[\rho,\Sigma_1]]\cup\cdots\cup[[\rho,\Sigma_k]]$, $k\ge 1$, are labelled Markov branching trees with
\begin{equation}\label{splitrul}\bP(\Pi^{[k]}=\pi)=\frac{1}{\lambda_k}\int_{S^\downarrow}\kappa_\fs(\cP^\pi)\nu_m(d\fs),\qquad\mbox{where $m=\min\pi_2-1$.}
\end{equation} 
Conditionally given $T_k=\ft_k$, $\Pi^B=\pi^B=(\pi_1^B,\ldots,\pi_r^B)$ with $m^B=\min\pi_2^B-1$, $B\in\ft_k$, 
\begin{itemize}\item the processes $((S^B(t),F_B(t),\cB^{\alpha\rightarrow 0}_B(t)),0\le t<\zeta_B^{\alpha\rightarrow 0})$, with distribution as in Corollary \ref{cor22},
  \item and the variables $(S^B,F_1^B,\ldots,F_r^B,\cB_B)$, with distribution
$$\frac{1}{\lambda_{\#B}\bP(\Pi^B=\pi^B)}\left(\sum_{i_1,\ldots,i_r\;\rm distinct}Q_{\widehat{\fs}^{(i_1,\ldots,i_r)}}(dB^\prime)\prod_{\ell=1}^r s_{i_\ell}^{\#\pi_\ell^B}\delta_{s_{i_\ell}}(dx_\ell)\right)\nu_{m^B}(d\fs),$$
\end{itemize}
$B\in\ft_k$, are independent. The tree $(\cT;\Sigma_1,\ldots,\Sigma_k)$ with $k$ leaves sampled via Procedure \ref{proc2} is a measurable function of  $(T_k;((F_B,S^B,\cB_B),((S^B(t),F_B(t),\cB^{\alpha\rightarrow 0}_B(t)),0\le t<\zeta_B^{\alpha\rightarrow 0})),B\in T_k)$.   
\end{cor}

\pagebreak[2]

\begin{pfofthm2} We will show that Procedure \ref{proc2} provides an embedding for 
  a RE hierarchy as in Corollary \ref{cor5}, provided that
  $\int_{S^\downarrow}(1-s_1)\nu(d\fs)<\infty$ and $\nu(s_0>0)=c_j=k_j=0$, $j\ge 1$.
  
  A RE hierarchy is uniquely determined by its restrictions
  to $[k]$, $k\ge 1$. But the formula for $\kappa$ in Corollary \ref{cor5} is
  identical to (\ref{splitrul}), hence the hierarchy constructed via Procedure 
  \ref{proc2} is a RE hierarchy associated with $(\nu_j,j\ge 1)$
  embedded in a CRT with characteristic pair $(\alpha,\nu)$.
\end{pfofthm2}

\section{Scaling limits, proof of Theorem \ref{thm3}}

\subsection{Asymptotics of block numbers in Gnedin's constrained partitions}

Before we describe Gnedin's framework and provide a slight extension of his asymptotic study, let us state the $p$th order renewal theory result that we need for this. 

\begin{lemm}[Gut \cite{Gut-74}, Theorem 2.3(b)]\label{renew} Let $N_t\!=\!\#\{n\ge 1\colon X_1\!+\cdots+\!X_n\le t\}$ be the renewal process
  associated with independent and identically distributed $X_j>0$. Then for all $p\in\mathbb{N}$\vspace{-0.2cm}
  $$\mathbb{E}\left[\frac{N_t^p}{t^p}\right]\rightarrow\frac{1}{\left(\bE[X_1]\right)^p}\in[0,\infty),\qquad\mbox{as $t\rightarrow\infty$.}\vspace{-0.1cm}$$
\end{lemm}

\noindent Gnedin \cite{Gne-06} introduced a \em constrained paintbox \em based on an $\bN$-valued  deterministic sequence $\psi=(\psi_k,k\ge 1)$ and a strictly decreasing random sequence $(G_k,k\ge 0)$ in $[0,1]$ with $G_0=1$ and $\lim_{k\rightarrow\infty}G_k=0$. Specifically, he considers a sequence $(I_n,n\ge 1)$ of
independent uniform random variables on $[0,1]$ independent of $(G_k)$, but then associates a modified sequence $(\overline{I}_n^{\psi},n\ge 1)$ that is constrained so that its lower records follow $(G_k,k\ge 1)$ with multiplicities given by $(\psi_k,k\ge 1)$:\vspace{-0.1cm}
\begin{itemize}\item Set $\overline{I}_1^\psi=\cdots=\overline{I}^\psi_{\psi_1}=G_1$; inductively, consider the number $K_n^\psi$ of records $G_k$ that
  have been attained $\psi_k$ times by $(\overline{I}_1^\psi,\ldots,\overline{I}_n^\psi)$,
  and the number $R_n^\psi$ of times that $G_{K_n^\psi+1}$ has been attained
  by $(\overline{I}_1^\psi,\ldots,\overline{I}_n^\psi)$; for $n=\psi_1$, we have 
  $K_n^\psi=1$ and $R_n^\psi=0$; this is the base case.\vspace{-0.1cm}
  \item Given $(\overline{I}_1^\psi,\ldots,\overline{I}_n^\psi)$, $K_n^\psi=k\ge 1$ and $R_n^\psi=r\in\{0,\ldots,\psi_{k+1}-1\}$, proceed as follows\vspace{-0.2cm}
    \begin{itemize}\item if $I_{n+1}\in[G_k,1]$, let $\overline{I}^\psi_{n+1}=I_{n+1}$, $K_{n+1}^\psi=K_n^\psi$ and $R_{n+1}^\psi=R_n^\psi$;\vspace{-0.1cm}
                   \item if $I_{n+1}\in[0,G_k)$ and $r\le\psi_{k+1}-2$, let $\overline{I}^\psi_{n+1}=G_{k+1}$, $K_{n+1}^\psi=K_n^\psi$ and $R_{n+1}^\psi=R_n^\psi+1$;\vspace{-0.1cm}
                   \item if $I_{n+1}\in[0,G_k)$ and $r=\psi_{k+1}-1$, let $\overline{I}^\psi_{n+1}=G_{k+1}$, $K_{n+1}^\psi=K_n^\psi+1$ and $R_{n+1}^\psi=0$.\vspace{-0.1cm} 
    \end{itemize}\pagebreak[2]
\end{itemize}
Eventually, each $G_k$ will appear $\psi_k$ times as lower record in $(\overline{I}_n^\psi,n\ge 1)$. Let $J_n^\psi=K_n^\psi+{\bf 1}_{\{R_n^\psi>0\}}$ be the number of records attained by the $n$ first terms of the sequence. Gnedin obtains the asymptotics of $J_n^\psi$ when
$G_k=Y_1\cdots Y_k$, where $Y_k$, $k\ge 1$, are i.i.d.\ in
$(0,1)$ with $\bE[-\log Y_1]<\infty$ and
$\textrm{Var}(-\log Y_1)<\infty$. Here we drop the requirement of finite logarithmic moments.

\begin{lemm}\label{Gcon}
Let $G_k=Y_1\cdots Y_k$, where $Y_k$, $k\ge 1$, are i.i.d.\ in $(0,1)$. If $\psi=(\psi_k,k\ge 1)$ is such that $\psi_k\in\bN$, $k\ge 1$ and\vspace{-0.3cm}
$$\log\left(\sum_{j=1}^k \psi_j\right)=o(k),\qquad\mbox{as $k\rightarrow\infty$,}\vspace{-0.1cm}$$ 
then\vspace{-0.2cm} $$\lim_{n\rightarrow\infty}
\frac{J_n^\psi}{\log n}=\frac{1}{\bE[-\log Y_1]}$$ in the sense that this limit vanishes when $\bE[-\log
Y_1]=\infty$. Furthermore, for every $p\geq 1$,\vspace{-0.2cm} $$
\limsup_{n\rightarrow\infty} \bE\left[\left(\frac{J_n^\psi}{\log
n}\right)^p\right]<\infty.\vspace{-0.1cm}$$
\end{lemm}

\begin{pf}
The case ${\rm Var}(-\log Y_1)<\infty$, and implicitly also $\bE[-\log Y_1]<\infty$, has been shown in
the proof of \cite[Proposition 8]{Gne-06}. Now let $\bE[-\log Y_1]\!=\!\infty$. Define $J^\prime_n=\#\{k\ge 1\colon 
G_k\ge 1/n\}=\#\{k\ge 1\colon  \sum_{i=1}^k (-\log Y_i)\leq \log n\}$.
By the Renewal Theorem \cite[Theorem 4.1, Chapter 3]{Db},\linebreak
$J^\prime_n/\log n\rightarrow 0$ a.s.\ when $\bE[-\log
Y_1]\!=\!\infty$. Let $I_{1,n}\!<\!\cdots\!<\! I_{n,n}$ be the order statistics
of $I_1,\ldots, I_n$. Define $\beta_n$ by $
I_{\beta_n,n}\!<\!1/n\!<\!I_{\beta_{n+1},n}$. According to Gnedin's
discussion, $J^\prime_n$ and $\beta_n$ are independent, $\beta_n$ is
binomial$(n, 1/n)$ and $J_n^\psi\!\leq\! J^\prime_n\!+\!\beta_n$. By Markov's inequality,
we have for all $\epsilon>0,$\vspace{-0.2cm}
\begin{equation}
   \bP(\beta_n>\epsilon\log n)=\bP\left({\rm
   e}^{2\beta_n/\epsilon}> n^2\right)\leq\frac{\bE\left[{\rm
   e}^{2\beta_n/\epsilon}\right]}{n^2}=\frac{1}{n^2}\left(1+\frac{{\rm
   e}^{2/\epsilon}-1}{n}\right)^n.\vspace{-0.2cm}\nonumber
\end{equation}
Hence,
$\sum_{n\ge 1}\bP(\beta_n>\epsilon\log n)<\infty$.
The Borel-Cantelli Lemma now yields
$\lim_{n\rightarrow\infty}\beta_n/\log n=0$ a.s. This gives us
$\limsup_{n\rightarrow\infty} J_n^\psi/\log n=0$ when $\bE[-\log
Y_1]=\infty$. Finally, for $p\geq 1$,
$$\bE\left[\left(\frac{J_n^\psi}{\log
n}\right)^p\right]\leq \bE\left[\left(\frac{J^\prime_n+\beta_n}{\log
n}\right)^p\right]\leq 2^{p-1}\left(\bE\left[\left(\frac{J^\prime_n}{\log
n}\right)^p\right]+\bE\left[\left(\frac{\beta_n}{\log
n}\right)^p\right]\right).$$ The first term is bounded (Lemma
\ref{renew}), the second tends to 0 ($\beta_n$ have bounded moments).
\end{pf}

\subsection{Special branch points and their asymptotics}\label{Gpb}

We consider the setting of Theorem \ref{thm3}, where $m=\inf\{n\ge 1\colon \nu_j=\nu_n\mbox{ for all $j\ge n$}\}<\infty$. In this setting, the selection probabilities of Section \ref{sample} for $k\geq m+1$ become
$$P^{\rm old}_k(\fs, s_i)= s_i\qquad\mbox{and}\quad P^{\rm new}_k(\fs, s_i, s_j)=s_j.$$
It is now easy to see that the sampling procedure in $(\cT,\mu)$ can be simplified in this setting so as to combine for each $k\ge m$ the steps until $\#B^\prime<m$ into a single selection according to $\mu$. 

\begin{proc}\label{proc3}\rm Use the steps of Procedure \ref{proc2}, but instead of
  steps ($k$,$[k-1]$)1.-2., use the following steps for $k\ge m$:
\begin{enumerate}\item[1$^\prime$.] Given $(\cT;\Sigma_i\in[k-1])$, sample $\Sigma_{k}^*\sim\mu$ independently. 
  \item[2$^\prime$.] We consider the spine $[[\rho,\Sigma_{k}^*[[=\{\Sigma_{k}^{*{\alpha\rightarrow 0}}(t),t\ge 0\}$ 
    and set 
    $$\cT_{k,B^\prime}=\cT^{\alpha\rightarrow 0}_{(\Sigma_{k}^*)}(\tau_k^*),\mbox{ where $\tau_k^*=\inf\{t\ge 0\colon\#\cL_{k-1}(\cT_{(\Sigma_{k}^*)}^{\alpha\rightarrow 0}(t))<m\}$, $B^\prime=\cL_{k-1}(\cT^{\alpha\rightarrow 0}_{(\Sigma_{k}^*)}(\tau_k^*))$,}$$
    and $\cL_{k-1}(\cS)=\{i\in[k-1]\colon\Sigma_i\in\cS\}$ is the set of labels in $\cS\subseteq\cT$. 
%
\end{enumerate}
\end{proc}
Theorem \ref{thm3} describes the convergence of unlabelled trees. In fact, more is true and it will be instructive to study
approximations of the spines $[[\rho,\Sigma_j[[$, $j\ge 1$, in $(\cT;\Sigma_i,i\in\bN)$ by discrete 
spines $\{B\in T_n\colon j\in B\}$, $n\ge j\ge 1$. In the proof of Theorem \ref{thm3} we will need to control these uniformly in
$j\ge 1$. In the exchangeable case, these spines can be regarded as independent uniform samples from a
strongly sampling consistent regenerative interval partition \cite{HMPW}. In the RE case here, the analogous
partitions will no longer be regenerative (except for $j=1$, and for $j=2$ if $m=2$) and the sampling is not independent uniform. However, 
both features are still present on parts of the spine and we will cut the spines at certain \em special branch points\em. 

Fix $j\ge 1$. A branch point $v\in[[\rho,\Sigma_j[[$ is called \em special in $(\cT;\Sigma_i,i\in\bN)$ for $[[\rho,\Sigma_j[[$ \em if some or all of the $m$ smallest labels $\cL(\cT^v)$ in the bush $\cT^v$ above $v$ are \em not \em included in the subtree $\cT_{(\Sigma_j)}({\rm d}(\rho,v))$ above $v$ containing $\Sigma_j$. Note that a branch point $v$ is special iff at $v$ the
$m$ smallest labels split or $j$ splits from the $m$ smallest labels. In particular, a branch point that is special for $[[\rho,\Sigma_j[[$ and an element of $[[\rho,\Sigma_{j^\prime}[[$ for some $j^\prime<j$ may not be special for $[[\rho,\Sigma_{j^\prime}[[$. For the analogous notion in $(\cT;\Sigma_i,i\in[n])$, for $n\ge j$, we write 
$$N_n^{(j)}=\#\{v\in[[\rho,\Sigma_j[[\colon \mbox{$v$ is a special branchpoint for $[[\rho,\Sigma_j[[$ in $(\cT;\Sigma_i,i\in[n])$}\}$$
for the number of special branch points, and $\tau_n^{(j)}=\inf\{t\ge 0\colon \#\cL_n(\cT_{(\Sigma_j)}^{\alpha\rightarrow 0}(t))<m\}$ for the time
when the label set first has fewer than $m$ elements. The significance of this time is that up to this time, all branch points that are special in $(\cT;\Sigma_i,i\in\bN)$ will also be special in $(\cT;\Sigma_i,i\in[n])$, but this fails afterwards.
We introduce $V_n^{(j)}=\inf\{t\ge 0\colon \Sigma_n\not\in\cT_{(\Sigma_j)}^{\alpha\rightarrow 0}(t)\}$, the
time when $\Sigma_n$ leaves the spine $[[\rho,\Sigma_j[[$.

%

\begin{prop}\label{sbp} Let $(\nu_j,j\ge 1)$ and $\nu$ be as in Theorem \ref{thm2} and $(\cT;\Sigma_i,i\in\bN)$ a sampling
  according to Procedure \ref{proc2}. Suppose furthermore that there is 
  $m\ge 1$ with $\nu_j=\nu_m$, $j\ge m$, and that $\nu_m(s_1\leq 1-\epsilon)=\epsilon^{-\alpha}\ell(1/\epsilon)$, which is equivalent to (\ref{regvar}) for $\nu$ as in (\ref{thm6form}). Then, 
\begin{enumerate}
\item[\rm (i)] for all $j\ge 1$, we have
$N_n^{(j)}/(n^\alpha\ell(n))\xrightarrow[n\rightarrow\infty]{a.s.}
0;$

\item[\rm (ii)] for every $p\geq 1$, we have $\displaystyle\limsup_{n\rightarrow\infty}\bE\left[\left(N_n^{(n)}/\log
n\right)^p\right]<\infty;$

\item[\rm (iii)]
 for every $p\geq 1$, there exists a constant $C_p^{\rm spec}$ such that for all $1\leq j\leq n$ and $x>0$
$$\bP\left(N_n^{(j)}>xn^\alpha\ell(n)\right)<\frac{C_p^{\rm spec}}{x^p n^{\alpha p-1}}.
$$
\end{enumerate}
\end{prop}

\begin{pf}
(i) Let us consider $N_n^{(1)}$ first. We will study the asymptotics by relating to the setting of Lemma \ref{Gcon}.
  Recall that 
    $X_{(\Sigma_1)}^{\alpha\rightarrow 0}(V_i^{(1)})$, $i\geq 2$, are the residual masses of the subtrees containing $\Sigma_1$ when $\Sigma_i$ has
    left the spine $[[\rho, \Sigma_1[[$. Let $Y_k$, $k\geq 1$, be independent copies of $X_{(\Sigma_1)}^{\alpha\rightarrow 0}(\tau_m^{(1)})$, the 
    residual mass of the subtree containing $\Sigma_1$ at the branch point separating $[m]$, and $G_k=Y_1\cdots Y_k$. 
    
Consider the filtration 
$\cF_{n}^{(1)}(t)=\sigma\left((S^{\Sigma_1}(s),F_{\Sigma_1}(s),\cB^{\alpha\rightarrow 0}_{(\Sigma_1)}(s),\cL_{n}(\cB^{\alpha\rightarrow 0}_{(\Sigma_1)}(s))),s\le t\right)$, $t\ge 0$,  
of the spinal Poisson point process $(S^{\Sigma_1},F_{\Sigma_1},\cB^{\alpha\rightarrow 0}_{(\Sigma_1)})$ studied in Proposition \ref{sleaf} augmented by
label sets of spinal bushes derived from sampled leaves $\Sigma_1,\ldots,\Sigma_n$. 
Let $H_n^{(1)}=\#\{\tau_i^{(1)},m\le i\le n\}$. Then $H_m^{(1)}=1$ is the initial state, we will also consider  $(\tau_m^{(1)},X_{(\Sigma_1)}^{\alpha\rightarrow 0}(\tau_m^{(1)}),\#\cL_m(\cT_{(\Sigma_1)}^{\alpha\rightarrow 0}(\tau_m^{(1)})))$. Now let $n\ge m+1$ and write $\overline{V}_n^{(j)}=\min\{\tau_n^{(j)},V_n^{(j)}\}$, $n\ge 1$. Conditionally given $\cF_{n-1}^{(1)}(\tau_{n-1}^{(1)})$, in particular $(X_{(\Sigma_1)}(\overline{V}_m^{(1)}),\ldots,X_{(\Sigma_1)}(\overline{V}_{n-1}^{(1)}))$,
$H_{n-1}^{(1)}=k$ and $\#\cL_{n-1}(\cT_{(\Sigma_1)}^{\alpha\rightarrow 0}(\tau_{n-1}^{(1)}))=\ell$, the argument to establish Procedure \ref{proc3} can be used to simplify Procedure \ref{proc2} slightly differently with modified steps 1$^\prime$.-2$^\prime$. combining the steps until $\#B^\prime<m$ or $1\not\in B^\prime$; specifically,  sample a leaf $\Sigma_n^*\sim\mu$, define $V^{(1)}_{n,*}=\inf\{t\ge 0\colon 1\not\in \cL_{n-1}(\cT_{(\Sigma_n^*)}^{\alpha\rightarrow 0}(t))\}$ and  
\begin{itemize}\item[--] if $V_{n,*}^{(1)}\le\tau_{n-1}^{(1)}$, set $\cT_{n-1,B^\prime}=\cT_{(\Sigma_1)}^{\alpha\rightarrow 0}(V_{n,*}^{(1)})$, note $H_n^{(1)}\!\!=\!k$, $\tau_n^{(1)}\!\!=\!\tau_{n-1}^{(1)}$, $\#\cL_n(\cT_{(\Sigma_1)}^{\alpha\rightarrow 0}(\tau_n^{(1)}))\!=\ell$;
  \item[--] if $V_{n,*}^{(1)}>\tau_{n-1}^{(1)}$ and $\ell<m-1$, set $\cT_{n-1,B^\prime}=\cT_{(\Sigma_1)}^{\alpha\rightarrow 0}(\tau_{n-1}^{(1)})$, note $H_n^{(1)}=k$, $\tau_n^{(1)}=\tau_{n-1}^{(1)}$, $\#\cL_n(\cT_{(\Sigma_1)}^{\alpha\rightarrow 0}(\tau_n^{(1)}))=\ell+1$;
  \item[--] if $V_{n,*}^{(1)}>\tau_{n-1}^{(1)}$ and $\ell=m-1$, then sampling of $\Sigma_n$ in the rescaled subtree $\cT_{(\Sigma_1)}^{\alpha\rightarrow 0}(\tau_{n-1}^{(1)})$ is independent of $\cF_{n-1}^{(1)}(\tau_{n-1}^{(1)})$ and by the same procedure as $\Sigma_m$ is sampled in $\cT$, therefore $$X_{(\Sigma_1)}^{\alpha\rightarrow 0}(\overline{V}_n^{(1)})\overset{d}{=}X_{(\Sigma_1)}^{\alpha\rightarrow 0}(\overline{V}_{n-1}^{(1)})Y_{k+1}\overset{d}{=}G_{k+1},$$
note $H_n^{(1)}\!\!=\!k+1$, $\tau_n^{(1)}-\tau_{n-1}^{(1)}\overset{d}{=}\tau_m^{(1)}$ independent of $\cF_{n-1}^{(1)}(\tau_{n-1}^{(1)})$, and $\#\cL_n(\cT_{(\Sigma_1)}^{\alpha\rightarrow 0}(\tau_n^{(1)}))<m$.
\end{itemize}
Independently of $(G_k,k\ge 1)$, consider $(\Psi_k,k\ge 1)\overset{d}{=}(m-\#\cL_{W_{k+1}}(\cT_{(\Sigma_1)}^{\alpha\rightarrow 0}(\tau_{W_{k+1}}^{(1)})),k\ge 1)$,
where $W_k=\inf\{n\ge 1\colon H_n^{(1)}=k\}$. 
As $(G_k,k\ge 1)\overset{d}{=}(X_{(\Sigma_1)}^{\alpha\rightarrow 0}(\overline{V}_{W_k}^{(1)}),k\ge 1)$, it is now straightforward to show that the dynamics of $H_n^{(1)}$ 
and $J_{n-m+1}^{\Psi}$ 
are the same, hence there exists a sequence $(I_i,i\ge 1)$ of independent uniform random variables on $[0,1]$ and an independent
random sequence $\Psi$, each member taking values in $[m]$ such that for all $n\geq m$
\begin{equation}\label{sbp1}
\left(H_m^{(1)},\ldots,H_n^{(1)}\right)\overset{d}{=}\left(J_1^{\Psi},\ldots,J_{n-m+1}^{\Psi}\right).
\end{equation} 
%
%
Now note that $n\in\bN$ with $W_k<n<W_{k+1}$ can only yield a new special branch point if  $V_{n,*}^{(1)}>\tau_{n-1}^{(1)}$, i.e. in the middle case of the procedure above, but after at most $m-1$ such steps, the
third case will apply and $H_n^{(1)}$ will increase. Therefore, 
\begin{equation}\label{sbp11} N_n^{(1)}\leq
m\#H_n^{(1)}.\end{equation} 
Lemma \ref{Gcon} ensures $H_n^{(1)}/\log n\rightarrow 1/\bE[-\log Y_1]$,
therefore $N_n^{(1)}/(n^\alpha\ell(n))\rightarrow 0$ a.s. as
$n\rightarrow \infty$.

The same argument, with $\Sigma_1$ replaced by $\Sigma^*\sim\mu$, yields $N_n^*/(n^\alpha\ell(n))\rightarrow 0$ a.s. as $n\rightarrow\infty$. For $j\ge 2$,
consider times $\chi_i^{(j)}=\inf\{t\ge 0\colon \#\cL_j(\cT^{\alpha\rightarrow 0}_{(\Sigma_j)}(t))<m-i\}$, $0\le i<m$, when $j$ changes rank below $m$ in the label set, s.th.\ $\chi_0^{(j)}\!\!=\!\tau_j^{(j)}$ and $\chi_{m-1}^{(j)}\!\!=\!\infty$. As $\#\cL_n(\cT_{(\Sigma_j)}^{\alpha\rightarrow 0}(\chi_i^{(j)}))\!\le\! n\!-\!j\!+\!m$, the number of special branch points between $\Sigma_j^{\alpha\rightarrow 0}(\chi_i^{(j)})$ and $\Sigma_{\ell_i}$,
where $\ell_i=\min\cL_j(\cT^{\alpha\rightarrow 0}_{(\Sigma_j)}(\chi_i^{(j)})$, will be no larger than $\widetilde{N}_{n-j+m}^{(1),i}$ where $(\widetilde{N}_k^{(1),i}, k\geq 1)$ are independent copies of $(N_k^{(1)},k\geq 1)$. Then
\begin{equation}\label{sbp2}
N_n^{(j)}\le \widetilde{N}_n^*+\sum_{i=1}^{m-1}\widetilde{N}_{n-j+m}^{(1),i},
\end{equation}
where $\widetilde{N}_n^*$ is the number of special branchpoints on $[[\rho,\Sigma_j^*[[$, so that
$\widetilde{N}^*\overset{d}{=}N^*$.  
%
Hence the convergence for $N_n^{(j)}/(n^\alpha \ell(n))$ follows from previous cases of $N^{(1)}$ and $N^*$.

(ii) To study $N^{(n)}_n$, we will identify new families $(G_k,k\ge 1)$ and $\psi$ different from the ones in (i) and again apply Lemma \ref{Gcon}. Let $b^{(j)}_1=\Sigma^{\alpha\rightarrow 0}_j(\chi^{(j)}_1)$ be the first special branch point in the spine
$[[\rho,\Sigma_j]]$. By Procedure \ref{proc3}, $X_{(\Sigma_j)}^{\alpha\rightarrow 0}(\chi_1^{(j)})=X_{(\Sigma_j^*)}^{\alpha\rightarrow 0}(\chi_1^{(j)})$
 for all $j\geq m+1$. Also, note that $\chi_1^{(j)}$ is determined by
 $(\cT;\Sigma_i,i\in[m];\Sigma_j^*)$. As $\Sigma_j^*$ is sampled
 according to $\mu$ in $\cT,$ we have
\begin{equation}\label{sbp22}X_{(\Sigma_j)}^{\alpha\rightarrow 0}(\chi_1^{(j)})=X_{(\Sigma_j^*)}^{\alpha\rightarrow 0}(\chi_1^{(j)})\overset{d}{=}X_{(\Sigma_{m+1}^*)}^{\alpha\rightarrow 0}(\chi_1^{(m+1)})=X_{(\Sigma_{m+1})}^{\alpha\rightarrow 0}(\chi^{(m+1)}_1).	
\end{equation}
Let $Y_k$, $k\geq 1$, be independent copies of
$X_{\Sigma_{m+1}}^{\alpha\rightarrow 0}(\chi^{(m+1)}_1)$ and consider a constrained painbox associated with $G_k=Y_1\cdots
Y_k$, $k\geq 1$, also $\psi_k=1$, $k\ge 1$. We claim that for all $n\geq m+1,$ and every
$x>0$,
\begin{equation}\label{sbp3}\bP(N_n^{(n)}-m+1>x)\leq
\bP(J_{n-m}^\psi>x).\end{equation}
This formula holds for $n=m+1$ as $N_{m+1}^{(m+1)}-m+1\leq J_1^\psi=1$.
Suppose (\ref{sbp3}) holds for all $n\leq j-1$. For $n=j$, the first
special branch point $b_1^{(j)}$ on the spine $[[\rho, \Sigma_j]]$ is
located on the spine $[[\rho, b_1^{(1)}]]$. For
$i=m+1,\ldots, j-1$, let $\cT_{(\Sigma_i)}^{\alpha\rightarrow 0}(V_i^{(1)}\wedge\chi_1^{(1)})$ be the
spinal subtree of $\cT$ containing $\Sigma_i$ rooted on a branch
point on the spine $[[\rho, b_1^{(1)}]]$, possibly at $b_1^{(1)}$ itself. 
By Procedure \ref{proc3}, $\Sigma_i^*\in
\cT_{(\Sigma_i)}^{\alpha\rightarrow 0}(V_i^{(1)}\wedge\chi_1^{(1)})$.
We can express the number $M_1^{(j)}$ of leaves in $\{\Sigma_{m+1},\ldots, \Sigma_{j-1}\}$ belonging to
the subtree containing $\Sigma_j$ above branch point $b_1^{(j)}$ as
$$M_1^{(j)}\!\!=\!\#\{i\in\{m+1,\ldots,j-1\}\colon \Sigma_i\in\cT_{(\Sigma_j)}^{\alpha\rightarrow 0}(\chi_1^{(j)})\}=\#\{i\in\{m+1,\ldots,j-1\}\colon \Sigma_i^*\in\cT_{(\Sigma_j)}^{\alpha\rightarrow 0}(\chi_1^{(j)})\}.$$ 
As $\Sigma_{m+1}^*,\ldots, \Sigma_{j-1}^*$ are sampled according to $\mu$ and
$X_{(\Sigma_j)}^{\alpha\rightarrow 0}(\chi_1^{(j)})\overset{d}{=}X_{(\Sigma_{m+1})}^{\alpha\rightarrow 0}(\chi^{(m+1)}_1)\overset{d}{=}Y_1$,
by (\ref{sbp22}),
\begin{eqnarray*}\bP\left(M_1^{(j)}=k\right)&=&\bE\left[{j-m-1\choose
k}\left(X_{(\Sigma_j)}^{\alpha\rightarrow 0}(\chi_1^{(j)})\right)^k\left(1-X_{(\Sigma_j)}^{\alpha\rightarrow 0}(\chi_1^{(j)})\right)^{j-m-k-1}\right]\\
&=& \bE\left[{j-m-1\choose
k}Y_1^k(1-Y_1)^{j-m-k-1}\right]=\bP\left(\overline{M}_{j-m}^\psi=k\right)
\end{eqnarray*} for all $0\leq k\leq j-m-1$, where
$\overline{M}_{j-m}^\psi$ is the number of $\overline{I}_1^\psi,\ldots,
\overline{I}_{j-m}^\psi$ hitting the interval $(0,G_1)$. 

Let $N_j^{(j)}(\chi_1^{(j)},\infty)=N_j^{(j)}-1$ be the number of special branch points
in $]]b_1^{(j)}, \Sigma_j]]$, and $J_{j-m}^\psi(0,Y_1)=J_{j-m}^\psi-1$. 
Given $M_1^{(j)}=k$, we have $\#\cL_j(\cT_{(\Sigma_j)}^{\alpha\rightarrow 0}(\chi_1^{(j)}))\le k+m\le j-1$. Hence, applying the induction hypothesis to the rescaled $(\cT_{(\Sigma_j)}^{\alpha\rightarrow 0}(\chi_1^{(j)});\Sigma_i,i\in\cL_j(\cT_{(\Sigma_j)}^{\alpha\rightarrow 0}(\chi_1^{(j)})))$
\begin{eqnarray*}
&&\hspace*{-1cm}\bP\left(\left.N_j^{(j)}(\chi_1^{(j)},\infty)-m+1>x\right|M_1^{(j)}=k,
\right)\\
&\leq&\bP\left(N_{k+m}^{(k+m)}-m+1>x\right)\leq \bP\left(J_k^\psi>x\right)=
\bP\left(J_{j-m}^\psi(0,Y_1)>x\left|\overline{M}_{j-m}^\psi=k\right.\right),
\end{eqnarray*}
and then
\begin{eqnarray*}
\bP\left(N_j^{(j)}-m+1>x\right)&=&\bE\left[\bP\left(\left.N_j^{(j)}(\chi_1^{(j)},\infty)-m+1>x-1\right|M_1^{(j)}\right)\right]\\
&\leq&\bE\left[\bP\left(J_{j-m}^\psi(0,Y_1)>x-1\left|\overline{M}_{j-m}^\psi\right.\right)\right]=\bP\left(J_{j-m}^\psi>x\right).
\end{eqnarray*}

Now (\ref{sbp3}) is clear and we deduce that 
$\bE\left[(N_n^{(n)}-m+1)^p\right]\leq
\bE\left[\left(J_{n-m}^\psi\right)^p\right]$ for every $p\geq
1$. 
The result in (ii) now follows from Lemma \ref{Gcon}.

(iii)  Formula (\ref{sbp2}) implies that for every $p\ge 1$ and $x>0$ and $z_n=xn^\alpha\ell(n)$
\begin{eqnarray*}
\bP\left(N_n^{(j)}>z_n\right)&\leq&
\bP\left(\widetilde{N}^*_n>z_n/m\right)+\sum_{i=1}^{m-1}\bP\left(\widetilde{N}_{n-j+m}^{(1),i}>z_n/m\right)\\
&\leq&\frac{\bE\left[(N_n^*)^p\right]}{z_n^p/m^p}+(m-1)
\frac{\bE\left[(N_{n-j+m}^{(1)})^p\right]}{z_n^p/m^p}\le\frac{C_p(\log n)^p}{z_n^p}.
\end{eqnarray*}
The last line is obtained by Markov's inequality. 
Formula (\ref{sbp11}) together with Lemma \ref{Gcon} gives the upper
bounds. 
The result in (iii) follows.
\end{pf}

Procedure \ref{proc3} and the notion of special branch points are also useful to show that the sampling uses the whole CRT $(\cT,\mu)$ and does not leave any
subtrees of positive mass unlabelled. One way of making this precise is to say that the reduced trees converge to the CRT:

\begin{prop}\label{fill} In the setting of Procedure \ref{proc3}, we have
  $$R(\cT;\Sigma_i,i\in[k])\rightarrow \cT\qquad\mbox{a.s. in the Gromov-Hausdorff sense as $k\rightarrow\infty$.}$$ 
\end{prop} 
\begin{pf} Let $\epsilon>0$. Consider $[[\rho,\Sigma_1[[$ and the associated spinal mass partition \cite{HPW}. 
  Here we denote by $\nu^{\rm sp}_\epsilon$ the distribution on $S^\downarrow$ of the masses of spinal subtrees that are greater than
  $\epsilon$. Let $\sigma^{(1)}_\epsilon=\inf\{t\ge 0\colon \mu(\cT^{\alpha\rightarrow 0}_{(\Sigma_1)}(t))<\epsilon\}$. Note that  
  $W_1:=\inf\{n\ge 1\colon \tau_n^{(1)}\ge\sigma^{(1)}_\epsilon\}<\infty$ a.s., by the 
  previous proof.  By Procedure \ref{proc3}, leaves
  $\Sigma_n^*$ and $\Sigma_n$ are in the same subtree of $[[\rho,\Sigma_1^{\alpha\rightarrow 0}(\sigma_\epsilon)]]$ for each
  $n>W_1$, in particular each subtree of mass greater than $\epsilon$ is selected with an asymptotic frequency greater than
  $\epsilon$. Inductively, we use Corollary \ref{cor22} and leaves selected according to Procedure \ref{proc3} to further
  split according to scaled $\nu^{\rm sp}_\epsilon$ each subtree of mass greater than $\epsilon$. 
 
 After a finite number of steps, all subtrees have mass less than $\epsilon$, e.g. because a homogeneous mass 
  fragmentation process $(F_t,t\ge 0)$ in $S^\downarrow$ with finite dislocation measure $\nu_\epsilon^{\rm sp}$ satisfies $F_t\rightarrow 0$ as $t\rightarrow\infty$, see e.g. \cite[Equation (4)]{Ber-asymp}, and so only has finitely
  many splits before $|F_1(t)|<\epsilon$.  
\end{pf}

Using arguments of \cite[Corollary 23]{PW2}, we can also show joint a.s.\ convergence in the Gromov-Prohorov sense of weighted
trees $(R(\cT;\Sigma_i,i\in[n]),n^{-1}\sum_{i=1}^n\delta_{\Sigma_i})\rightarrow(\cT,\mu)$.

\subsection{Convergence of reduced trees and large deviation estimates for spines} 

By Corollary \ref{cor24}, reduced trees $R(\cT;\Sigma_i,i\in[k])$ of self-similar CRTs with labelled leaves sampled according to Procedure \ref{proc3}, can be assigned subtree masses on edges (parts of spines) in terms of
Poisson point processes and associated spinal subordinators, and away from existing leaves, sampling of new leaves is according to subtree masses. To study the asymptotics of the number of spinal branchpoints, we will need the following refinement of results in \cite{GPYalpha06,HMPW}.

%

\begin{lemm}\label{pjs} Let $\xi=(\xi_t,t\geq 0)$ be a
pure jump subordinator with L\'{e}vy measure $\Lambda$ satisfying
$\Lambda([x,\infty))=x^{-\alpha}\ell_\Lambda(1/x)$, $x\downarrow 0$, for some $\alpha\in(0,1)$. Let $(\epsilon,\tau,\tau^\prime)$ be any random variables on $[0,\infty)^2\times[0,\infty]$ with $\tau\le\tau^\prime$. Let $(V_i,i\ge 1)$ be any random variables
conditionally independent given $(\xi,\epsilon,\tau)$ with 
$$\bP\left(\left.V_i\le\tau\right|\xi,\epsilon,\tau\right)=1-e^{-\epsilon}\quad\mbox{and}\quad\bP\left(\left.V_i>\tau+v\right|\xi,\epsilon,\tau\right)=e^{-\epsilon-\xi_v},\quad v\ge 0,$$ 
and $K_n(\epsilon,\tau,\tau^\prime)=\#\{V_i\colon 1\leq i\leq n,\tau<V_i\leq\tau^\prime\}$. Then
$$\lim_{n\rightarrow\infty}\frac{K_n(\epsilon,\tau,\tau^\prime)}{n^{\alpha}\ell_\Lambda(n)\Gamma(1-\alpha)}=\int_0^{\tau^\prime-\tau}\exp(-\alpha(\epsilon+\xi_v))dv\quad\mbox{a.s. as $n\rightarrow\infty$.}$$
If furthermore $\Lambda([xy,\infty))\le C_\Lambda y^{-\varrho}\Lambda([x,\infty))$ for all $y\ge 1$ and $0<x\le 1$, and some $\varrho>0$, 
then there is a constant $C_p$ for all $p>1/\alpha$, such that for all $x\ge 1$, $n\ge 1$ and all $(\epsilon,\tau,\tau^\prime)$ as above, but with the additional property that $\tau^\prime=\tau+\tau^{\prime\prime}$ for a stopping time $\tau^{\prime\prime}$ for a filtration in which $\xi$ is a subordinator,
\begin{equation}\label{tcp0}\bP\left(\frac{K_n(\epsilon,\tau,\tau^\prime)}{n^\alpha\ell_\Lambda(n)\Gamma(1-\alpha)}>(1+x)Y(\epsilon,\tau,\tau^\prime)\right)\le\frac{C_p}{x^pn^{\alpha p-1}},
\end{equation}
where $\displaystyle Y(\epsilon,\tau,\tau^\prime)=1+(1+A_\alpha)C_\Lambda\sum_{j=0}^{[\tau^\prime-\tau]}\exp(-\varrho(\epsilon+\xi_j))$ with $\displaystyle A_\alpha=2\sum_{j\ge 1}\frac{(j+1)^{\sqrt{\alpha}}}{j(j+1)}.$
\end{lemm}

This lemma is an extension of \cite[Lemmas 8 and 12]{HMPW}, which we recover as the special case $\tau=\epsilon=0$ and/or $\tau^\prime=\infty$. The proof is also essentially the same, but since this result is more general, we reproduce the proof
rewritten in the present generality in the appendix.

\pagebreak[2]

\begin{prop}\label{fdm}
Let $\nu_1,\ldots, \nu_m$ be conservative with
$\nu(s_1\leq 1-\epsilon)=\epsilon^{-\alpha}\ell(1/\epsilon)$, where
$\nu$ is as in Theorem \ref{thm2} with $\nu_j=\nu_m$, $j\geq m$. Let
$R(\cT, \Sigma_1,\ldots, \Sigma_n)$ be an $\bR$-tree sampled from an
$\alpha$-self-similar CRT $(\cT, \mu)$ with dislocation
measure $\nu$ by Procedure \ref{proc3}. Let $(T_n)_{n\geq 1}$ be the
associated labelled discrete RE
Markov branching trees with unit edge lengths. Then
$$
\frac{R(T_n,[k])}{n^{\alpha}\ell(n)\Gamma(1-\alpha)}\xrightarrow[n\rightarrow\infty]{a.s.}
 R(\cT,\Sigma_1,\ldots,\Sigma_k)\qquad\mbox{in the sense that all edge lengths converge.}$$
In particular, the delabelled trees $(R(T_n,[k]))^\circ$, $n\ge k$, converge in the Gromov-Hausdorff sense.
\end{prop}

\begin{pf} 
Consider $k=1$ and denote by $D_n^{(1)}$ the length of $R(T_n,
\{1\})$. 

If $\nu_1=\cdots=\nu_{m-1}=0$, then
$\Sigma_1,\ldots,\Sigma_m$ are always in the same subtree in $\cT$,
then $\tau_1^{(1)}=\cdots=\tau^{(1)}_{m}=\infty$. Conditionally on the
subordinator  $\xi^{\Sigma_1}$ associated with leaf $\Sigma_1$, cf. Proposition \ref{sleaf}(ii), the
leaves $\Sigma_{m+1},\ldots, \Sigma_n$ are sampled according to
$\mu$ along the spine $[[\rho, \Sigma_1[[$. Using Proposition \ref{sleaf}(ii), we see that the hypotheses of the first part of Lemma \ref{pjs} are satisfied, and the convergence result then follows. Specifically, it is easy to see that by (\ref{regvar}), as $x\downarrow 0$,
$$\Lambda_1([x,\infty))\sim\int_{\{1/2<s_1\le e^{-x}\}}P_0(\fs,s_1)\nu(d\fs)\sim\nu(s_1\le e^{-x})\sim x^\alpha\ell(1/x),$$
since by (\ref{finite}),
$\displaystyle\int_{\{1/2<s_1\le e^{-x}\}}(1-P_0(\fs,s_1))\nu(d\fs)\le 2\int_{S^\downarrow}\sum_{i\ge 2}(1-s_i)P_0(\fs,s_i)\nu(d\fs)<\infty$.

Now suppose that at least one of
$\nu_1,\ldots,\nu_{m-1}$ is non-zero. 
By Procedure \ref{proc3}, each $\Sigma_i$ is either placed in the same subtree of 
$[[\rho,\Sigma_1[[$ as $\Sigma_i^*\sim\mu$ or contributes a special branch point. Now 
\begin{equation}\label{fdm2}
D^{(1)}_n= \# \left\{V^{(1)}_i, 1\leq i\leq n\right\}\le 1+N_n^{(1)}+\#\left\{V^{(1)}_{i,*},2\leq i\leq n\right\},
\end{equation}
with $V_{i,*}^{(1)}=\inf\{t\ge 0\colon 1\not\in\cL_{n-1}(\cT_{(\Sigma_n^*)}^{\alpha\rightarrow 0}(t))\}$, where Lemma \ref{pjs} yields the asymptotics of $K_{n-1}^{(1)}(0,0,\infty)=\#\{V^{(1)}_{i,*},2\leq
i\leq n\}$.
Together with the asymptotics of $N_n^{(1)}$ obtained in Proposition \ref{sbp}, this yields
\begin{equation}\label{fdm3}\limsup_{n\rightarrow\infty}\frac{D^{(1)}_n}{n^\alpha\ell(n)\Gamma(1-\alpha)}\leq\int_0^\infty \exp(-\alpha\xi^{\Sigma_1}_t)dt\qquad\mbox{a.s.}\end{equation}
On the other hand, no special branch points are created for $n\ge l+1\ge m+2$ below $\tau_l^{(1)}$, so
$$D_n^{(1)}\geq 
\#\left\{V^{(1)}_i\colon  0<V^{(1)}_i\leq\tau^{(1)}_l,
l+1\leq i\leq n\right\}
=\#\left\{V^{(1)}_{i,*}\colon  0<V^{(1)}_{i,*}\leq\tau^{(1)}_l, l+1\leq i\leq
n\right\}.$$
At least one of $\nu_j\neq 0$, $j\le m-1$, so
$\tau_m^{(1)}<\infty$. By the proof of Proposition \ref{sbp}, $\tau_l^{(1)}\rightarrow\infty$, 
so
$$
\liminf_{n\rightarrow\infty}\frac{D_n^{(1)}}{n^\alpha\ell(n)\Gamma(1-\alpha)}\geq\sup_{l\geq m+1}
\liminf_{n\rightarrow\infty}\frac{\#\{V_{i,*}^{(1)}\colon 
V^{(1)}_{i,*}\leq\tau^{(1)}_l, l+1\leq i\leq
n\}}{n^\alpha\ell(n)\Gamma(1-\alpha)}
=\int_0^\infty\!\!\!\! \exp(-\alpha\xi^{\Sigma_1}_t) dt.
$$
Combining this with (\ref{fdm3}), 
the convergence for $D^{(1)}_n$
follows and establishes the result for $k=1$.\pagebreak[2]

\noindent Next, consider $k\ge 2$ assuming the result for $1,\ldots,k-1$. For the branch point $v_{k}$
adjacent to $\rho$ in $R(\cT, \Sigma_1,\ldots, \Sigma_{k})$, set
$D^{[k]}={\rm d}(\rho, v_{k})$, with time
$\zeta_{v_{k}}$ given by 
$$ D^{[k]}
=\int_0^{\zeta_{v_{k}}} \exp(-\alpha\xi_t^{\Sigma_1}) dt.$$
Let $D_n^{[k]}$ be the height of the branch point adjacent to the
root in $R(T_n,[k])$, then $D_n^{[k]}-1$ is the number of
distinct branch points of $R(\cT, \Sigma_1,\ldots, \Sigma_n)$
belonging to $[[\rho,v_{k}[[$, i.e.
$$D^{[k]}_n=1+\#\{V^{(1)}_i\colon  0<V_i^{(1)}<\zeta_{v_{k}},k+1\leq i\leq n\}.$$ 
If $2\leq k\leq m$, then $1\leq D^{[k]}_m\leq  m-1$ and, by the same argument as for $k=1$,
\begin{equation}\label{fdmkm}K_{n-m}^{(k)}(0,0,\zeta_{v_{k}})=\#\{V^{(1)}_{i,*}\colon  V^{(1)}_{i,*}<\zeta_{v_{k}}, m+1\leq i\leq n\}\leq
D^{[k]}_n\leq m+K_{n-m}^{(k)}(0,0,\zeta_{v_{k}}). 
\end{equation} 
If $k\geq m+1$, then $D^{[k]}_n=1+\#\{V^{(1)}_{i,*}\colon  V^{(1)}_{i,*}<\zeta_{v_{k}}, k+1\leq i\leq n\}$. In all cases, by Lemma \ref{pjs}
$$
\frac{D^{[k]}_n}{n^{\alpha}\ell(n)\Gamma(1-\alpha)}\xrightarrow[n\rightarrow
\infty]{a.s.}\int_0^{\zeta_{v_{k}}}
\exp(-\alpha\xi_s^{\Sigma_1}) ds=D^{[k]}.
$$
So the renormalized length of the root edge of $R(T_n,[k])$
converges as required. 

Now argue conditionally given that $[k]$ is first separated into $\Pi^{[k]}=(\pi_1,\ldots,\pi_r)$. For all $n\ge k+1$ and $1\le j\le r$, denote by $B_j(n)=\cL_n(\cT^{[k]}_j)\supset\pi_j$ the $j$th block of the partition at $v_{k}$ in $(\cT;\Sigma_i,i\in[n])$, and by $T^{[k]}_{n,j}$ the corresponding subtree of $T_n$. 
By Lemma \ref{lm23}, Procedure \ref{proc3} and the Strong Law of Large Numbers, 
$$\frac{\#B_j(n)}{n}\xrightarrow[n\rightarrow\infty]{a.s.}
\mu(\cT_j^{[k]}),\qquad 1\le j\le r,$$ 
and the Induction Hypothesis yields convergence of the remaining edge lengths, for $1\leq j\leq r$ 
\begin{eqnarray}\label{fdmp5}
\frac{R(T_{n,j}^{[k]},\pi_j)}{n^\alpha\ell(n)\Gamma(1-\alpha)}&=&
\frac{(\#B_j(n))^\alpha\ell(\#B_j(n))}{n^\alpha\ell(n)}\frac{R(T_{n,j}^{[k]},\pi_j)}{(\#B_j(n))^\alpha\ell(\#B_j(n))\Gamma(1-\alpha)}\nonumber\\
&\xrightarrow[n\rightarrow\infty]{a.s.}&(\mu(\cT_j^{[k]}))^\alpha R(\cT^{[k]}_j;\Sigma_i,i\in\pi_j),\nonumber
\end{eqnarray}
in the sense that all edge lengths converge, which implies Gromov-Hausdorff convergence.
\end{pf}

While the arguments of the analogous but much more specific \cite[Proposition 22]{PW2} do not apply here in cases where the
densities $f_k=d\nu_k/d\nu$ are degenerate, we can now deduce from our Proposition \ref{fill} that in the setting of 
Proposition \textrm{\ref{fdm}} here, delabelled trees converge a.s. when taking double limits
\begin{equation}\label{doublelim}\lim_{k\rightarrow
\infty}\lim_{n\rightarrow
\infty}\frac{\left(R(T_n,[k])\right)^\circ}{n^{\alpha}\ell(n)\Gamma(1-\alpha)}=\cT\qquad
\mbox{in the Gromov-Hausdorff sense a.s.}
\end{equation} 
Theorem \ref{thm3}, instead of restricting to $[k]$, then letting $n\rightarrow\infty$ and then $k\rightarrow\infty$, considers
$n\rightarrow\infty$ directly, at the cost of weakening the mode of convergence to convergence in probability. To prepare the proof of Theorem \ref{thm3}, we study the spines $[[\rho,\Sigma_j[[$, $j\ge 1$.

We denote by $\Lambda_1$ and $\Lambda^*$ the L\'{e}vy measures of the subordinators 
$\xi^{\Sigma_1}$ and $\xi^{\Sigma^*}$ generated, respectively, by the first sampled leaf $\Sigma_1$ and by a leaf $\Sigma^*$ sampled according 
to $\mu$. For $k\ge 1$ and $n\ge k$, denote by $D^{(k)}_n$ the length of $R(T_n,\{k\})$.

\begin{lemm}\label{Dj} For all $p\ge 0$, there is a constant $C_p^\prime>0$ such that for all $k\ge 1$, $n\ge k$ and $x\ge 1$
$$\bP\left(D_n^{(k)}>2(1+x)(2+Z_k)\max\{\overline{\Lambda}_1(n^{-1}),\overline{\Lambda}^*(n^{-1})\}\right)\leq
\frac{C^\prime_p}{x^p n^{\alpha p-1}},$$
where $\displaystyle Z_k=m+(1+A_\alpha)\max\{C_{\Lambda_1},C_{\Lambda^*}\}\left(m+\sum_{i=0}^{\infty} \left(X_{(\Sigma_k)}^{\alpha\rightarrow 0}(i)\right)^\varrho\right)$
has all moments finite.
\end{lemm}

\begin{pf} For $k=1$, we use (\ref{fdm2}) to write $D^{(1)}_n\leq 2(D^{(1)}_n-1)\leq 2N_n^{(1)}+2K_{n-1}^{(1)}(0,0,\infty)$ and
  deduce from Proposition \ref{sbp} and Lemma \ref{pjs} that for all $p\ge 0$ and all $n\ge 1$, $x\ge 1$, \vspace{-0.1cm}
\begin{eqnarray}
&&\hspace*{-0.5cm}\bP\left(D_n^{(1)}>2(1+x)(2+Z_1)\overline{\Lambda}_1(n^{-1})\right)\nonumber\\[-0.1cm]
 &&\leq\bP\left(N_n^{(1)}>(1+x)2\overline{\Lambda}_1(n^{-1})\right)+\bP\left(K_{n-1}^{(1)}(0,0,\infty)>(1+x)Z_1\overline{\Lambda}_1(n^{-1})\right)\nonumber
\leq\frac{C_p^{\rm spec}+C^{(1)}_p}{x^p n^{\alpha p-1}}.\vspace{-0.1cm}
\end{eqnarray}
Next, consider $2\leq k\leq m$. Recall that we denote by $\cL_k(\cS)=\{i\in[k]\colon \Sigma_i\in\cS\}$ the set of labels in a subtree $\cS\subseteq\cT$. We set $\gamma^{(k)}_k=0$ and split the spine $[[\rho,\Sigma_k[[$ at times $\gamma^{(k)}_j=\inf\{t\ge 0\colon \#\cL_k(\cT_{(\Sigma_k)}^{\alpha\rightarrow 0}(t))\le j\}$ for $k-1\ge j\ge 1$, some of which may coincide. Repeated application of Corollary \ref{cor22}, Lemma \ref{lm23} and arguing as for (\ref{fdm2}) yields 
that\vspace{-0.3cm}
$$D^{(k)}_n\le 2(D_n^{(k)}-1)\le 2N_n^{(k)}+2K_{n-k}^{(k,1)}\left(\xi^{\Sigma_k}(\gamma_1^{(k)}),\gamma_1^{(k)},\infty\right)+\sum_{j=2}^k 2K_{n-k}^{(k,j)}\left(\xi^{\Sigma_k}(\gamma_j^{(k)}),\gamma_j^{(k)},\gamma_{j-1}^{(k)}\right),\vspace{-0.3cm}$$
where $K_{n-k}^{(k,j)}\left(\xi^{\Sigma_k}(\gamma_j^{(k)}),\gamma_j^{(k)},\gamma_{j-1}^{(k)}\right)$ is as in Lemma 
\ref{pjs}, but here associated with the subordinator $\xi^{(k,j)}=\xi^{\Sigma_{k,j}}(\gamma_j^{(k)}+\cdot)-\xi^{\Sigma_{k,j}}(\gamma_j^{(k)})$ that has L\'evy measure $\Lambda^{(k,j)}=\Lambda_1$ and with random variables $V_i^{(k,j)}=\inf\{t\ge 0\colon \Sigma_{k+i}^*\not\in\cT^{\alpha\rightarrow 0}_{(\Sigma_{k,j})}(t)\}$, $i\ge 1$, 
where $\Sigma_{k,j}=\Sigma_\ell$ if $\ell=\min \cL_k(\cT^{\alpha\rightarrow 0}_{(\Sigma_k)}(\gamma_j^{(k)}))$.

\noindent Let $\displaystyle Z^{(k)}(\gamma_j^{(k)},\gamma_{j-1}^{(k)})=1+(1+A_\alpha)C_{\Lambda_1}\sum_{i=[\gamma_j^{(k)}]}^{[\gamma_{j-1}^{(k)}]}\exp(-\varrho\xi^{\Sigma_k}_i)$, noting $\displaystyle\sum_{j=1}^kZ^{(k)}(\gamma_j^{(k)},\gamma_{j-1}^{(k)})<Z_k$. Then\vspace{-0.1cm}
\begin{eqnarray*}
&&\hspace*{-0.5cm}\bP\left(D^{(k)}_n>2(1+x)(2+Z_k)\overline{\Lambda}_1(n^{-1})\right)
\leq\bP\left(N_n^{(k)}>2(1+x)\overline{\Lambda}_1(n^{-1})\right)\\[-0.1cm]
&&\hspace{3cm}+\sum_{j=1}^k\bP\left(K_{n-k}^{(k,j)}(\xi^{\Sigma_k}(\gamma_j^{(k)}),\gamma_j^{(k)},\gamma_{j-1}^{(k)})>(1+x)Z^{(k)}(\gamma_j^{(k)},\gamma_{j-1}^{(k)})\Lambda_1(n^{-1})\right)
\end{eqnarray*}
\vspace{-0.4cm}

\noindent and we can conclude again by Proposition \ref{sbp} and Lemma \ref{pjs} with constant $C_p^{\rm spec}+kC_p^{(1)}$.

Now consider $k\ge m+1$. We set $\gamma_{m+1}^{(k)}=0$ and $\gamma_0^{(k)}=\infty$. We split $[[\rho,\Sigma_k[[$ at times $\gamma_m^{(k)}=\inf\{t\ge 0\colon \Sigma_k^*\not\in\cT^{\alpha\rightarrow 0}_{(\Sigma_k)}(t)\}$ and
$\gamma^{(k)}_j=\inf\{t\ge\gamma_m^{(k)}\colon \#\cL_k(\cT_{(\Sigma_k)}^{\alpha\rightarrow 0}(t))\le j\}$ for $m-1\ge j\ge 1$.
Note that, by Procedure \ref{proc3}, $\#\cL_k(\cT_{(\Sigma_k)}^{\alpha\rightarrow 0}(\gamma_m^{(k)}))\le m$. Again \vspace{-0.2cm}
$$D^{(k)}_n\le 2(D_n^{(k)}-1)\le 2N_n^{(k)}+\sum_{j=1}^k 2K_{n-m}^{(k,j)}\left(\varepsilon_j^{(k)},\gamma_j^{(k)},\gamma_{j-1}^{(k)}\right)+2K_{n-m}^{(k,*)}(0,0,\gamma_m^{(k)}),\vspace{-0.2cm}$$
where $\varepsilon_j^{(k)}=\xi^{\Sigma_k}(\gamma_j^{(k)})$, other notation as for $k\le m$, and  
$K_{n-m}^{(k,*)}(0,0,\gamma_m^{(k)})$ is as in Lemma \ref{pjs}, here 
based on the subordinator $\xi^{\Sigma_k^*}$ with L\'evy measure $\Lambda^*$, and 
$V_i^{(k,*)}=\inf\{t\ge 0\colon \Sigma_i^*\not\in\cT^{\alpha\rightarrow 0}_{(\Sigma_k^*)}(t)\}$, the time when $\Sigma_k^*$ and $\Sigma_i^*$ are first in different subtrees. We get\vspace{-0.3cm}

\begin{eqnarray*}
&&\hspace*{-0.5cm}\bP\left(D^{(k)}_n>2(1+x)(2+Z_k)\max\{\overline{\Lambda}_1(n^{-1}),\overline{\Lambda}^*(n^{-1})\right)
\leq\bP\left(N_n^{(k)}>2(1+x)\overline{\Lambda}_1(n^{-1})\right)\\[-0.2cm]
&&\hspace{3cm}+\sum_{j=1}^k\bP\left(K_{n-m}^{(k,j)}(\xi^{\Sigma_k}(\gamma_j^{(k)}),\gamma_j^{(k)},\gamma_{j-1}^{(k)})>(1+x)Z^{(k)}(\gamma_j^{(k)},\gamma_{j-1}^{(k)})\Lambda_1(n^{-1})\right)\\
&&\hspace{3cm}+\bP\left(K_{n-m}^{(k,*)}(0,0,\gamma_m^{(k)})>(1+x)Z^{(k)}(0,\gamma_m^{(k)})\overline{\Lambda}^*(n^{-1})\right)
\end{eqnarray*}
\vspace{-0.2cm}

\noindent and conclude again by Proposition \ref{sbp} and Lemma \ref{pjs} with constant $C_p^\prime=C_p^{\rm spec}+mC_p^{(1)}+C_p^*$.

\noindent Let $H_{\cT}^\varrho$ be the height of the $\varrho$-self-similar CRT $(\cT^\varrho, \mu^\varrho)$ obtained from $(\cT,\mu)$ by $\varrho$-self-similar time-change. By \cite[Proposition 14]{H}, the height $H_{\cT}^\varrho$ has exponential moments and so does $Z_k$:
\begin{eqnarray*}\sup_{k\geq 1} Z_k&\leq&
m+(1+A_\alpha)\max\{C_{\Lambda_1},C_{\Lambda^*}\}\left(m+\sup_{k\geq 1}\int_0^\infty
\left(X_{(\Sigma_k)}^{\alpha\rightarrow 0}(t)\right)^\varrho dt\right)\\
&\leq&
m+(1+A_\alpha)\max\{C_{\Lambda_1},C_{\Lambda^*}\}\left(m+H_{\cT}^\varrho\right).
\end{eqnarray*}\nopagebreak
\vspace{-1cm}

\end{pf}

\subsection{Proof of Theorem \ref{thm3}}

The previous sections contain the new developments that we need to apply the techniques developed in \cite{HMPW} for the exchangeable case in the higher generality of Theorem \ref{thm3}. We only briefly retrace this argument here so as to identify the places where a result in the previous sections here replaces a more specific result of \cite{HMPW}.

\begin{lemm}[Lemma 10 and Corollary 11 of \cite{HMPW}]\label{Hn}
Let $H_n=\max_{1\leq k\leq n}D_n^{(k)}$ be the height of $T_n$. Then there is a
constant $C_{p,a}$ for all $a>0$, $p\geq 2/\alpha$, such that for all $x\geq 1$ and $n\ge 1$
$$\bP\left(\frac{H_n}{n^\alpha\ell(n)}>ax\right)\leq
\frac{C_{p,a}}{x^p}.$$
\end{lemm}
The proof is based on Lemma \ref{Dj} replacing \cite[Lemma 12]{HMPW}, and $\overline{\Lambda}_1(n^{-1})\sim\overline{\Lambda}^*(n^{-1})\sim n^\alpha\ell(n)$.

\begin{lemm}[Proposition 9 of \cite{HMPW}]\label{TC} Under the hypotheses of Theorem \ref{thm3}, let for $n\geq k$
$$\Delta(n,k):=\max_{1\leq i\leq n} {\rm d}_n(\{i\}, R(T_n,[k])),$$
${\rm d}_n$ being the metric associated with $T_n$. Then for each
$\eta>0$,
$$\lim_{k\rightarrow
\infty}\limsup_{n\rightarrow\infty}\bP\left(\frac{\Delta(n,k)}{n^\alpha\ell(n)}>\eta\right)=0.$$
\end{lemm}
The proof is based on Proposition \ref{fill} or (\ref{doublelim}) replacing \cite[``Clearly, $\lambda_{\rm max}^k:=\max_{j\ge 1}\lambda_j^k\rightarrow 0$ a.s.'' on page 1819]{HMPW}, Corollary \ref{cor24} replacing \cite[reference {[10]} there, Lemma 3.14]{HMPW}, and Lemma \ref{Hn} replacing \cite[Corollary 11]{HMPW}.

\begin{pfofthm3} This proof is now based on (\ref{doublelim}) replacing \cite[reference {[29]}]{HMPW}, Lemma \ref{TC} replacing \cite[Proposition 9]{HMPW}, Proposition \ref{fdm} replacing \cite[Proposition 7]{HMPW}.
\end{pfofthm3}


\section{Skewed PD model; proofs of Propositions \ref{propagpd} and \ref{skewsc}}
\label{alphagamma}


Recall that Proposition \ref{propagpd} asserts that the alpha-gamma model for $\alpha\in[0,1)$ and $\gamma\in[0,\alpha]$ is a RE Markov branching model with dislocation measures of the form identified in Corollary \ref{cor5} with $\nu_1=(1-\alpha){\rm PD}_{\alpha,-\alpha-\gamma}^*$ and $\nu_j=\gamma{\rm PD}_{\alpha,-\alpha-\gamma}^*$, $j\ge 2$.

\begin{pfofprop8} We focus on the multifurcating case $\alpha\in(0,1)$ and $\gamma\in[0,\alpha)$, the binary case being easier. We claim that the distribution of the partition $\Pi_n$ of $T_n\sim Q_n^{\alpha,\gamma}$ at $[n]$ is given by
  $$\bP(\Pi_n=\pi)=\left\{\begin{array}{ll}p^1_n(n_1,\ldots,n_k)=(1-\alpha)\frac{\Gamma(2-\alpha)}{\Gamma(n+1-\alpha)}\frac{\alpha^{k-2}\Gamma(k-1-\gamma/\alpha)}{\Gamma(1-\gamma/\alpha)}\prod_{i=1}^k\frac{\Gamma(n_i-\alpha)}{\Gamma(1-\alpha)},&\pi\in\cK^{{\bf 0}_{[2]}}\cap\cP_n,\\  p^2_n(n_1,\ldots,n_k)=\gamma\frac{\Gamma(2-\alpha)}{\Gamma(n+1-\alpha)}\frac{\alpha^{k-2}\Gamma(k-1-\gamma/\alpha)}{\Gamma(1-\gamma/\alpha)}\prod_{i=1}^k\frac{\Gamma(n_i-\alpha)}{\Gamma(1-\alpha)},&\pi\in\cK^{{\bf 1}_{[2]}}\cap\cP_n,\end{array}\right.$$
and that $(Q_n^{\alpha,\gamma},n\ge 2)$ has the labelled Markov branching property
$$\bP(\Pi_n=\pi,S_1^n=\fs_1,\ldots,S_k^n=\fs_k)=p^j_n(\#\pi_1,\ldots,\#\pi_k)\prod_{i=1}^kQ_{\pi_i}^{\alpha,\gamma}(\{\fs_i\}),\qquad\pi\in\cP^j_n,$$
where $S_i^n$ is the $i$th subtree of $T_n$ above the first branchpoint, and $Q_{\pi_i}^{\alpha,\gamma}$ is the push-forward of $Q^{\alpha,\gamma}_{\#\pi_i}$ under the natural bijection on the set of hierarchies induced by the 
increasing bijection from $[\#\pi_i]$ to $\pi_i$.

We show this by induction on $n$. Specifically, for $n=2$, this is trivial, for $n=3$ we have e.g.
$$\bP(\Pi_3=\{\{1,3\},\{2\}\})=\bP(\Pi_3=\{\{1\},\{2,3\}\})=\frac{1-\alpha}{2-\alpha},\quad
  \bP(\Pi_3=\{\{1,2\},\{3\}\})=\frac{\gamma}{2-\alpha}.$$
If the claim holds for $n$, we can apply the growth rules and the induction hypothesis to see
$$\bP(\Pi_{n+1}=\{[n],\{n+1\}\},S_1^{n+1}=\fs_1,S_2^{n+1}=\{\{n+1\}\})=\frac{\gamma}{n-\alpha}Q_n(\{\fs_1\})Q_{\{n+1\}}(\{n+1\}),$$
and for $\pi=(\pi_1,\ldots,\pi_k)\in\cP^j_n$, $j=1,2$, and hierarchies $\fs_i$ of $\pi_i$, $i\neq i^\prime$, and $\fs_{i^\prime}$ of $\pi_{i^\prime}\cup\{n+1\}$,
\begin{eqnarray*}&&\hspace{-0.5cm}\bP(\Pi_{n+1}=(\pi_1,\ldots,\pi_k,\{n+1\}),S_1^{n+1}=\fs_1,\ldots,S_k^{n+1}=\fs_k,S_{k+1}^{n+1}=\{\{n+1\}\})\\&&=\frac{(k-1)\alpha-\gamma}{n-\alpha}p_n^j(\#\pi_1,\ldots,\#\pi_k)Q_{\{n+1\}}(\{\{n+1\}\})\prod_{i=1}^kQ_{\pi_i}(\{\fs_i\}),\\
&&\hspace{-0.5cm}\bP(\Pi_{n+1}=(\pi_1,\ldots,\pi_{i^\prime}\cup\{n+1\},\ldots,\pi_k),S_1^{n+1}=\fs_1,\ldots,S_k^{n+1}=\fs_k)\\&&=\frac{n_{i^\prime}-\alpha}{n-\alpha}p_n^j(\#\pi_1,\ldots,\#\pi_k)Q_{\pi_{i^\prime}\cup\{n+1\}}(\{\fs_{i^\prime}\})\prod_{i\neq i^\prime}Q_{\pi_i}(\{\fs_i\}),
\end{eqnarray*} 
as conditionally given that the insertion of $n+1$ is in subtree $S_{i^\prime}^n$, it is just as an insertion of $\#\pi_{i^\prime}+1$ into $T_{\#\pi_{i^\prime}}$, pushed forward from $[\#\pi_{i^\prime}+1]$ to $\pi_{i^{\prime}}\cup\{n+1\}$. The result follows. 
\end{pfofprop8}

Recall that Proposition \ref{skewsc} asserts that the skewed Poisson-Dirichlet model is sampling consistent only for 
  parameters that reduce it to the exchangeable Poisson-Dirichlet model or to the alpha-gamma 
  model.

\begin{pfofprop10} By Corollary \ref{cor5}, the skewed Poisson-Dirichlet model has dislocation measure
  $$\kappa=\int_{S^\downarrow}\left(\lambda\kappa_\fs\left(\cdot\cap\cP^{{\bf 0}_{[2]}}\right)+\left(1-\lambda\right)\kappa_\fs\left(\cdot\cap\cP^{{\bf 1}_{[2]}}\right)\right){\rm PD}^*_{\alpha,\theta}(d\fs).$$
  From this, we can calculate splitting rules. Specifically, we can calculate the distribution of the ranked sequence 
  $S_n=(\#\Pi_{n,1},\ldots,\#\Pi_{n,K_n})^\downarrow$ of block sizes of $\Pi_n=(\Pi_{n,1},\ldots,\Pi_{n,K_n})$ by summing
  (\ref{split}) over partitions of equal ranked sequence of block sizes and obtain
  $$\bP(S_2=(1,1))=1,\quad \bP(S_3=(1,1,1))=\frac{\lambda(2\alpha+\theta)}{D_3},\quad\bP(S_3=(2,1))=\frac{(1+\lambda)(1-\alpha)}{D_3}$$
  \begin{eqnarray*}\bP(S_4=(1,1,1,1))=\frac{\lambda(3\alpha+\theta)(2\alpha+\theta)}{D_4},&&\bP(S_4=(2,1,1))=\frac{(1+4\lambda)(2\alpha+\theta)(1-\alpha)}{D_4}\\
 \bP(S_4=(2,2))=\frac{(1+\lambda)(1-\alpha)^2}{D_4}&&\bP(S_4=(3,1))=\frac{2(1-\alpha)(2-\alpha)}{D_4},
  \end{eqnarray*}
  where $D_3$ and $D_4$ are normalisation constants of the form $a_3\lambda+b_3$ and $a_4\lambda+b_4$. 
  Using the criterion of \cite{HMPW}, sampling consistency requires, in particular, that
  $$\bP(S_3=(1,1,1))=\bP(S_4=(1,1,1,1))+\frac{1}{2}\bP(S_4=(2,1,1))+\frac{1}{4}\bP(S_4=(3,1))\bP(S_3=(1,1,1)),$$
  which upon multiplication by $D_3D_4$ is a quadratic equation in $\lambda$. Common coefficients of all terms include
  $(1-\alpha)$ and $(\theta+2\alpha)$. For $\alpha<1$ and $\theta>-2\alpha$, the quadratic equation has the two solutions
  $\lambda=1/2$ and $\lambda=(1-\alpha)/(1-\theta-2\alpha)$ corresponding, respectively, to the Poisson-Dirichlet and 
  alpha-gamma models, so no other models can be sampling consistent. 

  The exchangeable Poisson-Dirichlet model is trivially sampling consistent. The alpha-gamma model was shown in \cite{CFW} to be sampling consistent. In the excluded case $\alpha=1$ 
  models for all $\theta$ collapse to the same deterministic model where all leaves are connected directly to a single branch 
  point \cite{MPW}. For the binary case $\theta=-2\alpha$, which we also had to exclude for our argument here, we need to consider $S_5$. This gives  
  similar quadratic equations, but also leads to the required conclusion that only the alpha model $\lambda=1-\alpha$ and the 
  beta-splitting model $\lambda=1/2$ are sampling consistent. We leave details to the reader.  
\end{pfofprop10}

\begin{appendix}
\section{Proof of Lemma \ref{pjs}}

The first part of Lemma \ref{pjs} is a straightforward consequence of \cite{GPYalpha06}, see also \cite[Lemma 8]{HMPW}. The second part generalises \cite[Lemma 12]{HMPW}. 
In the following, we indicate the most relevant changes that needed for our higher generality.

Let $N_y(t_1,t_2)$ denote the number of jumps of $\xi$ of size at least
$y$ in the time interval $[t_1,t_2]$,
$\widetilde{N}_y^{\epsilon,\tau}(t_1,t_2)$ denote the number of
jumps of $\exp(-\epsilon)(1-\exp(-\xi))$ of size at least $y$
in the same time interval. 

\textbf{Step 1. Large deviations for
$\widetilde{N}_y^{\epsilon,\tau}(0,\tau^{\prime\prime})$.} 

\begin{lemm}\label{lm13} For all $x>0$ and $0<y\le 1$, 
  $$\bP\left(\widetilde{N}_y^{\epsilon,\tau}(0,\tau^{\prime\prime})>(1+x)C_\Lambda\sum_{i=0}^{[\tau^{\prime\prime}]}\exp(-\varrho(\epsilon+\xi_i))\overline{\Lambda}(y)\right)\leq\exp(-a_x\overline{\Lambda}(y)),$$
  where $a_x:=(1+x)\ln(1+x)-x>0$.
\end{lemm}
\begin{pf} We adapt the proof of \cite[Lemma 36]{HMPW}. Let $\cF_t^{\epsilon, \tau}$ denote the $\sigma$-field generated by $(\epsilon, \tau,\tau^{\prime\prime}\wedge t)$ and $\xi$ until time $t$, and $\cF_\infty^{\epsilon, \tau}$ the one generated by $(\epsilon,\tau,\tau^{\prime\prime})$ and $\xi$, and
observe that
$$\widetilde{N}_y^{\epsilon,\tau}(0,\tau^{\prime\prime})\leq
\sum_{i=0}^{[\tau^{\prime\prime}]}\widetilde{N}_y^{\epsilon,\tau}(i,i+1)\leq
\sum_{i=0}^{[\tau^{\prime\prime}]}N_{y\exp(\epsilon+\xi_i)}(i,i+1).$$
Conditional on $\cF_i^{\epsilon, \tau}$,
$N_{y\exp(\epsilon+\xi_i)}(i,i+1)$ is a Poisson random variable
with mean $\overline{\Lambda}(y\exp(\epsilon+\xi_i))$. 
The remainder of the proof of \cite[Lemma 36]{HMPW} now applies to give
$$\bP\left(\sum_{i=0}^{[\tau^{\prime\prime}]\wedge n}N_{y\exp(\epsilon+\xi_i)}(i,i+1)\geq(1+x)C_\Lambda\sum_{i=0}^{[\tau^{\prime\prime}]\wedge n}\exp(-\varrho(\epsilon+\xi_i))\overline{\Lambda}(y)\right)
\leq\exp(-a_x\overline{\Lambda}(y)),$$
and we can let $n\rightarrow\infty$ and apply Fatou's lemma to complete the proof. 
\end{pf}

\textbf{Step 2. Large deviations for $\bE[K_n(\epsilon,
\tau,\tau^\prime)|\cF_{\tau^{\prime\prime}}^{\epsilon, \tau}]$.} 

\begin{lemm}\label{lm14} Let $B_\alpha:=\sum_{k\ge 1}\exp(-4^{-1}a_1k^{\alpha/2})$ with $a_1=2\ln 2-1$. Then for all $x\ge 1$ and all integers $n$ large enough,
  $$\bP\left(\bE\left[K_n(\epsilon,\tau,\tau^\prime)\left|\cF_{\tau^{\prime\prime}}^{\epsilon,\tau}\right.\right]>(1+x)(Y(\epsilon,\tau,\tau^\prime)-1)\overline{\Lambda}(n^{-1})\right)\le(1+B_\alpha)\exp(-4^{-1}a_1x\overline{\Lambda}(n^{-1})).$$
\end{lemm}
\begin{pf} We adapt the proof of \cite[Lemma 14]{HMPW}. 
According to
formula (4) of \cite{GPYalpha06},
$$ \bE\left[K_n(\epsilon,
\tau,\tau^\prime)\left| \cF_{\tau^{\prime\prime}}^{\epsilon,
\tau}\right.\right]=n\int_0^1(1-y)^{n-1}\widetilde{N}_y^{\epsilon,\tau}(0,\tau^{\prime\prime})
dy\leq
\widetilde{N}_{1/n}^{\epsilon,\tau}(0,\tau^{\prime\prime})+n\int_0^{1/n}\widetilde{N}_y^{\epsilon,\tau}(0,\tau^{\prime\prime})
dy.$$ 
Hence, setting
$S:=C_\Lambda\sum_{i=0}^{[\tau^{\prime\prime}]}\exp(-\varrho(\epsilon+\xi_i)),$
\begin{eqnarray}
&&\hspace*{-1cm}\bP\left(\bE\left[K_n(\epsilon, \tau,\tau^{\prime})\left |
\cF_{\tau^{\prime\prime}}^{\epsilon,
\tau}\right.\right]>(1+x)(1+A_\alpha)S\overline{\Lambda}(n^{-1})\right)\nonumber\\
&\leq
&\bP\left(\widetilde{N}_{1/n}^{\epsilon,\tau}(0,\tau^{\prime\prime})>(1+x)S\overline{\Lambda}(n^{-1})\right)+\bP\left(n\int_0^{1/n}\widetilde{N}_y^{\epsilon,\tau}(0,\tau^{\prime\prime})
dy>(1+x)A_\alpha S\overline{\Lambda}(n^{-1})\right)\nonumber.
\end{eqnarray}
The first probability in the RHS is smaller than
$\exp(-a_x\overline{\Lambda}(n^{-1}))$ by Lemma \ref{lm13}. To bound the
second probability, we use
$n\int_{1/(k+1)n}^{1/kn}\widetilde{N}_y^{\epsilon,\tau}(0,\tau^{\prime\prime})dy\leq
\widetilde{N}_{1/(n(k+1))}^{\epsilon,\tau}(0,\tau^{\prime\prime})\frac{1}{k(k+1)},$
which gives
\begin{eqnarray}&&\hspace*{-1cm}\bP\left(n\int_0^{1/n}\widetilde{N}_y^{\epsilon,\tau}(0,\tau^{\prime\prime})
dy>A_\alpha(1+x)S\overline{\Lambda}(n^{-1})\right)\nonumber\\ &\leq&
\sum_{k\ge 1}\bP\left(\widetilde{N}_{1/(n(k+1))}^{\epsilon,\tau}(0,\tau^{\prime\prime})>2(k+1)^{\sqrt{\alpha}}(1+x)S\overline{\Lambda}(n^{-1})\right),\nonumber
\end{eqnarray}
and we proceed as in \cite[Lemma 14]{HMPW} to see that this is bounded by
$\exp(-4^{-1}a_1x\overline{\Lambda}(n^{-1}))B_\alpha$ for all $x\ge 1$ 
and $n$ large enough. 
\end{pf}

\textbf{Step 3. Proof of inequality (\ref{tcp0}).} We adapt the proof of \cite[(28)]{HMPW}. To start with, fix $x\geq 1$, $n\in\bN$, and note that
\begin{eqnarray}\label{tcp3}&&\hspace*{-1cm}\bP\left(K_n(\epsilon, \tau,
\tau^\prime)>(1+x)Y(\epsilon, \tau,
\tau^\prime)\overline{\Lambda}(n^{-1})\right)\nonumber\\
&\leq& \bP\left(\bE\left[K_n(\epsilon, \tau,\tau^\prime)\left
|\cF_{\tau^{\prime\prime}}^{\epsilon,\tau}\right.\right]>(1+x)\left(Y(\epsilon, \tau,
\tau^\prime)-1\right)\overline{\Lambda}(n^{-1})\right)\\
&&+\bP\left(K_n(\epsilon, \tau,\tau^\prime)-\bE\left[K_n(\epsilon,
\tau,\tau^\prime)\left|\cF_{\tau^{\prime\prime}}^{\epsilon,\tau}\right.\right]>(1+x)\overline{\Lambda}(n^{-1})\right).\nonumber
\end{eqnarray}
Lemma \ref{lm14} gives an upper bound for the first probability provided $n$
is large enough. To get an upper bound for the second probability, we proceed as for
\cite[(28)]{HMPW} 
to find that for all $m\ge 1$, 
there exists some deterministic constant $B_m$ depending only on  $m$ such that
\begin{eqnarray}\label{tcp4}
&&\hspace*{-1cm}\bP\left(K_n(\epsilon, \tau,
\tau^\prime)-\bE\left[K_n(\epsilon, \tau,\tau^\prime)\left|\cF_{\tau^{\prime\prime}}^{\epsilon,\tau}\right.\right]>(1+x)\overline{\Lambda}(n^{-1})\left|
\cF_{\tau^{\prime\prime}}^{\epsilon,\tau}\right.\right)\nonumber\\
&\leq& 2^{m-1}B_m\frac{\bE\left[(K_n(\epsilon, \tau, \tau^{\prime}))^m\left|\cF_{\tau^{\prime\prime}}^{\epsilon,\tau}\right.\right]+\left((1+x)\overline{\Lambda}(n^{-1})\right)^m}{\left((1+x)\overline{\Lambda}(n^{-1})\right)^{2m}},\nonumber
\end{eqnarray}
We then take expectations on both sides of the resulting inequality. Theorem 6.3 of \cite{GPYalpha06} ensures that $\bE[(K_n(\epsilon,\tau,\tau^\prime))^m|\epsilon,\tau]\le\bE[(K_n(0,0,\infty)^m]\sim(\overline{\Lambda}(n^{-1}))^m$, up to a constant. Hence, we have
\begin{equation}\label{thirty}
\bP(K_n(\epsilon,\tau,\tau^\prime)-\bE[K_n(\epsilon,\tau,\tau^\prime)|\cF_{\tau^{\prime\prime}}^{\epsilon,\tau}]>(1+x)\overline{\Lambda}(n^{-1}))\leq B_{m,\Lambda}\left((1+x)\overline{\Lambda}(n^{-1})\right)^{-m},
\end{equation}
where  $B_{m,\Lambda}$ depends only on $m$ and $\Lambda$. The proof of (\ref{tcp0}) now follows the proof of \cite[(28)]{HMPW}.

%
%
This completes the proof of Lemma \ref{pjs}.\hspace*{\fill} $\square$

\end{appendix}

\subsection*{Acknowledgement}

We would like to thank Sasha Gnedin for drawing Kerov's work to our attention and for further stimulating discussions. Thanks are also due to Gerold Alsmeyer for directing us to Allan Gut's work \cite{Gut-74}. The first author would like to acknowledge support from the K C Wong Education Foundation. We would like to thank an anonymous referee for carefully reading the paper and for making suggestions that led to an improvement of the presentation.

\bibliographystyle{abbrv}
\bibliography{cfw}

\end{document}